%% file: main_OO.tex
\documentclass[ijoc,nonblindrev]{informs3_OO} 

\OneAndAHalfSpacedXII 

\usepackage{natbib}
 \bibpunct[, ]{(}{)}{,}{a}{}{,}%
 \def\bibfont{\small}%
 \def\bibsep{\smallskipamount}%
\TheoremsNumberedThrough     

\EquationsNumberedThrough    

\input{Defs}

\begin{document}






\RUNAUTHOR{\c{S}eker et al.}

\RUNTITLE{A Multiobjective Approach for Sector Duration Optimization}


\TITLE{A Multiobjective Approach for Sector Duration Optimization in Stereotactic Radiosurgery Treatment Planning}


\ARTICLEAUTHORS{%
\AUTHOR{Oylum \c{S}eker}
\AFF{Department of Mechanical and Industrial Engineering, University of Toronto, Toronto, Ontario, Canada, 
\EMAIL{oylum.seker@utoronto.ca} 
}
\AUTHOR{Mucahit Cevik}
\AFF{Department of Mechanical and Industrial Engineering, Ryerson University, Toronto, Ontario, Canada, \\
\EMAIL{mcevik@ryerson.ca}
}
\AUTHOR{Merve Bodur}
\AFF{Department of Mechanical and Industrial Engineering, University of Toronto, Toronto, Ontario, Canada, \\
\EMAIL{bodur@mie.utoronto.ca}
}
\AUTHOR{Young Lee-Bartlett}
\AFF{Elekta Instrument AB, Toronto, Ontario, Canada, \\
\EMAIL{young.lee-bartlett@elekta.com}
}
\AUTHOR{Mark Ruschin}
\AFF{Department of Radiation Oncology, Sunnybrook Health Sciences Centre, University of Toronto, Ontario, Canada,\\
\EMAIL{mark.ruschin@sunnybrook.ca}
}
}

\ABSTRACT{%
\input{abstract}
}%


\maketitle


\input{introduction_OO}

\input{literature}

\input{model}

\input{eps_cons_method}

\input{enhancements}

\input{twophase}

\input{results_OO}

\input{conclusion}

\ACKNOWLEDGMENT{%
This work is in part supported by NSERC Discovery Grant RGPIN-2019-05588.
We thank the two anonymous referees for their careful review and constructive comments that helped us improve the content and exposition of the paper. We are also grateful to Berk G\"{o}rg\"{u}l\"{u} for his valuable suggestions regarding the machine learning-related portion of this study.
}

\setlength\bibsep{-2pt}
\bibliographystyle{informs2014} 
\bibliography{references}

\newpage
\begin{APPENDICES}

\EquationsNumberedThrough  
\renewcommand{\theequation}{A.\arabic{equation}}

\input{appendix_proofs_OO}
\newpage
\input{appendix_results}

\end{APPENDICES}

\end{document}

%% file: Defs.tex
\usepackage{makecell}
\usepackage{booktabs}
\usepackage{array, tabularx}
\usepackage{comment}
\usepackage{xcolor}
\usepackage{subcaption}
\usepackage{multirow}
\usepackage{rotating}
\usepackage{pdflscape}
\usepackage{afterpage}
\usepackage{colortbl}
\usepackage{soul}

\usepackage{paralist}
\usepackage{hyperref}
\hypersetup{
    colorlinks = true,
    citecolor = {blue},
    linkbordercolor = {white},
}
\usepackage{xr-hyper}
\usepackage{cleveref}

\usepackage{algorithm}
\newcommand{\algoand}{\textbf{ and }}
\usepackage{varwidth}

\usepackage[noend]{algpseudocode}
\makeatletter
\def\algbackskip{\hskip-\ALG@thistlm}
\makeatother

\algnewcommand\Or{\textbf{or~}}

\usepackage{tikz}
\usetikzlibrary{shapes,arrows.meta}
\usetikzlibrary{calc}

\usepackage[nocomma]{optidef}

\usepackage{siunitx}

\newcolumntype{M}[1]{>{\centering\arraybackslash}m{#1}}
\DeclarePairedDelimiterX{\norm}[1]{\lVert}{\rVert}{#1}

\makeatletter
\newcommand*{\addFileDependency}[1]{
  \typeout{(#1)}
  \@addtofilelist{#1}
  \IfFileExists{#1}{}{\typeout{No file #1.}}
}
\makeatother

\newcommand{\R}{\mathbb{R}}
\newcommand{\Z}{\mathbb{Z}}
\newcommand{\Rplus}{\R_{+}}
\newcommand{\Zplus}{\Z_{+}}

\newcommand{\cred}{\color{red!85!black}}

\newcommand{\cblue}{\color{black}}

\colorlet{crevi}{purple!80!blue}
\colorlet{crevii}{red}
\newcommand{\crevi}{\color{black}}
\newcommand{\crevii}{\color{black}}

\newcommand{\BI}{\begin{itemize}}
\newcommand{\EI}{\end{itemize}}
\newcommand{\I}{\item}
\newcommand{\BE}{\begin{enumerate}}
\newcommand{\EE}{\end{enumerate}}

\newcommand{\mop}{\text{MOP}}
\newcommand{\mopeps}{\mop(\eps)}

\newcommand{\sdo}{SDO}
\newcommand{\srs}{SRS}
\newcommand{\gk}{Gamma Knife}
\newcommand{\lgk}{LGK}
\newcommand{\icon}{Icon{\texttrademark}}
\newcommand{\perfexion}{Perfexion{\texttrademark}}
\newcommand{\lgkicon}{{\lgk} \icon}

\newcommand{\structure}{\ell}
\newcommand{\setStructure}{\mathcal{L}}
\newcommand{\dose}{d}

\newcommand{\beamWeight}{g}
\newcommand{\doseMat}{G}
\newcommand{\doseOver}{o}
\newcommand{\doseUnder}{u}
\newcommand{\BOT}{b}
\newcommand{\targetDose}{D}

\newcommand{\setVoxel}{\mathcal{V}}
\newcommand{\voxel}{v}
\newcommand{\setIsocenter}{\mathcal{I}}
\newcommand{\isocenter}{\theta}
\newcommand{\setSector}{\mathcal{S}}
\newcommand{\sector}{s}
\newcommand{\setSectorSize}{\mathcal{K}}
\newcommand{\sectorSize}{k}

\newcommand{\setTumor}{\mathcal{T}}
\newcommand{\tumor}{t}
\newcommand{\setOAR}{\mathcal{C}}
\newcommand{\oar}{c}
\newcommand{\setRing}{\mathcal{R}} 
\newcommand{\ring}{r} 

\newcommand{\epsSet}{\mathcal{E}}

\newcommand{\epsinfeas}{\eps^{\text{infeas}}}

\newcommand{\numGrids}{r}

\newcommand{\objVar}{f}
\newcommand{\objVal}{z}
\newcommand{\objVarOptVal}{\hat{\objVar}}

\newcommand{\objVarFive}{h}

\newcommand{\slack}{y}
\newcommand{\objRange}{\delta}
\newcommand{\objIndex}{i}
\newcommand{\numObjs}{p}
\newcommand{\epsCoefObj}{\beta}

\newcommand{\feasSpace}{\mathcal{X}}

\newcommand{\allVars}{x}

\newcommand{\payoffModel}{P}
\newcommand{\payoffIndex}{k}

\newcommand{\objRangeLB}{lb}
\newcommand{\objRangeUB}{ub}

\newcommand{\doseUnderSol}{\hat{\doseUnder}}
\newcommand{\covSol}{\widehat{\cov}}
\newcommand{\beamWeightSol}{\hat{\beamWeight}}
\newcommand{\BOTSol}{\hat{\BOT}}
\newcommand{\objVarFiveSol}{\hat{\objVarFive}}

\newcommand{\cov}{\text{Cov}}

\newcommand{\epsi}{\eps^1}
\newcommand{\epsii}{\eps^2}
\newcommand{\nonbindingSubset}{N}

\newcommand{\nonbindingSubsetEps}{\nonbindingSubset_{\eps}}
\newcommand{\nonbindingSubsetEpsi}{\nonbindingSubset_{\epsi}}

\newcommand{\objVarOptVali}{\objVarOptVal^1}
\newcommand{\objVarOptValii}{\objVarOptVal^2}
\newcommand{\bindingSubsetIndex}{k}
\newcommand{\optObjSet}{\hat{F}}
\newcommand{\optSolSet}{\hat{X}}
\newcommand{\optSol}{\hat{\allVars}}
\newcommand{\optSoli}{\hat{\allVars}^1}

\newcommand{\optSlack}{\hat{\slack}}

\newcommand{\eps}{\varepsilon}

\newcommand{\sel}{Select}

\newcommand{\epsSetSelI}{\epsSet^{\text{I},\text{\sc \sel}}}
\newcommand{\epsSetSelII}{\epsSet^{\text{II},\text{\sc \sel}}}
\newcommand{\epsSetI}{\epsSet^{\text{I}}}
\newcommand{\epsSetII}{\epsSet^{\text{II}}}
\newcommand{\mlprob}{\rho}
\newcommand{\trainFnc}{\texttt{Train}}

\newcommand{\predFnc}{\text{Pred}}

\newcommand{\epsSetFeasTrain}{{\epsSet}^\text{TrFeas}}
\newcommand{\epsSetPerfTrain}{{\epsSet}^\text{TrPred}}
\newcommand{\notdom}{\mathcal{ND}}

\newcommand{\coverage}{\text{Cov}}
\newcommand{\PI}{\text{PCI}}
\newcommand{\BOTVal}{\text{BOT}}
\newcommand{\feas}{\text{Feas}}

\newcommand{\coverageThres}{\coverage^{\text{min}}}
\newcommand{\PIThres}{\PI^{\text{min}}}
\newcommand{\BOTThres}{\BOTVal^{\text{max}}}

\newcommand{\epsSetPromising}{\epsSet^{\text{I}, \text{\sc \sel More}}}

\newcommand{\TV}{\text{TV}}
\newcommand{\TVPIV}{\text{TV}_\text{PIV}}
\newcommand{\PIV}{\text{PIV}}

\newcommand{\cevikModel}{SOLP}

\newcommand{\twophase}{\scshape 2phas$\eps$}
\newcommand{\twophasereg}{\twophase-r}
\newcommand{\twophaseml}{\twophase-ml}

\newcommand{\tuplesetfeastr}{\tupleset^\text{TrFeas}}
\newcommand{\tuplesetperftr}{\tupleset^\text{TrPerf}}

%% file: abstract.tex
Sector duration optimization (\sdo) is a problem arising in treatment planning for stereotactic radiosurgery on Gamma Knife.
Given a set of isocenter locations, {\sdo} aims to select collimator size configurations and irradiation times thereof such that target tissues receive prescribed doses in a reasonable amount of treatment time, while healthy tissues nearby are spared.
We present a multiobjective linear programming model for {\sdo} 
to generate a diverse collection of solutions so that 
{\crevii clinicians can select the most appropriate treatment.}
We develop a generic two-phase solution strategy based on the $\eps$-constraint method for solving multiobjective optimization models, {\twophase}, which aims to systematically increase the number of high-quality solutions obtained, instead of conducting a traditional uniform search. 
To improve solution quality further and to accelerate the procedure, we incorporate some general and problem-specific enhancements. 
Moreover, we propose an alternative version of {\twophase}, which makes use of machine learning tools to
{\crevii reduce the computational effort.}
In our computational study on eight previously treated real test cases, a significant portion of {\twophase} solutions outperformed clinical results and those from a single-objective model from the literature. 
In addition to significant benefits of the algorithmic enhancements, 
our experiments illustrate the usefulness of machine learning strategies 
{\crevii to reduce the overall run times nearly by half while maintaining or besting the clinical practice.}

%% file: introduction_OO.tex
\section{Introduction}

Radiotherapy is one of the most commonly used cancer treatment methods.
The technique uses high-energy radiation to kill or damage cancer cells,
and is employed in approximately half of the treatments, sometimes in conjunction with chemotherapy and surgery~\citep{breedveld2019multi}.
There are different forms of radiotherapy, one of which is {\it stereotactic radiosurgery} ({\srs}), {\crevi a technique specifically designed to deliver high doses of radiation from multiple positions to precisely defined targets in the brain.} 
{\gk} was the first {\srs} delivery system, but the state of the art is Leksell Gamma Knife\textsuperscript{\textregistered} ({\lgk}) Icon{\texttrademark} by Elekta (Stockholm, Sweden)~\citep{zeverino2017}.

{\lgkicon} (and also its predecessor {\perfexion}) has eight independently controlled {\it sectors} and contains 192 cobalt-60 (a radioactive isotope) sources, 24 per sector.
Each sector can be robotically moved to one of the three circular {\it collimator sizes} (with 4, 8, and 16 mm diameters) or be blocked. 
All the sources in a sector deliver radiation beams through the chosen size, or are completely blocked. 
A patient lies still on a treatment couch and radiation beams are emitted from the sectors to the target, where their focal point is termed an {\it isocenter}.
The couch remains stationary while radiation is being delivered to a specific isocenter and is repositioned automatically for other isocenter locations defined in the treatment plan. 
{\crevii All beams are blocked as the patient is repositioned.}
Each isocenter is irradiated for a certain amount of time through one or more collimator configurations, and the treatment session
{\crevii is complete once the planned radiation amounts have been delivered to each isocenter.}

{\crevi 
Different types of radiotherapy techniques require similar decisions to be made to determine a treatment plan. 
{\crevii Common decisions include the irradiation times as well as the positions and aperture shapes and/or sizes through which radiation beams target cancerous regions.}
For {\gk} systems, a set of aperture sizes at every sector is fixed on the device, and the planner needs to decide on the irradiation durations for each of~them.
}

{\srs} treatment planning 
has two main steps, those being 
(i) determination of the number and locations of isocenters, and (ii) finding the irradiation times of each isocenter. 
Existing studies mostly {\crevii handle} 
these two problems separately.
For the former, heuristic methods are commonly used; for instance, \cite{ghobadi2012automated} combine the well-known grassfire algorithm \citep{wagner2000geometrically} with sphere packing, and \cite{doudareva2015skeletonization} introduce a method based on skeletonization techniques.
We focus on the latter problem in this study.

Given a set of isocenter locations, {\crevii the} {\it sector duration optimization} ({\sdo}) problem {\crevii aims to find} the irradiation durations of every isocenter from each sector through each collimator size, such that a quality solution is achieved in terms of {\crevi clinical goals}. 
{\crevi The problem 
has high dimensionality and incorporates a large number of clinical restrictions.} 
Consequently, it is not generally possible to generate an optimal {\crevi or clinically acceptable} treatment plan simply by inspection, 
thus mathematical models and optimization tools are needed.

{\crevi As in other radiotherapy problems, multiple criteria are considered in {\sdo}. Target sites should receive the prescribed doses, while surrounding tissues should be spared.}
The duration of the treatment is another important consideration, because longer treatment times may cause patient discomfort and increases the possibility of patient movement leading to dose uncertainties \citep{srsbook2019ch7}.
These and other criteria are {\crevi commonly considered in radiotherapy treatment planning. They are} conflicting in nature, and it is generally not possible to simultaneously attain individual best values for each. 
{\crevi Moreover,} clinicians are often interested in evaluating trade-offs among competing criteria 
and choosing a plan most suitable to the specifics of a case.
{\crevii Motivated by this,}
we {\crevii propose} a multiobjective approach for {\sdo} arising in {\srs} treatment planning {\crevi on Gamma Knife systems}, particularly with {\lgkicon}.
Our aim is to provide a collection of solutions addressing the trade-offs between conflicting clinical goals, and thereby offer flexibility to decision makers {\crevii as they select a plan for an individual patient.}

{\crevi We present a multiobjective linear programming model for {\sdo}, 
with the the objectives being to minimize the underdosing of tumors, the overdosing of tumors and healthy tissues, and the treatment time.}
We introduce a two-phase solution approach that employs {\crevii a widely-used technique}
to solve multiobjective optimization models, namely, the $\eps$-constraint method,
{\crevi which consists of optimizing one objective while bounding the others.}
We perform the first phase to explore solutions from all parts of the feasible space, and the second phase to restrict the search to achieve a denser representation of clinically desirable solutions. 
{\crevi Generated solutions yield Pareto-optimal points in terms of the model objectives; however, we continue the evaluative process by evaluating a treatment plan with some widely-used clinical measures, namely, tumor coverage and Paddick conformity index. Note that these clinical metrics are not directly optimized by the model, and 
in that regard, our objectives, other than the treatment time, serve as surrogates to derive solutions that yield clinically preferable options.}
We conduct computational experiments using anonymized data from eight previously treated cases, and we compare our results with the clinical plans and the treatment plans from the single-objective model by \cite{cevik2019simultaneous}.
{\crevi
We note that even though we demonstrate our solution scheme on the SDO problem with {\lgkicon}, it is a generic framework that can be applied to other radiosurgery modalities, and even other multiobjective problems.
}

Our main contributions {\crevii can be summarized as follows}:
\begin{itemize}
    \item We present a multiobjective linear programming model for {\sdo} 
    and generate a diverse collection of solutions, offering flexibility to decision makers to select one that best suits patient-specific needs.
    We also derive inequalities that improve our model and that aid the solution process.
    %
    \item We provide a self-contained alternative proof showing that a properly constructed $\eps$-constraint model yields an efficient solution for the multiobjective model.
    \item We prevent generation of repeating solutions in the $\eps$-constraint method by extending an existing strategy to linear programs, applicable regardless of problem context. 
    \item We propose a generic and efficacious two-phase $\eps$-constraint method, {\twophase}, that systematically focuses the search on more favorable regions. We additionally develop an alternative version of {\twophase} that makes use of machine learning tools to predict performance, which thereby reduces the computational burden.
\end{itemize}

{\crevi
{\crevii The rest of the paper is organized as follows.} We survey the related literature in Section~\ref{sec:literature} and present our multiobjective model for SDO in Section~\ref{sec:model}.
We then review the $\varepsilon$-constraint method and propose problem-specific improvements in Section~\ref{sec:epscons}, and we provide algorithmic enhancements in Section~\ref{sec:enhancements}.
We introduce {\twophase} in Section~\ref{sec:twophase}, report our experimental results in Section~\ref{sec:expresults}, and conclude our work in Section~\ref{sec:conc}.
We provide all the proofs in Appendix~\ref{app:proofs}. 
}

%% file: literature.tex
\section{Literature Review}
\label{sec:literature}

{\crevii We first review the related {\srs} treatment planning literature, and then provide an overview of the related multiobjective optimization studies.}

\subsection{{\crevii {\srs} Treatment Planning}}

Several {\srs} technologies have been developed since the introduction of {\gk} devices.
Linear accelerator-based systems comprise a commonly-used class of such devices, which are widely used for conventional radiotherapy, starting in the early 1980s \citep{chen2005stereotactic}.
{\crevi There are various studies focusing on {\srs} treatment planning with such delivery systems, e.g., \citep{yu1997multiobjective,yu2000multi,lee2000optimization}.
Linear accelerator-based systems are also used in intensity modulated radiation therapy (IMRT), and there exist various studies presenting multiobjective approaches for IMRT treatment planning, e.g., \citep{romeijn2004unifying,breedveld2012icycle,cabrera2018matheuristic,cabrera2018pareto}.}
Though {\sdo} has similarities with IMRT treatment planning, the differences in delivery devices precludes a direct correspondence between them.
We refer the interested reader to the survey studies by \cite{ehrgott2010mathematical} and \cite{breedveld2019multi} on optimization for IMRT treatment planning and multiobjective optimization for radiotherapy in general, respectively.
{\crevi 
We restrict our review to studies {\crevii focusing on} 
{\sdo} on {\gk} systems, sometimes in conjunction with isocenter placement problem.
}

{\crevi Relatively early studies on {\gk} treatment planning for older devices required manual interventions to change the collimator size and couch position, leading to limitations on the number of isocenters and collimator sizes.
These studies use mixed integer and nonlinear programming models, as well as heuristics to determine the isocenter locations and irradiation durations, e.g., \citep{ferris2000optimization,shepard2000inverse,shepard2003clinical,zhang2001optimization,zhang2003plug}.
A viable strategy is to determine a treatment plan in multiple steps to distribute the computational burden, starting with a coarse grid approximation by considering only a subset of unit volumes and refining it in subsequent steps, e.g., \citep{ferris2002optimization,ferris2003radiosurgery,cheek2005relationship}.
\cite{ferris2002optimization} solve a series of nonlinear and mixed integer programming models, where solution quality gradually improves {\crevii while maintaining clinical efficacy}.
A similar approach is followed by \cite{ferris2003radiosurgery}.
An earlier work by \cite{shepard1999optimizing} presents radiotherapy treatment planning models that impose bounds on one or more criteria while optimizing a weighted sum of them. 
The $\eps$-constraint model is built similarly, and the $\eps$-constraint method we employ systematically varies the bounds on the objectives and solves the corresponding models to obtain different solutions. 
Another relevant study in that regard is by \cite{schlaefer2008stepwise}, which presents a step-wise multiobjective approach that iteratively reaches a desired trade-off between clinical goals by optimizing one objective and bounding the others at each step. 
The method by \cite{lee2006automated} derives an initial plan by finding a closely matching treatment volume from a database of past treatments, and it improves the treatment with a nonlinear objective function that seeks to improve the conformity.
}
In the studies {\crevii handling} 
{\crevi treatment planning for more recent {\gk} systems,} 
{\crevi different mathematical optimization models are presented with multiple objectives that are mostly combined into one as a weighted sum.
{\crevii Heuristic and exact methods then solve the resulting single objective model.}} 
Proposed models for {\sdo} include a second-order cone optimization model \citep{oskoorouchi2011interior}; 
convex models with quadratic and/or linear penalty terms \citep{ghobadi2012automated,ghobadi2013automated,ghaffari2017incorporation}; 
linear, quadratic, and unconstrained convex moment-based penalty models \citep{cevik2018modeling}; 
and more standard linear models \citep{cevik2019simultaneous,sjolund2019linear,tian2020preliminary}.
%

{\crevi We lastly review} multiobjective approaches for {\gk} treatment planning that explore the trade-offs among the {\crevi clinical} goals 
by offering an assortment of solutions. 
One such study by \cite{giller2011feasibility} conducts Pareto analysis via a heuristic, namely, a genetic algorithm, by considering the trade-off between the percentage of tumor covered and the volume of normal tissue exposed. 
Another work by \cite{svensson2014multiobjective} adapts an existing Pareto-surface approximation algorithm
for {\sdo}, 
{\crevi and it develops a graphical user interface to aid the selection of a treatment.}
{\crevi The work by \cite{ripsman2015interactive} also develops a graphical user interface for {\srs}. 
The underlying algorithm solves an approximate linearized model that minimizes a weighted sum of doses delivered to ring tissues. 
}

\subsection{\crevii Multiobjective Optimization}
\label{sec:lit_rev_eps}

{\crevii 

A goal in multiobjective optimization is to generate the set of all Pareto-optimal points (Pareto set), although if this is not desirable or reasonable, then a representative subset is sought.
The $\eps$-constraint method is a generic technique for handling multiobjective optimization problems, introduced by \cite{haimes1971bicriterion} and later discussed in detail by \cite{cohon1978multiobjective} and \cite{chankong1983multiobjective}.
It has been widely used for different model classes, such as linear and integer programs, to find the Pareto set (e.g., see \cite{ozlen2009multi,kirlik2014new,tamby2021enumeration}), or a representation of it (e.g., see \cite{eusebio2014finding,jindal2017multi}) which we also aim to find in this study.

The idea of the $\eps$-constraint method is to optimize one objective while bounding the remaining ones, and to explore different Pareto-optimal points by varying the bounds on the objectives. 
Selection of the objective bounds determines what portion of the Pareto set the model is confined to, and it is thus an important design component of the $\eps$-constraint method.
The work by \cite{cohon1978multiobjective} suggests to place an evenly spaced grid over the objective space and solve the model at each point on the grid which serves as objective bounds. This technique has been widely used in the literature, e.g, by \cite{mavrotas2009effective}.
We also employ this technique within each phase of our two-phase $\eps$-constraint method, {\twophase}, which is a generic framework particularly suitable for problems where at least some of the model objectives serve as surrogates to the actual performance measures, those being widely-used clinical metrics in our case.
In the regular version of {\twophase}, we carry out an evenly-distributed search in the first phase, and based on the values of the clinical metrics that the obtained solutions yield, we determine the promising regions to conduct another evenly-distributed search therein in the second phase.
In the machine learning-guided version of {\twophase}, we make predictions to further guide the search in both phases in order to reduce the overall computational effort.

The $\eps$-constraint method we employ in this study is most closely related to the works of \cite{cohon1978multiobjective,mavrotas2009effective,mavrotas2013improved,zhang2014simple} and \cite{nikas2020robust}, which we discuss in more detail in Sections~\ref{sec:epscons} and \ref{sec:enhancements}.
We mention a few other $\eps$-constraint method related studies here, which are, to the best of our knowledge, the remaining relevant works for this survey.
The study by \cite{laumanns2006efficient} develops an adaptive scheme, where the objective bounds are generated dynamically based on the previously obtained solutions, and infeasible and already-discovered regions are updated during the search.
Another study by \cite{sahebjamnia2015integrated} presents a modified $\eps$-constraint model that incorporates decision maker's preferences as objective weights, 
and it uses that model to generate solutions while avoiding infeasibilities and repeating solutions.
The work by \cite{jindal2017multi} presents an interactive $\eps$-constraint method that produces sets of Pareto-optimal points until the decision maker is satisfied with one of them. If the decision maker does not pick any of the presented points, the most preferred segments of objective ranges are determined and the next search is conducted there.
The similarity of our study with the work of \cite{jindal2017multi} is that we also identify the preferred objective ranges and ultimately confine our search to those regions, but using certain clinical metrics within a two-phase framework rather than interactively consulting to a decision maker. 


We note that some attributes have been introduced to assess the quality of discrete representations of the Pareto set.
For instance, \cite{sayin2000measuring} defines coverage, uniformity, and cardinality as the three attributes of quality, while \cite{shao2016discrete} propose a method that provides guarantees to the extent coverage and uniformity are achieved. 
Previous quality measures from the literature do not align well with our goal of producing solutions from specific parts of the feasible space,
particularly those that yield clinically desirable performance.
Therefore, we chose a general-purpose solution method for multiobjective optimization models, namely, the $\eps$-constraint method, and we tailor it for our problem within a two-phase framework, explained in Section~\ref{sec:twophase}.




}

%% file: model.tex
\section{A Multiobjective Model for {\sdo}}
\label{sec:model}

We assume a predetermined set of isocenter positions, with the goal of SDO being to determine the irradiation times through every sector for each collimator size and isocenter location.
{\crevii The received radiation should conform to prescribed doses for the targets, and exposure of nearby healthy tissues should be as low as possible.}
These nearby tissues involve those surrounding the targets, which we label as {\it rings}, and sometimes specific organs as well, called {\it organs-at-risk} (OARs) or {\it critical structures}, such as the brain stem or cochlea, which are especially radiosensitive.
Also, the {\it beam-on time} (BOT), i.e., the total radiation delivery time of the treatment, should be as low as possible.
We note that BOT approximately indicates the overall treatment time, because automatic collimator adjustment takes less than three seconds \citep{novotny2012leksell} and does not exceed a few minutes in total.
All structures 
are discretized into a collection of small cubes, called {\it voxels} (volumetric pixels), for planning purposes, which are typically 1 mm or less in length.

\subsection{Notation}
\label{subsec:model_notation}

Let $\setStructure$ be the set of all structures, comprised of the sets of tumors, rings, and OARs denoted by $\setTumor, \setRing,$ and $\setOAR$, respectively.  
Each structure $\structure \in \setStructure$ consists of a set $\setVoxel_{\structure}$ of voxels. 
The sets of isocenters, sectors, and collimator sizes are respectively $\setIsocenter, \setSector,$ and $\setSectorSize$.
The $\doseMat_{\structure \voxel \isocenter \sector \sectorSize}$ parameters represent the dose rate (dose per unit time) absorbed by voxel $\voxel \in \setVoxel_{\structure}$ of structure $\structure \in \setStructure$ via isocenter $\isocenter \in \setIsocenter$ through sector $\sector \in \setSector$ at collimator size $\sectorSize \in \setSectorSize$.
The $\targetDose_{\tumor}$ and $\targetDose_{\structure}^{\max}$ parameters respectively denote the prescribed dose for tumor $\tumor \in \setTumor$ and the desired maximum dose for structure $\structure \in \setStructure$, measured in units of gray (Gy).

The variable $\beamWeight_{\isocenter \sector \sectorSize}$ is the duration isocenter $\isocenter \in \setIsocenter$ is irradiated through sector $\sector \in \setSector$ at collimator size $\sectorSize \in \setSectorSize$, and $\dose_{\structure \voxel}$ denotes the dose received by voxel $\voxel \in \setVoxel_{\structure}$ of structure $\structure \in \setStructure$.
Since it may not always be possible to deliver the required dose to each voxel of the targets without violating prescribed limitations, we define auxiliary under- and overdose variables as well.
The variable $\doseUnder_{\tumor \voxel}$ denotes the underdose amount to voxel $\voxel \in \setVoxel_{\tumor}$ of tumor $\tumor \in \setTumor$, and $\doseOver_{\structure \voxel}$ denotes the overdose to voxel $\voxel \in \setVoxel_{\structure}$ of structure $\structure \in \setStructure$.  
Finally, $\BOT_{\isocenter}$ is the BOT value of isocenter $\isocenter \in \setIsocenter$.


\subsection{Five-Objective Linear Program}
\label{subsec:model_5objmodel}

Our model is an extension of the single-objective model by \cite{cevik2018modeling} to our multiobjective setting:
%
\input{Models/model_5obj}


We denote the five objectives of our MOLP with variables $\objVarFive_1, \ldots, \objVarFive_5$ in \eqref{obj5}, and we set their associated expressions in constraints \eqref{cons:objfix1}--\eqref{cons:objfix5}.
The objective function denoted with 
$\objVarFive_1$ is
the total overdose amount to rings, 
$\objVarFive_2$ is the total dose amount to OARs and rings, 
$\objVarFive_3$ is the total overdose amount to tumors, 
$\objVarFive_4$ is the total underdose amount to tumors,
and $\objVarFive_5$ is the BOT of the treatment plan. 

Constraints \eqref{cons:dose} determine the dose delivered to each voxel in each structure. 
Constraints \eqref{cons:underdose} through \eqref{cons:dose_oar} aim to keep the dose delivered to voxels of each structure within desired limits. 
In particular, constraints \eqref{cons:underdose} enforce the received dose and underdose amount to add up to at least the required dose for each tumor voxel.
Similarly, constraints \eqref{cons:overdose_tumor} and \eqref{cons:overdose_ring} respectively link the under- and overdose amounts with their received doses and dose limits for tumors and rings.
Since OARs are subject to strict dose restrictions, we do not define overdose variables for OARs, and thus constraints \eqref{cons:dose_oar} simply ensure that OARs do not receive more than their maximum allowed doses. 
For the sake of finding a feasible treatment plan, the model is allowed to go beyond the maximum dose limits for rings, in which case the excessive amount of radiation is measured by overdose variables.

At each isocenter location, radiation is simultaneously delivered from all sectors that are not blocked, through one or more collimator sizes. 
The amount of time spent at each sector $\sector$ is the total time used over all collimator sizes $\sectorSize \in \setSectorSize$ in it, and the time spent at each isocenter location $\isocenter$ is the maximum time spent over all sectors $\sector \in \setSector$.
Constraints \eqref{cons:bot} force the BOT value at each isocenter to be at least as much as the total time radiation is delivered to it. 
The sum of $\BOT_{\isocenter}$ variables then corresponds to total BOT.
Constraints \eqref{cons:qnonneg}--\eqref{cons:ononneg} are nonnegativity restrictions on the decision variables. 

The model in \citep{cevik2019simultaneous} has constraints as \eqref{cons:dose}--\eqref{cons:qnonneg}, and its objective is a weighted sum of the objectives that we treat separately with the difference being that we do not allow overdosing of OARs to protect the critical structures from excessive damage. 
{\crevi We note that these models are not specific to {\lgkicon} and that they apply to other versions of the device. 
} 

{\crevi Identifying the weights to combine the objectives as a weighted sum is a nontrivial task. For instance, \cite{cevik2019simultaneous} and \cite{babier2018inverse} determine the objective weights via a simulated annealing algorithm and inverse optimization, respectively. Our solution strategy is based on treating the different objectives separately rather than combined, as we explain in the next section.
}

%% file: Models/model_5obj.tex
\begin{subequations}
\label{model:5obj}
\begin{alignat}{3}
& \min_{\objVarFive, \dose, \beamWeight, \BOT, \doseUnder, \doseOver} \qquad && \left(\objVarFive_{1}, \objVarFive_{2}, \objVarFive_{3}, \objVarFive_{4}, \objVarFive_{5} \right) \label{obj5} \\[0.1cm]
& \text{\qquad s.t.} && \objVarFive_1 \ = \ \sum_{\ring \in \setRing} \sum_{\voxel \in \setVoxel_{\ring}} \doseOver_{\ring \voxel} \label{cons:objfix1} \\[0.05cm]
&&& \objVarFive_2 \ = \ \sum_{\oar \in \setOAR} \sum_{\voxel \in \setVoxel_{\oar}} \dose_{\oar \voxel} \ + \ \sum_{\ring \in \setRing} \sum_{\voxel \in \setVoxel_{\ring}} \dose_{\ring \voxel} \label{cons:objfix2} \\[0.05cm]
&&& \objVarFive_3 \ = \ \sum_{\tumor \in \setTumor} \sum_{\voxel \in \setVoxel_{\tumor}} \doseOver_{\tumor \voxel} \label{cons:objfix3} \\[0.05cm]
&&& \objVarFive_4 \ = \ \sum_{\tumor \in \setTumor} \sum_{\voxel \in \setVoxel_{\tumor}} \doseUnder_{\tumor \voxel} \label{cons:objfix4} \\[0.05cm]
&&& \objVarFive_5 \ = \ \sum_{\isocenter \in \setIsocenter} \BOT_{\isocenter}  \label{cons:objfix5} \\[0.05cm]
&&& \dose_{\structure \voxel} \ = \ \sum_{\isocenter \in \setIsocenter} \sum_{\sector \in \setSector} \sum_{\sectorSize \in \setSectorSize} \doseMat_{\structure \voxel \isocenter \sector \sectorSize} \ \beamWeight_{\isocenter \sector \sectorSize} && \qquad \structure \in \setStructure, \ \voxel \in \setVoxel_{\structure} \label{cons:dose} \\[0.05cm]
&&&\dose_{\tumor \voxel} \ + \ \doseUnder_{\tumor \voxel} \ \geq \ \targetDose_{\tumor}  && \qquad \tumor \in \setTumor, \ \voxel \in \setVoxel_{\tumor} \label{cons:underdose} \\[0.05cm]
%
%
&&& \dose_{\tumor \voxel} \ - \ \doseOver_{\tumor \voxel} \ \leq \ \targetDose_{\tumor}^{\max}  && \qquad \tumor \in \setTumor, \ \voxel \in \setVoxel_{\tumor}  \label{cons:overdose_tumor} \\[0.05cm]
&&& \dose_{\ring \voxel} \ - \ \doseOver_{\ring \voxel} \ \leq \ \targetDose_{\ring}^{\max}  && \qquad \ring \in \setRing, \ \voxel \in \setVoxel_{\ring}  \label{cons:overdose_ring} \\[0.05cm]
&&& \dose_{\oar \voxel} \ \leq \ \targetDose_{\oar}^{\max}  && \qquad \oar \in \setOAR, \ \voxel \in \setVoxel_{\oar}  \label{cons:dose_oar} \\[0.05cm]
&&& \BOT_{\isocenter} \ \geq \ \sum_{\sectorSize \in \setSectorSize} \beamWeight_{\isocenter \sector \sectorSize}  && \qquad \isocenter \in \setIsocenter, \ \sector \in \setSector \label{cons:bot} \\[0.05cm]
&&& \BOT_{\isocenter} \ \geq \ 0 && \qquad  \isocenter \in \setIsocenter \label{cons:qnonneg} \\
&&& \beamWeight_{\isocenter \sector \sectorSize} \ \geq \ 0 && \qquad  \isocenter \in \setIsocenter, \ \sector \in \setSector, \ \sectorSize \in \setSectorSize  \label{cons:wnonneg} \\
&&& \doseUnder_{\tumor \voxel} \ \geq \ 0 && \qquad  \tumor \in \setTumor, \ \voxel \in \setVoxel_{\tumor}  \label{cons:unonneg} \\
&&& \doseOver_{\structure \voxel} \ \geq \ 0 && \qquad  \structure \in \setStructure, \ \voxel \in \setVoxel_{\structure}  \label{cons:ononneg} 
%
%
\end{alignat}
\end{subequations}


%% file: eps_cons_method.tex

\section{The $\boldsymbol{\eps}$-Constraint Method}
\label{sec:epscons}

{\crevii The traditional approach to solving multiobjective optimization models is \emph{scalarization}, i.e., formulating a single-objective model related to it \citep{ehrgott2005multicriteria}.
The $\eps$-constraint method is one of the most popular scalarization techniques for general multiobjective optimization models.}
We review the $\eps$-constraint method
and present some problem-specific improvements.

\subsection{Preliminaries}
\label{subsec:epscons_prelim}

A generic multiobjective program has the form 
\begin{mini}|s|[2] 
    {\substack{\allVars \in \feasSpace}}{ 
    \left(\objVar_{1}(\allVars), \ldots, \objVar_{\numObjs}(\allVars) \right) \label{model:mop} }{}{\hspace{-3cm}\text{\mop}: \quad}
\end{mini}

\noindent where $\numObjs \in \Zplus$ is the number of objectives, $\allVars$ is the vector of decision variables, $\feasSpace \subseteq \R^n$ 
is the set of feasible solutions, and $\objVar \colon \feasSpace \rightarrow \R^{\numObjs}$ maps each solution $\allVars \in \feasSpace$ into a $\numObjs$-dimensional {objective vector} $\objVar(\allVars) \coloneqq \left(\objVar_{1}(\allVars), \ldots, \objVar_{\numObjs}(\allVars) \right)$, with $\objVar_{\objIndex} \colon \feasSpace \rightarrow \R $ for $\objIndex \in  \{ 1, \ldots, \numObjs \}$.

{
\newcommand{\allVarsSpec}{\Tilde{\allVars}}
\newcommand{\objIndexSpec}{\objIndex^{\prime}}
\newcommand{\objValSpec}{\Tilde{\objVal}}

A solution $\allVars \in \feasSpace$ is {\it weakly efficient} if there is no $\allVarsSpec \in \feasSpace$ satisfying $\objVar_{\objIndex}(\allVarsSpec) < \objVar_{\objIndex}(\allVars) $ for all $\objIndex \in \{1, \ldots, \numObjs\}$.
A solution $\allVars \in \feasSpace$ is {\it efficient} if there exists no other solution $\allVarsSpec \in \feasSpace$ such that $\objVar_{\objIndex}(\allVarsSpec) \leq \objVar_{\objIndex}(\allVars) $ for all $\objIndex \in \{1, \ldots, \numObjs\}$, with at least one inequality being strict, 
i.e., $\objVar_{\objIndexSpec}(\allVarsSpec) < \objVar_{\objIndexSpec}(\allVars) $ for some index $\objIndexSpec$.
Similarly, 
$\objVar(\allVars)$ {\it dominates} $\objVar(\allVarsSpec)$ if $\objVar_{\objIndex}(\allVars) \leq \objVar_{\objIndex}(\allVarsSpec)$ for all $\objIndex \in \{1, \ldots, \numObjs\}$, with at least one index $\objIndexSpec$ such that $\objVar_{\objIndexSpec}(\allVars) < \objVar_{\objIndexSpec}(\allVarsSpec)$ holds.
An objective vector $\objVar(\allVars)$ is {\it nondominated}, or {\it Pareto-optimal}, if there is no other solution $\allVarsSpec \in \feasSpace$ such that $\objVar(\allVarsSpec)$ dominates $\objVar(\allVars)$.
{\crevii The set of all Pareto-optimal points is the {\it Pareto set}.}
}

\subsection{The Formulation}
\label{subsec:epscons_method}

The \emph{$\eps$-constraint method} first appeared in \citep{haimes1971bicriterion}, but it is discussed in detail in \citep{cohon1978multiobjective} and \citep{chankong1983multiobjective}. 
In this method, an \emph{$\eps$-constraint formulation} is constructed for {\mop} by selecting one of the objectives as being primary, then minimizing it while constraining the remaining $\numObjs-1$ objectives.
An $\eps$-constraint formulation for {\mop} is 
\begin{mini}|s|[2] 
    {\substack{\allVars \in \feasSpace}}{\crevi \objVar_{k}(\allVars) \label{model:eps_cons}}{}{}
    \addConstraint{\objVar_{\objIndex}(\allVars)}{\ \leq \ \eps_{\objIndex} \qquad}{\crevi \objIndex \in \{1, \ldots, \numObjs\}, \ \objIndex \neq k}
\end{mini}

\noindent {\crevii where $\eps_{\objIndex}$ bounds the $\objIndex^{th}$ objective.}
{\crevi 
A solution $\hat{\allVars} \in \feasSpace$ is efficient if and only if it is an optimal solution of the $\eps$-constraint model in \eqref{model:eps_cons} for every $k = 1, \ldots, \numObjs$, with $ \eps_{\objIndex} = \objVar_\objIndex(\hat{\allVars})$ for $\objIndex=1,\ldots,\numObjs, \ \objIndex \neq k$. 
This means that with appropriate choices of $\eps$, all efficient solutions are identifiable in theory, independent of the structure of $\feasSpace$ 
\citep{chankong1983multiobjective,ehrgott2008improved}.}

{\crevi 
The $\eps$-constraint model relates to a \emph{weighted sum} model $P(w)$ that combines objectives as a weighted sum, 
i.e., $P(w) \colon \min \left\{ \sum_{\objIndex = 1}^{\numObjs} w_\objIndex \objVar_\objIndex(\allVars) \colon {\allVars \in \feasSpace} \right\}$, where $w_\objIndex \geq 0$ for all $\objIndex \in \{1,\ldots,\numObjs\}$ {\crevii and $\sum\limits_{\objIndex=1}^\numObjs w_\objIndex = 1$}.
If $\hat{\allVars}$ is an optimal solution to \eqref{model:eps_cons} with $\eps_\objIndex = \objVar_\objIndex(\hat{\allVars}), \ {\objIndex \neq k}$ and if $\feasSpace$ and $\objVar_\objIndex(\allVars)$'s are convex, then $\hat{\allVars}$ is optimal for $P(w)$ for some nonnegative $w$. 
Conversely, if $\hat{\allVars}$ is an optimal solution for $P(w)$, then either (i) if $w_k > 0$, $\hat{\allVars}$ is also optimal for \eqref{model:eps_cons} with $\eps_\objIndex = \objVar_\objIndex(\hat{\allVars}), \ {\objIndex \neq k}$, or (ii) if $\hat{\allVars}$ is the unique minimizer of $P(w)$, then $\hat{\allVars}$ is optimal for \eqref{model:eps_cons} with $\eps_\objIndex = \objVar_\objIndex(\hat{\allVars}), \ {\objIndex \neq k}$ for all $k$ \citep{chankong1983multiobjective}. 
}

The $\eps$-constraint method solves \eqref{model:eps_cons} (or a modified version of it, which we present in the sequel) for different \emph{$\eps$-vectors} to generate a diverse set of efficient solutions. 
The collection of $\eps$-vectors is built by determining a range of values for each individual $\eps_{\objIndex}$ parameter and systematically combining them.

\subsubsection*{Guaranteed efficiency of solutions.}
An optimal solution to the $\eps$-constraint formulation in \eqref{model:eps_cons} is not guaranteed to be an efficient solution for {\mop};
it may be weakly efficient  \citep{cohon1978multiobjective, chankong1983multiobjective}. 
To ensure that optimal solutions of the $\eps$-constraint model are efficient, 
the following modified formulation is proposed by \cite{ehrgott2008improved} and \cite{mavrotas2009effective}. {\crevi Letting $\objVar_1$ denote the primary objective for simplicity,}
\begin{mini}|s|[2] 
  {\substack{\allVars \in \feasSpace, \ \slack \in \Rplus^{\numObjs-1}}}{\objVar_{1}(\allVars) \ - \  \sum_{\objIndex=2}^{\numObjs} \epsCoefObj_{\objIndex} \slack_{\objIndex} }
  {\label{model:eps_cons_eff}}{\mopeps: \quad}
    \addConstraint{\objVar_{\objIndex}(\allVars) \ + \ \slack_{\objIndex}}{\ = \ \eps_{\objIndex}, \qquad}{\objIndex \in \{2, \ldots \numObjs\} }
\end{mini}

\noindent where $\epsCoefObj_{\objIndex}$'s are strictly positive weight parameters for the slack variables $\slack_{\objIndex}$'s. 
{\crevi 
We place no particular hierarchy among the objectives in our problem and scale the slack terms in the objective function with the widths of the ranges that the $\objVar_{\objIndex}(\allVars)$'s can achieve, as suggested by \cite{mavrotas2009effective}. Namely, we set $\epsCoefObj_{\objIndex} = \epsCoefObj / \objRange_{\objIndex}$, where $\objVar_{\objIndex}(\allVars)$ spans a width of $\objRange_{\objIndex}$, and $\epsCoefObj > 0$, which is chosen to be sufficiently small.}

\subsubsection*{Alternative proof of correctness.}
The optimal solutions of \eqref{model:eps_cons_eff} being efficient solutions for {\mop} follows from Lemma 3.1 and Theorem 3.1 in the work of \cite{ehrgott2008improved}, proven by making use of other findings and definitions presented in that paper. 
The study of \cite{mavrotas2009effective} also provides a proposition to show this result, but the proof 
does not fully show the generality of the argument. 
We provide a simple and self-contained alternative proof to the following in the Online Supplement. 

\begin{proposition}
\label{prop:efficiency}
\input{Propositions/prop_eps_cons_correctness}
\end{proposition}

{\crevii We note that not only an optimal solution of the $\eps$-constraint model in \eqref{model:eps_cons_eff} is an efficient solution for {\mop}, but also any efficient solution for {\mop} is in theory attainable through \eqref{model:eps_cons_eff} with appropriate choices of $\eps$. A weighted sum model, on the other hand, is not able to attain some efficient solutions (\emph{unsupported} solutions) for nonconvex problems \citep{ehrgott2006discussion}.}

\subsection{Procedure to Set $\boldsymbol{\eps}$-vectors}
\label{subsec:epscons_procedure}

\cite{cohon1978multiobjective} proposed a method to build a set of $\eps$-vectors that was later implemented by \cite{mavrotas2009effective}. 
The scheme has two main steps:
\begin{enumerate}[(i)]
  \item Construct a payoff table to set the objective ranges, and
  \item Divide each range into equidistant points and use their combinations as $\eps$-vectors.
\end{enumerate}

\subsubsection{Setting the objective ranges.}

Let $\objRangeLB_{\objIndex}$ and $\objRangeUB_{\objIndex}$ respectively denote the lower and upper limits for $\objVar_{\objIndex}(\allVars)$, and $\objRangeLB_{\objIndex}^e$ and $\objRangeUB_{\objIndex}^e$ denote the minimum and maximum values $\objVar_{\objIndex}(\allVars)$ attains over the set of efficient points of {\mop}.
Ideally,
$\objRangeLB_{\objIndex} = \objRangeLB_{\objIndex}^e$ and $\objRangeUB_{\objIndex} = \objRangeUB_{\objIndex}^e$. 
The $\objRangeLB_{\objIndex}^e$ values can be obtained by minimizing $\objVar_{\objIndex}(\allVars)$ over the {\mop} feasible region, but it is nontrivial to determine the $\objRangeUB_{\objIndex}^e$ values \citep{weistroffer1985careful, isermann1988computational}.  
A commonly used alternative is to use surrogates to construct a {\it payoff table}.
{\crevi To this end, we employ {\it lexicographic optimization}, which has been tested by \cite{isermann1988computational} and later used by \cite{mavrotas2009effective}. 
This method takes an ordering of $\numObjs$ objectives and sequentially optimizes each with the optimal values of the previous objectives being fixed to obtain a nondominated point in the end. 
The nondominated point populates one row in the payoff table.
This process is carried out $\numObjs$ times with a different objective at the start of the ordering, as shown in Algorithm~\ref{algo:objranges}.
Then, $\objRangeLB_{\objIndex}$ and $\objRangeUB_{\objIndex}$ values are specified as the lowest and highest values of the associated objectives, respectively.
{\crevii This way, we have $\objRangeLB_{\objIndex} = \objRangeLB_{\objIndex}^e$, but the $\objRangeUB_{\objIndex}$ values may fall below or exceed the $\objRangeUB_{\objIndex}^e$ values \citep{cohon1978multiobjective, weistroffer1985careful, isermann1988computational}. We note that we use nondominated points to populate the payoff table in order to obtain possibly smaller $\objRangeUB_{\objIndex}$ values, because our preliminary experiments showed that higher values do not lead to clinically desirable solutions.}
}

\input{PseudoCodes/algo_set_obj_ranges_payoff}


\subsubsection{Constructing the $\boldsymbol{\eps}$-vectors.}

We divide the objective ranges $[\objRangeLB_{\objIndex}, \objRangeUB_{\objIndex}]$, for ${\objIndex \in \{2, \ldots, \numObjs\}}$, into equal-width intervals to create a set $\epsSet_{\objIndex}$ of equidistant $\eps_{\objIndex}$ values, a technique suggested by \cite{cohon1978multiobjective} and adopted by \cite{mavrotas2009effective}. 
We essentially place a grid over the space of objectives in the constraint set as described by \cite{cohon1978multiobjective}.
%
{\crevii We set the number $\numGrids_{\objIndex}$ of equidistant points for each objective $\objIndex \in \{2, \ldots, \numObjs\}$, and enumerate the set $\epsSet$ of all $\eps$-vectors beforehand as 
{
\newcommand{\increment}{\Delta}
\begin{gather}
\epsSet \coloneqq \epsSet_{2} \times \ldots \times \epsSet_{\numObjs}  \label{eq:grid_pts1} 
\shortintertext{where}
\epsSet_{\objIndex} \coloneqq \left\{  \objRangeLB_{\objIndex} + j \cdot \frac{\objRangeUB_{\objIndex}-\objRangeLB_{\objIndex}}{\numGrids_{\objIndex} - 1} \right\}_{j = 0,1,\hdots,\numGrids_{\objIndex} - 1}, \ \ \forall \objIndex \in \{2, \ldots, \numObjs\}.  \label{eq:grid_pts2}
\end{gather}%
}%
}%
\subsection{SDO-specific Improvements}
\label{subsec:epscons_improvements}

The procedure above sets the objective ranges regardless of the problem context.
Some portions of these ranges, however, may yield solutions not useful in practice. 
One may make use of the problem characteristics or some desired levels in the outputs to reduce the objective ranges.

While constructing the payoff table, recall that each objective is minimized once before fixing the value of any objective. 
When we start with the BOT objective, the model simply chooses not to deliver any radiation.
As a result, the upper limit of the range for total underdose is always the sum of the prescribed doses, and the lower limit of the range for BOT is zero, which is unrealistic.

There are certain clinical performance criteria used to assess the quality of a treatment plan (discussed in more detail in Section~\ref{sec:twophase}). 
One of the most important measures is tumor coverage, which is the proportion of tumor voxels receiving at least the prescribed doses. 
Depending on the ranges from which $\eps_{\objIndex}$'s are chosen, a considerable number of $\eps$-vectors may yield coverage values that are not clinically acceptable. 
We can, however, use some minimum clinically desirable coverage value to confine the objective ranges and to derive additional constraints to tighten the model. We apply this idea for tumor underdosing and BOT objectives.

\subsubsection{Upper bound for tumor underdose.}
\label{subsec:ub_underdose}

We first show that there is a direct relation between the coverage achieved by a solution and the total underdose amounts for tumors.

\begin{proposition}[Underdose upper bound]
\label{prop:underdose_UB}
\input{Propositions/prop_underdose_UB}
\end{proposition}

Proposition~\ref{prop:underdose_UB} implies that if we set a minimum desired coverage value $\coverageThres$ and use it instead of $\cov_{\tumor}$, then the inequalities in \eqref{eq:underdose_extra_cons} are necessary for the MOLP to yield solutions with at least $\coverageThres$ coverage for each tumor, albeit they are not sufficient generally. 
Therefore, as an improvement to our model, we incorporate \eqref{eq:underdose_extra_cons} as an additional set of constraints, with $\coverageThres$ in place of $\cov_{\tumor}$.

We further use these inequalities to tighten the upper limit $\objRangeUB_4$ of the range of the underdose objective $\objVarFive_4$.
Summing each side over $\tumor \in \setTumor$, we obtain
\begin{align}
\objVarFive_4 \ = \ \sum\limits_{\tumor \in \setTumor} \sum\limits_{\voxel \in \setVoxel_{\tumor}} \doseUnder_{\tumor \voxel} \ \leq \  \sum\limits_{\tumor \in \setTumor} \targetDose_{\tumor} \ |\setVoxel_{\tumor}| \ (1-\cov_{\tumor}). \label{eq:underdose_UB}
\end{align}
\noindent  
If we replace $\cov_{\tumor}$'s in \eqref{eq:underdose_UB} with a minimum desired coverage value $\coverageThres$, the new upper limit for the range of the underdose objective $\objVarFive_4$ is smaller than the one we obtain in our payoff table, as the minimum desired coverage would never be zero.

\subsubsection{Lower bound for BOT.}
\label{subsec:lb_bot}

Using the constraints of our MOLP and Proposition \ref{prop:underdose_UB}, we improve the lower limit for the range of the BOT objective. 

\begin{proposition}[BOT range tightening]
\label{prop:BOT_LB}
\input{Propositions/prop_BOT_LB}
\end{proposition}

For a minimum desired coverage value $\coverageThres$ to replace $\cov_{\tumor}$ in \eqref{eq:BOT_LB}, this new bound replaces the zero lower limit for the range of the BOT objective obtained from the payoff table.
{\cblue When using \eqref{eq:underdose_extra_cons}, \eqref{eq:underdose_UB}, and \eqref{eq:BOT_LB}, we take $\coverageThres $ as 0.98.}

%% file: Propositions/prop_eps_cons_correctness.tex
An optimal solution to  \eqref{model:eps_cons_eff} with $\epsCoefObj_{\objIndex} > 0$ for all $\objIndex \in \{2, \ldots, \numObjs\}$ is an efficient solution of {\mop} in \eqref{model:mop}.

%% file: PseudoCodes/algo_set_obj_ranges_payoff.tex
\begin{algorithm}
\small
\vspace*{0.15cm}
    \caption{Setting objective ranges}\label{algo:objranges}  
    {\bf Input:} Model $\payoffModel_{\payoffIndex} \crevii \colon \min\limits_{\allVars \in \feasSpace} \objVar_\payoffIndex(\allVars)$ for $\payoffIndex \in \{1, \ldots, \numObjs\}$ \\
    {\bf Output:} Lower and upper limit pairs $(\objRangeLB_{\objIndex}, \objRangeUB_{\objIndex})$ of $\objVar_{\objIndex} $'s $ \forall \objIndex \in \{1, \ldots, \numObjs\}$
    \begin{algorithmic}[1]
        \Procedure{SetObjRanges}{}
        \For {{\crevii$\payoffIndex \coloneqq 1, \ldots, \numObjs$} }
            \For {{\crevii$\objIndex \coloneqq \payoffIndex, \ldots, \numObjs, 1, \ldots, \payoffIndex-1$} }
            \State Set the objective of $\payoffModel_{\payoffIndex}$ as $\objVar_{\objIndex}(\allVars)$ 
            \State Solve $\payoffModel_{\payoffIndex}$ to get the optimal objective value $\objVal_{\objIndex}^{\payoffIndex}$ and add $ \objVar_{\objIndex}(\allVars) \leq \objVal_{\objIndex}^{\payoffIndex}$ to $\payoffModel_{\payoffIndex}$ as a constraint
            \EndFor
        \EndFor
        \ForAll {$\objIndex \in \{1, \ldots, \numObjs\}$}
            \State $\objRangeLB_{\objIndex} \gets \min_{\payoffIndex}\{\objVal_{\objIndex}^{\payoffIndex}\} $, $\objRangeUB_{\objIndex} \gets \max_{\payoffIndex}\{\objVal_{\objIndex}^{\payoffIndex}\} $
        \EndFor
        \EndProcedure
    \end{algorithmic}
\end{algorithm}%
%

%% file: Propositions/prop_underdose_UB.tex
Any feasible solution to the MOLP in \eqref{model:5obj} yielding $\emph{\cov}_{\tumor}$'s as the tumor coverage values satisfies the following inequalities:
\begin{align}
\sum\limits_{\voxel \in \setVoxel_{\tumor}} \doseUnder_{\tumor \voxel} \ \leq \  \targetDose_{\tumor} \ |\setVoxel_{\tumor}| \ (1-\emph{\cov}_{\tumor}) \qquad  \tumor \in \setTumor \label{eq:underdose_extra_cons}
\end{align}
%
%

%% file: Propositions/prop_BOT_LB.tex
Any feasible solution to the MOLP in \eqref{model:5obj} yielding $\emph{\cov}_{\tumor}$'s as the coverage values satisfies the following inequality:
\begin{align}
\max\limits_{\tumor \in \setTumor} \left\{ \frac{\targetDose_{\tumor} \ \emph{\cov}_{\tumor}}{\max\limits_{\voxel, \isocenter, \sector, \sectorSize} \{ \doseMat_{\tumor \voxel \isocenter \sector \sectorSize} \} \ |\setSector| \ } \right\} \ \leq \ \objVarFive_5. \label{eq:BOT_LB}
\end{align}
%
%

%% file: enhancements.tex
\section{Algorithmic Enhancements}
\label{sec:enhancements}
  
We next discuss two enhancements that accelerate the $\eps$-constraint method, those being  
(i) early detection of infeasibilities, and (ii) early detection of repeating Pareto-optimal points.

\subsection{Early Detection of Infeasibilities}

Even though each individual objective range $(\objRangeLB_{\objIndex}, \objRangeUB_{\objIndex})$ contains attainable values for $\objVar_{\objIndex}$, coexistence of $\eps_{\objIndex}$'s drawn from these intervals may lead to infeasible $\mopeps$ models.
{\crevii This fact risks computational resources to identify infeasibilities, and we posit the omission of $\eps$-vectors that lead to infeasibility.}

Early detection of infeasibilities has been used in the augmented $\eps$-constraint method (AUGMECON) of \cite{mavrotas2009effective} and was later improved in SAUGMECON by \cite{zhang2014simple}.
SAUGMECON discards the $\eps$-vectors that imply infeasibility, and we do the same.
Once we encounter a vector $\epsinfeas$ such that $\mop(\epsinfeas)$ is infeasible, we search through the unexplored $\eps$-vectors to find those that are more restrictive than $\epsinfeas$ and remove them from consideration.
That is, we extract 
%
\begin{align}
    \epsSet^{\prime} \coloneqq \left\{\eps \in \epsSet \colon \exists \objIndex \in \{2, \ldots, \numObjs\} \text{~with~} \eps_{\objIndex} > \epsinfeas_{\objIndex} \right\}  \subseteq \epsSet. \label{eq:detectinfeas}
\end{align}

\subsection{Early Detection of Repeating Pareto-optimal Points}
We also detect $\eps$-vectors that yield an already discovered Pareto-optimal point. 
%
%
{\crevi 
This scheme, which is due to \cite{mavrotas2013improved}, is based on generating $\eps$-vectors via nested loops and bypassing some iterates in the innermost loop. 
SAUGMECON extends this method by forwarding some of the iterates of the outer loops, and it also uses the observation that if an optimal solution remains feasible in a tighter solution space, it remains optimal too.
{\crevii \cite{nikas2020robust} also improves upon the work of \cite{mavrotas2013improved} in a similar way.}
We exploit some properties of linear programs and propose a generalized method to identify repeating solutions.
}

Let $\optSolSet(\eps)$ be the set of optimal solution vectors $\optSol$'s for $\mopeps$, 
$\nonbindingSubsetEps(\optSol)$ be the index set of nonbinding $\eps_{\objIndex}$-constraints at $\optSol$, 
and $\optObjSet(\eps)$ be the set of
objective vectors, $\objVar(\optSol)$'s, for all $\optSol \in \optSolSet(\eps)$.

\begin{proposition}[Identifying repeating solutions]
\label{prop:equiv1} 
\input{Propositions/prop_repeating_sols_1}
\end{proposition}

We note that there may be some $\objVarOptValii \in \optObjSet(\epsii)$ with $\objVarOptValii \neq \objVarOptVali $, but we ignore these because they are likely similar to a discovered point.
We scan the unexplored vectors and discard those satisfying the conditions in Proposition~\ref{prop:equiv1}. 
Namely, given a nondominated point $\objVarOptVal \coloneqq \objVar(\optSol)$ of MOLP, we extract 
%
\begin{align}
    \epsSet^{\prime\prime} \coloneqq \left\{\eps \in \epsSet \colon \ \exists \objIndex \in \nonbindingSubsetEps(\optSol) \text{~with~} \eps_{\objIndex} < \objVarOptVal_{\objIndex}  \ \text{~or~} \ \exists \bindingSubsetIndex \in  \{2, \ldots, \numObjs\} \setminus \nonbindingSubsetEps(\optSol) \text{~with~} \eps_{\bindingSubsetIndex} \neq \objVarOptVal_{\bindingSubsetIndex} \  \right\}  \subseteq \epsSet. \label{eq:detectrepeating}
\end{align}

%% file: Propositions/prop_repeating_sols_1.tex
Suppose that $\mopeps$ is a linear program.
Let $\optSoli \in \optSolSet(\epsi)$, $\nonbindingSubsetEpsi(\optSoli) \neq \varnothing$, and $\objVarOptVali \coloneqq \objVar(\optSoli)$.  
For a given $\epsii$ with $\epsi \neq \eps^2$, if $\objVarOptVali_{\objIndex} \leq \epsii_{\objIndex} \ \ \forall \objIndex \in \nonbindingSubsetEpsi(\optSoli)$ and $\epsi_\bindingSubsetIndex = \epsii_\bindingSubsetIndex, \ \forall \bindingSubsetIndex \in \{2,\ldots,\numObjs\} \setminus \nonbindingSubsetEpsi(\optSoli)$, 
then $\objVarOptVali \in \optObjSet(\epsii)$. 
%

%% file: twophase.tex
\section{Two-phase $\boldsymbol{\eps}$-constraint Method}
\label{sec:twophase}

Determining an attainable range of values for each objective, whether using the procedure in Section~\ref{subsec:epscons_procedure} or another, does not ensure that $\eps_{\objIndex}$ values within their identified ranges yield a feasible solution.
Moreover, it is not generally obvious which subintervals of the objective ranges give desirable solutions with regard to representing the set of Pareto-optimal points.
{\crevii This fact is further challenged in our model because solutions are assessed by metrics other than our objectives, which are clinical proxies.}

A straightforward way to reduce the likelihood of missing desirable solutions is to finely divide the objective space and to increase the density of the $\eps$-vectors.
{\crevii However, this scheme burdens the computation effort with several objectives.}
Moreover, we would rather devote most of our computational effort on obtaining desirable solutions, which is unlikely if we simply increase $\eps$-vector granularity without focusing on specific regions.

We introduce a two-phase $\eps$-constraint method, called {\twophase} that helps focus our search of the $\eps$-vectors. 
We explore in Phase I the overall objective space, generate solutions, and identify the objective ranges on which desirable solutions can be obtained.
Phase II confines the search to those regions to diversify such solutions. 
We propose two versions of {\twophase}: the regular version, {\twophasereg} (Section~\ref{subsec:twophaseRegular}), and the machine learning-guided version, {\twophaseml} (Section~\ref{subsec:twophaseML}).

A flow chart for {\twophase} is in Figure~\ref{fig:flowchart}.
The sequence of boxes with solid borders and solid arcs show the flow for {\twophasereg}, while filled boxes with dashed borders display the digressions related to {\twophaseml}. 
The top and bottom parts of the figure are dedicated to Phase I and Phase II, respectively, as indicated on the left panel (Phase I boxes in teal blue, Phase II boxes in pink). 
We adopt such distinctions (achieved by colors and fillings) in the plots of Section~\ref{sec:expresults}.

\input{Figures/flowchart}

We now generally present {\twophase} and subroutines tailored for SDO.

\subsubsection*{Performance criteria.}
The five objectives we consider can have a wide range of nonnegative values.
Even though BOT corresponds to an actual amount of time (in minutes), and is thus easy to interpret, the remaining objectives are only proxies for clinical benchmarks.
To this end,
we consider two well-accepted clinical performance criteria, those being {\it coverage} (Cov) and {\it Paddick conformity index} (PCI) \citep{paddick2000simple}.
Cov is the proportion of target volume receiving at least the prescription doses, 
and PCI additionally accounts for healthy tissue exposure. 
PCI and Cov relate as
\begin{align}
    \PI \ = \ \frac{(\TVPIV)^2}{\TV \times \PIV} \ = \ \coverage \times \frac{\TVPIV}{\PIV}, 
\end{align}
where $\TV$ is the target 
volume, $\TVPIV$ is the target volume that receives at least the prescribed dose, and $\PIV$ is the total volume that receives at least the prescribed dose.
Both Cov and PCI range between zero and one, with higher values indicating better performance. 
{\crevi 
We heuristically infer ranges for the model objectives that yield desirable values for these clinical measures, 
} 
but we note that these clinical measures and the MOLP objectives each provide valuable information not implied by the other. 
{\crevi We also note that {\twophase} is a generic strategy; one that can incorporate any set of such criteria to guide the search process.} 
For our purposes, we use Cov, PCI, and BOT as our performance criteria henceforth.

\subsection{Regular Version}
\label{subsec:twophaseRegular}
We start Phase I of {\twophasereg} by specifying the objective ranges following the procedure in Section~\ref{subsec:epscons_procedure}, in conjunction with the problem-specific improvements presented in Section~\ref{subsec:epscons_improvements}.
We divide the ranges to generate $\eps$-vectors and solve the corresponding $\mopeps$ models.
In Phase II we take the Pareto-optimal points from Phase I and filter those deemed desirable, an assessment that can depend on the {\mop} objectives, or other measures such as Cov and PCI.
We let the minimum and maximum values of the objectives over the filtered sets be the bounds of their ranges in Phase II, and we divide each range into equidistant points (of not necessarily the same granularity as in Phase I).
We then enumerate the $\eps$-vectors and solve the corresponding $\eps$-constraint models.
The union of the outputs from the two phases comprise the final set of Pareto-optimal points (see Figure~\ref{fig:flowchart}).

Let $\coverage(\allVars), \PI(\allVars)$, and $\BOTVal(\allVars)$ respectively denote the Cov, PCI, and BOT values of a feasible solution $\allVars$.
We use set
threshold values $\coverageThres, \PIThres$, and $\BOTThres$ and determine objective ranges for Phase II as stated in Algorithm~\ref{algo:objranges_phase2}. 

\input{PseudoCodes/algo_set_obj_ranges}

We set the threshold values of $\coverageThres$ and $\PIThres$ to 0.98 and 0.75, respectively.
Since BOT values vary considerably, we do not set a fixed value for $\BOTThres$.
We first filter the performance criteria tuples from Phase I using $\coverageThres$ and $\PIThres$, and then set $\BOTThres$ to be the first quartile of BOT values in the remaining subset.

\subsection{Machine Learning-guided Version}
\label{subsec:twophaseML}
{\crevii We adapt {\twophase} by incorporating machine learning (ML) methods to reduce the computational effort.}
Phase I is divided into two rounds.
We first randomly select $\eps$-vectors enumerated at the beginning and solve the corresponding $\mopeps$ models.
We then train a prediction model for Phase I to predict the values of our performance criteria 
over the unexplored $\eps$-vectors, and we solve $\mopeps$ for promising $\eps$-vectors in the second round of Phase I. 
We use all the data obtained in Phase I and generate the $\eps$-vectors as in {\twophasereg}. 
We continue by training another prediction model in Phase II to forecast quality $\eps$-vectors and solve $\mopeps$ for them.

\subsubsection{Data collection.}
\label{subsubsec:twophase_datacollect}
An $\eps$-vector deemed unworthy to explore is one that either makes $\mopeps$ infeasible or produces a solution that fails to satisfy desirable criteria. 
Hence, our prediction process is twofold: (i) feasibility prediction, and (ii) performance criteria prediction. 
To this end, we collect two data sets.
We first record if an $\eps$-vector makes $\mopeps$ infeasible, and 
we then store the performance criteria for the $\eps$-vectors that produce a solution.
We employ the filters provided in \eqref{eq:detectinfeas} and \eqref{eq:detectrepeating} to detect infeasibility and repeating solutions. 
The time required by the filtering process is negligible (see Section~\ref{sec:expresults}), therefore we search the entire set of $\eps$-vectors instead of confining ourselves to a subset. 
We associate each $\eps$-vector identified to yield repeating results with a copy of the performance criteria obtained previously. 
As such, we expand our data sets considerably at no additional cost, and we use them to train prediction models for both phases. 
Note that data collection is not needed in Phase II as we do not predict anything afterwards.

\subsubsection{Prediction models.}
\label{subsubsec:twophase_predictmodels}
We now detail our prediction models (see the book by \cite{statlearningbook2013} for more details).
We consider an $\eps$-vector as input (features), and either the feasibility status or each one of the performance criterion values as output (response).
Feasibility status is {\it categorical} (feasible or infeasible), whereas the performance criteria are numerical.
Both for classification and regression, we utilize {\it decision tree-based} models, namely, {\it random forests}.

A decision tree is a graphical structure built through a series of rules that splits data.
Starting from a root node, observations are partitioned into two or more groups based on typically one of the features until a stopping condition is satisfied. 
The resulting tree partitions the feature space into a number of regions. 
The response of a given feature vector in a region is predicted to be the most commonly encountered or average response in that region, depending on if the tree is a classification or regression tree, respectively.
A random forest is a collection of decision trees, each being built on a sample drawn from the training data and a sample of features. 
For classification, given an observation of feature values, the predicted class is the most commonly encountered predicted class over all trees.
Similarly for regression, the predicted response is the average of the predicted responses from the individual trees.

We note that other ML models could be used in {\twophaseml}, but in our preliminary experiments, random forests showed 
{\crevi the best overall performance among several others. The models we tested were random forests, decision trees, nearest neighbors, ridge, support vector machines, and logistic \& linear regression (see Section~\ref{subsubsec:expresults_2phase_ml} for details). }

\subsubsection{Finding promising $\boldsymbol{\eps}$-vectors.}
\label{subsubsec:twophase_findpromising}

We train random forest classification and regression models for each performance criterion separately after data collection. 
Prediction steps explore feasibility of all candidate $\eps$-vectors.
Feasibility then prompts prediction of the performance criteria after which $\eps$-vectors failing to satisfy our desirability conditions are removed. 
We set $\BOTThres$ to a redundantly large value to explore all candidate solutions with desirable Cov and PCI values in Phase I.
Phase II is the same as in {\twophasereg}.

\input{PseudoCodes/algo_predict_select}

\subsubsection{Domination analysis.}
It is possible to solve $\mopeps$ for every $\eps$-vector estimated to produce desirable performance in each phase, but we use an additional filter if the objectives are not deemed the ultimate performance criteria.
Efficient solutions do not necessarily yield nondominated points with regard to our clinical assessments, and we filter $\eps$-vectors whose estimated performance criteria are not dominated by those of another $\eps$-vector.
We solve $\mopeps$ only for the remaining $\eps$-vectors.
We define the $\notdom$ operator in \eqref{eq:notdom} to represent this filtering,
%
\begin{align}
    \epsSet^\text{ND} \coloneqq 
    \left\{\eps^{\tupleindex} \right\}_{\tupleindex \in \tupleset^\text{ND}} = 
    \notdom \big( \left\{\eps^{\tupleindex}, \ \coverage^{\tupleindex}, \ \PI^{\tupleindex}, \ \BOTVal^{\tupleindex} \right\}_{\tupleindex \in \tupleset} \big). \label{eq:notdom}
\end{align}

We note that the performance criteria values after the domination filtering are not necessarily nondominated. 
Conversely, some performance criteria values might be dominated due to prediction errors and be eliminated.
We also note that this domination analysis is not meant for cases in which the {\mop} objectives are to be predicted, because $\mopeps$ models are guaranteed to yield nondominated objective values.

%% file: Figures/flowchart.tex
\afterpage{
\clearpage
\begin{landscape}
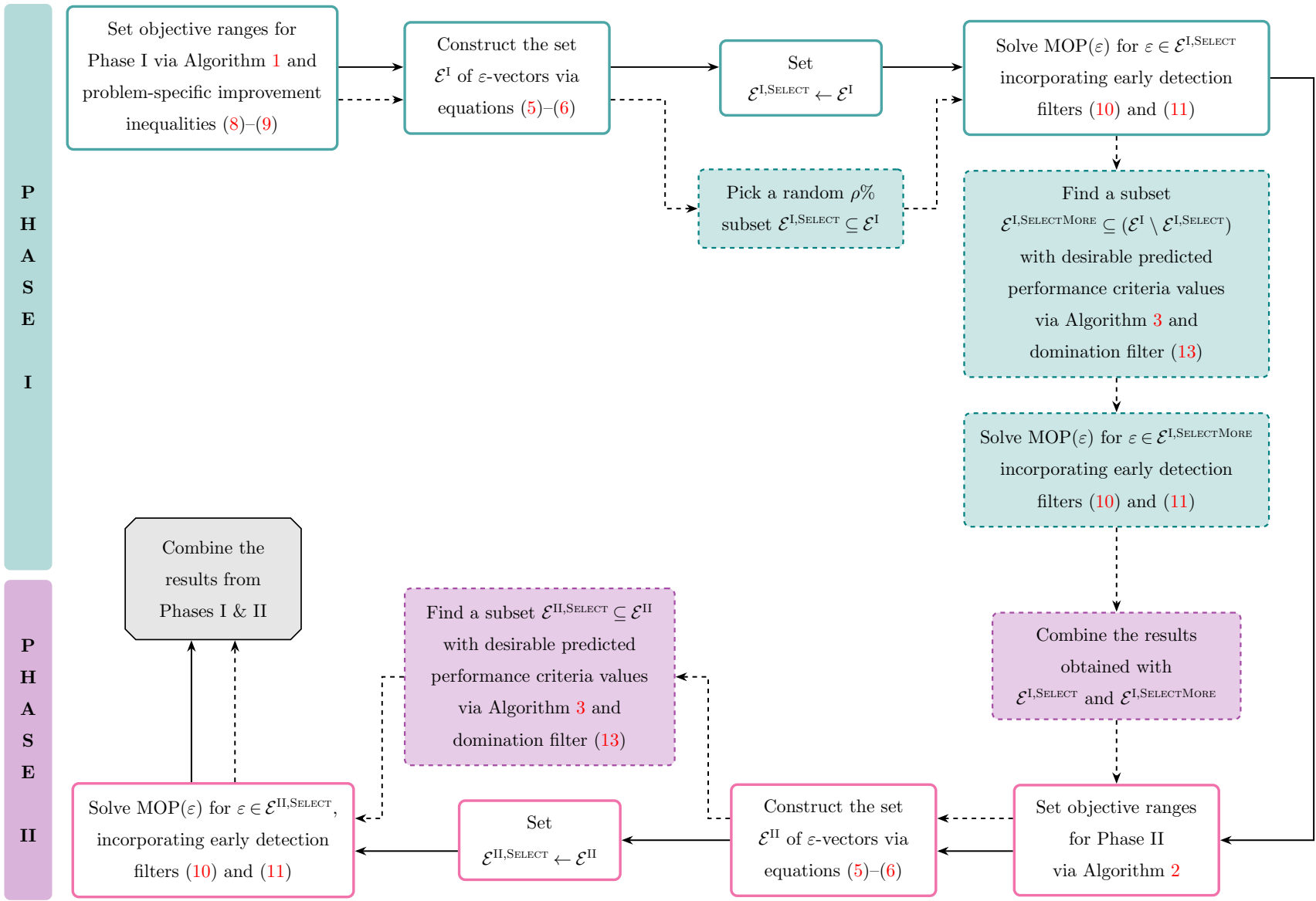
\begin{figure}
\scalebox{0.75}{

\begin{tikzpicture}[node distance = 5cm and 3cm, auto]

\tikzstyle{infonode1} = [rectangle, rounded corners, minimum height=13cm, text width=0.4cm, text centered, inner sep=3.5mm, fill=teal!30]
\tikzstyle{infonode2} = [rectangle, rounded corners, minimum height=7.35cm, text width=0.4cm, text centered, inner sep=3.5mm, fill=violet!30]
\tikzstyle{regnode1} = [rectangle, draw=teal!70, line width=0.6mm, rounded corners, minimum height=1.5cm, text width=4cm, text centered, inner sep=3.5mm]
\tikzstyle{mlnode1} = [rectangle, draw=teal!100, line width=0.4mm, dashed, fill=teal!20, rounded corners, minimum height=1.5cm, text width=3.2cm, text centered, inner sep=3.5mm]
\tikzstyle{regnode2} = [rectangle, draw=magenta!70, line width=0.6mm, rounded corners, minimum height=1.5cm, text width=4cm, text centered, inner sep=3.5mm]
\tikzstyle{mlnode2} = [rectangle, draw=violet!70, line width=0.4mm, dashed, fill=violet!20, rounded corners, minimum height=1.5cm, text width=3.2cm, text centered, inner sep=3.5mm]
\tikzstyle{finishnode} = [chamfered rectangle, text width=4cm, draw=black, line width=0.4mm, fill=gray!20, minimum height=1.25cm, text width=3.2cm, text centered, inner sep=4mm]
\tikzstyle{regarrow1} = [line width=0.35mm, ->, >=Stealth]
\tikzstyle{mlarrow1} = [line width=0.35mm, ->, >=Stealth, dashed] 
\tikzstyle{regarrow2} = [line width=0.35mm, ->, >=Stealth]
\tikzstyle{mlarrow2} = [line width=0.35mm, ->, >=Stealth, dashed]

    \node [regnode1, text width=5.5cm] (r1) 
    {Set objective ranges for Phase I via Algorithm \ref{algo:objranges} and \\ problem-specific improvement \\ inequalities \eqref{eq:underdose_UB}--\eqref{eq:BOT_LB}};
    \node [regnode1, right of=r1, xshift=2cm, text width=4cm] (r2) 
    { Construct the set $\epsSetI$ of $\eps$-vectors via equations \eqref{eq:grid_pts1}--\eqref{eq:grid_pts2} 
    };
    \node [regnode1, right of=r2, xshift=1.75cm, text width=3cm] (r3) 
    { Set \\ $\epsSetSelI \gets \epsSetI $ };
    \node [mlnode1, below of=r3, yshift=2cm, text width=4cm] (ml1) 
    { Pick a random $\mlprob \%$ subset $\epsSetSelI \subseteq \epsSetI$
    };
    \node [regnode1, right of=r3, xshift=2.25cm, text width=6.3cm] (r4) 
    { Solve $\mopeps$ for $\eps \in \epsSetSelI$ \\ incorporating early detection \\ filters \eqref{eq:detectinfeas} and \eqref{eq:detectrepeating}   \\[0.1cm]
    };
    \node [mlnode1, below of=r4, yshift=0.5cm, text width=6.3cm] (ml2) 
    { Find a subset \\ $\epsSetPromising \subseteq (\epsSetI \setminus \epsSetSelI)$ \\ with desirable predicted performance criteria values \\ via Algorithm \ref{algo:predictselect} and domination filter \eqref{eq:notdom} 
    };
    \node [mlnode1, below of=ml2, yshift=0.5cm, text width=6.3cm] (ml3) 
    { Solve $\mopeps$ for $\eps \in \epsSetPromising $ \\ incorporating early detection \\ filters \eqref{eq:detectinfeas} and \eqref{eq:detectrepeating}
    };
    \node [mlnode2, below of=ml3, yshift=0.5cm, text width=5cm] (ml4) 
    { Combine the results obtained with \\ $\epsSetSelI $ and $\epsSetPromising $
    };
    \node [regnode2, below of=ml4, xshift=0.0cm, yshift=1cm] (r5) 
    { Set objective ranges \\ for Phase II \\via Algorithm \ref{algo:objranges_phase2}
    };
    \node [regnode2, left of=r5, xshift=-1.5cm, text width=4cm] (r6) 
    { Construct the set $\epsSetII$ of $\eps$-vectors via equations \eqref{eq:grid_pts1}--\eqref{eq:grid_pts2}
    };
    \node [regnode2, left of=r6, xshift=-1.75cm, text width=3cm] (r7) 
    { Set \\[0.1cm]
    $\epsSetSelII \gets \epsSetII$
    };
    %
    \node [mlnode2, above of=r7, yshift=-1.25cm, text width=5.5cm] (ml5) 
    { Find a subset $\epsSetSelII \subseteq \epsSetII $ with desirable predicted performance criteria values \\ via Algorithm \ref{algo:predictselect} and domination filter \eqref{eq:notdom}
    };
    %
    %
    \node [regnode2, left of=r7, xshift=-2.5cm, text width=5.75cm] (r8) 
    { Solve $\mopeps$ for $\eps \in \epsSetSelII$, \\ incorporating early detection \\ filters \eqref{eq:detectinfeas} and \eqref{eq:detectrepeating}
    };
    \node [finishnode, above of=r8, yshift=1cm, text width=3cm, draw=black, line width=0.3mm, fill=gray!20] (fin) 
    { Combine the results from \\ Phases I \& II
    };
    \node [infonode1, left of=r1, xshift=1cm, yshift=-4.8cm] (i1) 
    {\bf P\\H\\A\\S\\E\\ \phantom{o} I
    };
    \node [infonode2, below of=i1, yshift=-5.4cm] (i2) 
    { \bf P\\H\\A\\S\\E\\ \phantom{o} II
    };
    \draw [regarrow1] ([yshift=0.25cm]r1.east) -- ([yshift=0.25cm]r2.west);
    \draw [regarrow1] ([yshift=0.25cm]r2.east) -- ([yshift=0.25cm]r3.west);
    \draw [regarrow1] ([yshift=0.25cm]r3.east) -- ([yshift=0.25cm]r4.west);
    \draw [regarrow1] (r4.east) |- ++(1cm,0mm) -| ++(0mm,0mm) |- (r5.east);
    \draw [regarrow2] ([yshift=-0.25cm]r5.west) -- ([yshift=-0.25cm]r6.east);
    \draw [regarrow2] (r6.west) -- (r7.east);
    \draw [regarrow2] ([yshift=-0.25cm]r7.west) -- ([yshift=-0.25cm]r8.east);
    \draw [regarrow2, draw=black] ([xshift=-0.5cm]r8.north) -- ([xshift=-0.5cm]fin.south);
    \draw [mlarrow1] ([yshift=-0.5cm]r2.east) |- ++(1.25cm,0mm) -| ++(0mm,0mm) |- (ml1);  
    \draw [mlarrow1] (ml1.east) -| ++(0.75cm,0mm) -| ++(0mm,0mm) |- ([yshift=-0.5cm]r4.west);
    \draw [mlarrow1] ([yshift=-0.5cm]r1.east) -- ([yshift=-0.5cm]r2.west);
    \draw [mlarrow2] ([yshift=0.5cm]r5.west) -- ([yshift=0.5cm]r6.east);
    \draw [mlarrow1] (r4.south) -- (ml2);
    \draw [mlarrow1] (ml2) -- (ml3);
    \draw [mlarrow1] (ml3) -- (ml4);
    \draw [mlarrow2] (ml4.south) -- (r5.north);
    \draw [mlarrow2] ([yshift=0.5cm]r6.west) -| ++(-0.5cm,0mm) |- ++(0mm,0mm)|- (ml5.east);
    \draw [mlarrow2] (ml5.west) -| ++(-0.5cm,0mm) |- ++(0mm,0mm) |- ([yshift=0.5cm]r8.east);
    \draw [mlarrow2, draw=black] ([xshift=0.5cm]r8.north) -- ([xshift=0.5cm]fin.south);
\end{tikzpicture}
}
\caption{Flow chart of our two-phase method, where boxes with solid borders show the regular version, filled boxes with dashed borders represent the digressions for the machine learning-guided version. (Note that \eqref{eq:detectrepeating} only applies to the case of MOLP.)} 
\label{fig:flowchart}
\end{figure}

\end{landscape}
}

%% file: PseudoCodes/algo_set_obj_ranges.tex
\newcommand{\objValSetAboveThres}{F^{\text{sel}}}

\newcommand{\tupleindex}{j}
\newcommand{\tupleset}{\mathcal{J}}
\begin{algorithm}
\small
\vspace*{0.15cm}
    \caption{Setting the objective ranges for Phase II}
    \label{algo:objranges_phase2}
    {\bf Input:} 
    {\cblue A set $\epsSet \coloneqq \{ \eps^\tupleindex\}_{\tupleindex \in \tupleset}$ of $\eps$-vectors, a tuple set 
    $\big\{ \big( \objVar^{\tupleindex}, \coverage^{\tupleindex}, \PI^{\tupleindex}, \BOTVal^{\tupleindex} \big) \big\}_{\tupleindex \in \tupleset}$,  where ${Op}^{\tupleindex} \coloneqq {Op}(\optSol(\eps^\tupleindex))$ for all operators $Op \in \{\objVar, \ \coverage, \ \PI, \ \BOTVal \}$ and  $\tupleindex \in \tupleset$,}
    and threshold values $\coverageThres, \ \PIThres, \ \BOTThres$ \\
    {\bf Output:} Lower and upper limit pairs $(\objRangeLB_{\objIndex}, \objRangeUB_{\objIndex})$ of $\objVar_{\objIndex} $'s for all $ \objIndex \in \{2, \ldots, \numObjs\}$
    \begin{algorithmic}[1]
        \Procedure{SetObjRangesPhaseII}{}
        \ForAll {$\objIndex \in \{2, \ldots, \numObjs\}$} 
            \State $\objValSetAboveThres_{\objIndex} \gets \big\{ \objVar^{\tupleindex}_{\objIndex} \colon \coverage^{\tupleindex} \geq \coverageThres, \  \PI^{\tupleindex} \geq \PIThres, \  \BOTVal^{\tupleindex} \leq \BOTThres \big\}$ 
            \State $\objRangeLB_{\objIndex} \gets \min\{\objValSetAboveThres_{\objIndex}\}, \ \objRangeUB_{\objIndex} \gets \max\{\objValSetAboveThres_{\objIndex}\} $
        \EndFor
        \EndProcedure
    \end{algorithmic}
\end{algorithm}

%% file: PseudoCodes/algo_predict_select.tex
\begin{algorithm}
\small
\vspace*{0.15cm}
    \caption{Selecting promising $\eps$-vectors via prediction}
    \label{algo:predictselect}  
    {\bf Input:} 
    {\cblue Two sets 
    $\epsSetFeasTrain \coloneqq \{ \eps^{\tupleindex_1}\}_{\tupleindex_1 \in \tuplesetfeastr}$ and $\epsSetPerfTrain \coloneqq \{ \eps^{\tupleindex_2}\}_{\tupleindex_2 \in \tuplesetperftr}$; 
    two training sets 
    $\left\{\left(\eps^{\tupleindex_1}, \ \feas^{\tupleindex_1} \right) \right\}_{\tupleindex_1 \in \tuplesetfeastr}$ and $\left\{\left(\eps^{\tupleindex_2}, \ \coverage^{\tupleindex_2}, \ \PI^{\tupleindex_2}, \ \BOTVal^{\tupleindex_2} \right) \right\}_{\tupleindex_2 \in \tuplesetperftr}$ where ${Op}^{\tupleindex_i} \coloneqq {Op}(\optSol(\eps^{\tupleindex_i}))$ for all operators $Op \in \{\feas, \ \coverage, \  \PI, \ \BOTVal \}$ and  $\tupleindex_1 \in \tuplesetfeastr, \ \tupleindex_2 \in \tuplesetperftr$;
    training functions $\trainFnc_{\text{Crit}}(.)$ for $\text{Crit} \in \{\feas, \coverage, \ \PI, \ \BOTVal \}$; 
    a set $\epsSet$ to make predictions for;
    and threshold values $\coverageThres, \PIThres, \BOTThres$}
    \\
    {\bf Output:} A subset $\epsSet^{\prime} \subseteq \epsSet $ of selected $\eps$-vectors
    \begin{algorithmic}[1]
        \Procedure{PredictAndSelect}{}
        \State $\predFnc_{\feas}(.) \gets \trainFnc_{\feas} 
        \left( \left\{\left(\eps^{\tupleindex_1}, \ \feas^{\tupleindex_1} \right) \right\}_{\tupleindex_1 \in \tuplesetfeastr} \right)$ 
        \State $\predFnc_{\coverage}(.) \gets \trainFnc_{\coverage} 
        \left( \left\{\left(\eps^{\tupleindex_2}, \ \coverage^{\tupleindex_2} \right) \right\}_{\tupleindex_2 \in \tuplesetperftr} \right)$ 
        \State $\predFnc_{\PI}(.) \gets \trainFnc_{\PI} 
        \left( \left\{\left(\eps^{\tupleindex_2}, \ \PI^{\tupleindex_2} \right) \right\}_{\tupleindex_2 \in \tuplesetperftr} \right)$ 
        \State $\predFnc_{\BOTVal}(.) \gets \trainFnc_{\BOTVal} 
        \left( \left\{\left(\eps^{\tupleindex_2}, \ \BOTVal^{\tupleindex_2} \right) \right\}_{\tupleindex_2 \in \tuplesetperftr} \right)$ 
        \State $\epsSet^{\prime} \gets \varnothing$
        \ForAll {$\eps \in \epsSet$}
            \If{$\predFnc_{\feas}(\eps)$ = 1}
            \Comment{If $\eps$ is predicted to yield a feasible solution}
                \If{$\predFnc_{\coverage}(\eps) \geq \coverageThres \algoand \predFnc_{\PI}(\eps) \geq \PIThres \algoand \predFnc_{\BOTVal}(\eps) \leq \BOTThres$}
                \State $\epsSet^{\prime} \gets \epsSet^{\prime} \cup \{\eps\}$
                \EndIf
            \EndIf
        \EndFor
        \EndProcedure
    \end{algorithmic}
\end{algorithm}%

%% file: results_OO.tex
\section{Experimental Results}
\label{sec:expresults}

We conducted our experiments on a MacOS computer with a 3 GHz Intel Core i5 CPU and 16 GB memory.
We used CPLEX 12.10 \citep{cplex1210} as linear programming solver,
and \texttt{scikit-learn} Python library \citep{scikit-learn} for our ML-related implementations.

{\crevii 
We considered combinations of the primal and dual formulations of the $\eps$-constraint formulation of our MOLP with all the applicable algorithms offered by CPLEX (including the ``automatic" option where the solver decides which algorithm(s) to use) to determine the fastest option. 
The average solution times of 77 $\varepsilon$-constraint models of a representative instance (Case 3) are provided in Table~\ref{tab:lp_algos}, which shows that solving the dual model with the primal simplex and sifting algorithms give the best overall performance. 
Since the two average times are quite close, we further compared the two options on another instance (Case 6) via 71 $\varepsilon$-constraint models. 
The average times were 4.53 and 5.67 seconds for primal simplex and sifting algorithms, respectively. 
Hence, we conducted our computational experiments by solving the dual formulation using the primal simplex method. }

\setlength{\tabcolsep}{9pt} 
\renewcommand{\arraystretch}{1.01} 
\begin{table}[!ht]
\crevii
    \centering
    \caption{\crevii Average solution times (in seconds).}
    \label{tab:lp_algos}%
    \resizebox{0.47\textwidth}{!}{
    \begin{tabular}{lrr}
    \toprule
    Algorithm & Primal model & Dual model\\
    \midrule
	Automatic &  10.16 &  7.74\\
	Primal simplex &  7.00 &  2.36\\
	Dual simplex &  3.91 &  4.96\\
	Barrier	&  27.89 &  21.08\\
	Concurrent	&  10.74 &  8.20\\
	Sifting	&  49.96 &  2.13\\
    \bottomrule
    \end{tabular}
    }
\end{table}


We begin our empirical analysis by providing descriptive information on our test bed in Section~\ref{subsec:expresults_instinfo}.
Afterwards, we present summary statistics on the performance of {\twophase} in Section~\ref{subsec:expresults_2phase}, including the impact of the early detection strategies. 
Lastly, in Section~\ref{subsec:expresults_comparison}, we evaluate the solution quality of {\twophase}, and compare it to that of the single-objective linear program ({\cevikModel}) of \cite{cevik2019simultaneous}, and the actual clinical results obtained by applying the plans generated via the planning module Leksell GammaPlan\textsuperscript{\textregistered} \citep{leksellGammaPlan2010}.

\subsection{Test Bed}
\label{subsec:expresults_instinfo}

Our test instances consist of anonymized data of eight previously treated clinical cases with a single tumor, which are also used {\crevii by} \cite{cevik2019simultaneous}.
Table~\ref{tab:inst_info} lists descriptive information about the instances.


\input{Tables/table_instance_info}
Isocenter locations are typically determined manually by the planners so that the overall tumor volume is adequately covered.
In theory, it is possible to use each tumor voxel as an isocenter, but the size of the model then grows prohibitively large and cannot be handled in practice.
Isocenters in our cases are positioned automatically with the grassfire and sphere packing algorithms~\citep{cevik2019simultaneous}.
We note that the quality of the generated treatment plans is quite robust to the number of isocenters, and potential effect of fewer isocenters is usually compensated for by the use of larger collimator sizes in delivering radiation, see for instance \citep{cevik2019simultaneous}.


\subsection{Efficacy of Our Approach}
\label{subsec:expresults_2phase}
Infeasibilities and repeating solutions are prevalent and the filters to avoid them give substantial benefit.
Phase II subsequently leads to fewer infeasibilities, which is one of our two main incentives ---the other being to focus on regions expected to yield desirable solutions.

\subsubsection{Results from the regular version.}

Table~\ref{tab:stats_regular} summarizes the results from {\twophasereg}. 
{\crevii The column labels are defined as follows:
\BI 
\I ``\# $\eps$-vectors" is the number of $\eps$-vectors constructed at the beginning of the phase (i.e., the cardinality of $\epsSetI$ or $\epsSetII$ in Figure~\ref{fig:flowchart}),
\I ``\# sol" is the number of obtained Pareto-optimal points,
\I ``\% sol" is the percentage of Pareto-optimal points with respect to the number of $\eps$-vectors,
\I ``\% infeas" is the percentage of $\eps$-vectors for which $\mopeps$ is infeasible,
\I ``\% omit in infeas" is the percentage of early-detected $\eps$-vectors among the infeasible cases,
\I ``\% omit for infeas" is the percentage of early-detected $\eps$-vectors for infeasibility (among all $\eps$-vectors),
\I ``\% omit for repeat" is the percentage of early-detected $\eps$-vectors for repeating solutions,
\I ``Unit sol time" is the average time to obtain a solution,
\I ``Total sol time" is the overall time spent solving all $\mopeps$ models.
\EI
}

\input{Tables/table_stats_regular}

We observe that early detection of infeasibilities leaves a small proportion of $\eps$-vectors to be evaluated.
For Phase I, on average, 66.4\% of $\eps$-vectors are discarded for infeasibility, and only 0.9\% of the infeasible models are later considered. 
Together with an average 28.5\% of $\eps$-vectors discarded to eliminate repeating solutions, around 95\% of all $\eps$-vectors are omitted {\crevii in a negligible computation effort requiring around $10^{-4}$ seconds (thus not explicitly reported).}
For Phase II, the average infeasibility percentage drops to 10.5\%, with more than half of those being detected early.
This shows that revision of objective ranges for Phase II concentrates the search as intended.  
The increased percentage of solutions in Phase II also indicates this fact, and the forthcoming Section~\ref{subsec:expresults_comparison} illustrates that the overall solution quality improves as planned.
Combined with the $\eps$-vectors omitted to avoid repeating solutions, an average of nearly 48\% of all $\eps$-vectors are discarded in Phase II, and we observe that the unit solution times are decreased. 
We note that, except for Cases 2 and 8, limits for tumor overdosing collapsed to zero, making $\numGrids^3$ $\eps$-vectors rather than $\numGrids^4$, with $\numGrids$ denoting the number of distinct $\eps_{\objIndex}$ values used for each objective, set to 10 and 5 for Phases I and II, respectively.

We tested {\twophasereg} without the repeating solution filter to identify its impact.
Among the $\eps$-vectors that yield a solution in Phase I, only 0.34\% gives a solution not previously obtained but discarded with the early detection filter (due to existence of alternative optima).
We encountered no such solutions in Phase II.
Early detection eliminated around 93\% and 55\% of $\eps$-vectors that produced a solution in Phases I and II, with no repeating solutions, 
so the early detection strategies significantly reduce calculations without jeopardizing Pareto representation.
Run times of some cases (e.g., Case 4) may nevertheless be high due to the size and/or inherent difficulty of those instances.

\subsubsection{Results from the ML-guided version.}
\label{subsubsec:expresults_2phase_ml}
{\crevi 
We selected the classification and regression models from several alternatives, 
by taking the results from (the first round of) Phase I for each one of the eight cases, and using them as training data.
We performed a 5-fold cross validation 10 times for each model and case combination separately, and we compared the mean accuracy (proportion of correct classifications) and the area under the receiver operating characteristic curve for classifiers, and mean absolute error for regressors.
Random forests yielded the best overall performance with decision trees doing almost as well, probably because we only have four features (which is the length of an $\eps$-vector for our 5-objective model). For details, see Table~\ref{tab:ml_models_comparison} 
in Appendix~\ref{app:results}. 
}

{\crevi 
We find upon comparing predictions with actual results that} feasibility classification is at least 96\% accurate on the average, and that mean absolute errors for Cov and PCI are approximately no greater than 2\%.
BOT values vary on a larger scale, which occasionally lead to relatively higher prediction errors.
Nevertheless, with a mean absolute error of less than seven (minutes), prediction performance for BOT remains useful. 
We note that it is possible to hedge against prediction inaccuracies by slightly relaxing the filtering threshold values, 
see Table~\ref{tab:prediction_scores}  
in Appendix~\ref{app:results}. 

{\twophaseml} decreases the run times by nearly half.
Specifically, in Phase I, with 8.6\% of the initially constructed $\eps$-vectors for {\twophasereg} on the average, {\twophaseml} 
averages a {\crevi 47.6\%}
time saving. 
The number of $\eps$-vectors selected for Phase II of {\twophaseml} is 14.5\% of those built for the same phase of the regular version, and the time savings range from {\crevi 24\% to 79\%, averaging 46.5\%},
see Table~\ref{tab:stats_ml} 
in Appendix~\ref{app:results}. 



\subsection{Solution Quality Assessment}
\label{subsec:expresults_comparison}

Table~\ref{tab:comparison} contains summary statistics for Cov, PCI, BOT, and maximum OAR dose values of the two different versions of {\twophase} compared with those of {\cevikModel} and the clinical outcomes.
Our solution approach generates a collection of solutions for each case, and we report the number of Pareto-optimal 
points obtained and the number that remained after filtering dominated Cov, PCI, and BOT values (in parentheses next to the total number of solutions).
For Cov, PCI and BOT values, we report the minimums and maximums obtained. 
We also list the average maximum dose each OAR received over the efficient solutions.

We see from Table~\ref{tab:comparison} that the maximum Cov and PCI, and the minimum BOT values are always better than the point values achieved clinically and by {\cevikModel}.
Comparing the performances of {\twophasereg} and {\twophaseml}, we observe that the latter yields average and minimum Cov and PCI values at least as good as those of the former. 
This is because {\twophaseml} has an extra filtering step to only generate desirable solutions in the second round of Phase I,  
as a consequence of which it also yields fewer solutions than {\twophasereg}.
{\crevii The minimum and maximum BOT values from {\twophasereg} and {\twophaseml} are almost identical --- the maximum BOT value allowed was three hours, as also reported in \citep{maxBOTref}.}

\input{Tables/table_comparison}

\input{Figures/fig_plots_cases1to4}

We next investigate the distribution of Pareto solutions
over the space of the three performance criteria in Figure~\ref{fig:cases1to4}.
For brevity we consider Cases 1 to 4 as representatives to show the different ways solutions are dispersed. 
Results for the remaining cases are provided in 
the Appendix. 
Each row of Figure~\ref{fig:cases1to4} is dedicated to one case and contains two plots.
The plot on the left shows Cov versus PCI and the plot on the right displays BOT versus PCI for solutions yielding at least 99.7\% Cov.
This limit is higher than the maximum coverage achieved clinically and by {\cevikModel} over all eight cases, guaranteeing that our results on the right-hand side plots outperform clinical outcomes and {\cevikModel} results. 

Figure~\ref{fig:cases1to4} illustrates that the two phases of {\twophase} work as intended. 
While the solutions from Phase I are relatively scattered throughout the space of the three performance criteria, those from Phase II are typically clustered around more desirable regions.
The left-hand side plots demonstrate that {\twophasereg} is able to produce a considerably large collection of solutions outperforming the others in terms of both Cov and PCI, and {\twophaseml} identifies many such solutions.
The right-hand side plots show that several of our solutions are significantly better than those from the clinical design or from {\cevikModel}.
Cases 2 and 4 are exceptions.
A Cov of only 82\% has been achieved previously for Case 2, potentially due to clinicians' preference to keep the dose to OAR1 much lower than the predefined OAR1 dose limit (9.7 Gy vs 15.0 Gy). 
Our approach, however, produces solutions attaining essentially full coverage with higher, yet still realizable, BOT values, and hence, we provide the clinician options. 
Moreover, we obtain many solutions with strictly better Cov, PCI, and BOT values than the others, which are not contained in the right-hand side plot because points with less than 99.7\% Cov are not shown.
Almost all solutions for Case 4
are better than the others in terms of Cov and PCI, yet they yield relatively higher BOT values than that of {\cevikModel}.
The general distribution of solutions from Cases 5 to 8 are similar to that of Case 3, and thus our conclusions remain valid for them. 

\input{Figures/fig_DVH_isodose}

{\crevi
We lastly examine the dose-volume histograms (DVHs) and isodose curves for a representative case, Case 3, see Figure~\ref{fig:DVH_isodose}. 
A DVH displays the fraction of each structure's volume receiving at least a certain amount of dose. We plot the relative dose percentages on the $x$-axis by proportioning the received doses to the prescription dose for the tumor and the maximum allowed doses for OARs. 
Among the set of solutions obtained from {\twophasereg}, {\crevii we first remove solutions whose Cov is less than 99.7\% and whose PCI is less than 0.75.
We then discard solutions dominated in terms of the Cov, PCI and BOT.}
We plot the DVH curves for the remaining solutions that provide virtually full coverage, PCI values between 0.87 and 0.96, and BOT from 31 to 115 minutes.
{\crevii We then select a solution with 0.96 PCI and 101-minute BOT}
and plot the associated isodose curves for several cross sections of the tumor and OAR tissues.
For each cross section, we contour the areas that receive at least 50\% and 100\% of the prescribed dose for the tumor (12.5~Gy for Case~3), and indicate the tumor and nearby OAR tissues with shaded areas.
}

{\crevi A DVH curve for the tumor should ideally mimic a vertical decay at the prescription dose, and in our case it is deviating only slightly from the reference line (shown with the gray vertical line).
Due to the proximity of the tumor to the brainstem and cochlea, which is evident from the cross sections depicted in the first three rows of the plots on the right in Figure~\ref{fig:DVH_isodose}, the 50\% isodose curves for the target intersect these organs. 
As the intersecting volume percentage of the cochlea is larger, we observe that its DVH curve starts decaying later than that of the brainstem.
However, the clinically prioritized measure for OARs is the maximum doses they receive, and they are less than 85--90\% of the allowed values.
}

%% file: Tables/table_instance_info.tex
\setlength{\tabcolsep}{3pt}
\begin{table}[htbp]
  \centering
  \caption{Instance information.}
\resizebox{0.95\textwidth}{!}{
    \begin{tabular}{c S[table-format=2.0] S[table-format=2.1] c S[table-format=2.2] S[table-format=5.0] S[table-format=5.0] S[table-format=5.0] S[table-format=5.0]}
    \toprule\\[-0.5cm]
    \parbox{1.5cm}{\centering Case} & \parbox{2cm}{\centering \# isocenters} & \parbox{2.3cm}{\centering Prescribed\\ dose (Gy)} & \parbox{2.5cm}{\centering OAR dose limit (Gy)} & \parbox{2.5cm}{\centering Tumor volume (cm$^3$)} & \parbox{2cm}{\centering \# tumor voxels} & \parbox{2cm}{\centering \# ring voxels} & \parbox{2cm}{\centering \# OAR voxels} & \parbox{2cm}{\centering Total \# voxels} \\[0.4cm]
    \midrule\\[-0.7cm]
    1     & 24    & 14.0    & $\{15, 15\}$   & 4.52  & 4524  & 695   & 10068 & 15287 \\
    2     & 25    & 24.0    & $\{15, 8\}$    & 5.23  & 5234  & 6657  & 8934  & 20825 \\
    3     & 19    & 12.5    & $\{15, 15\}$   & 1.35  & 1347  & 3089  & 5844  & 10280 \\
    4     & 34    & 17.0    & {--}           & 12.83 & 12829 & 18172 & {--}  & 31001 \\
    5     & 8     & 20.0    & {--}           & 1.63  & 1628  & 3809  & {--}  & 5437 \\
    6     & 20    & 12.0    & $\{15, 11.5\}$ & 2.60  & 2601  & 4077  & 8108  & 14786 \\
    7     & 33    & 12.5    & $\{15, 7\}$    & 2.52  & 2521  & 4235  & 7643  & 14399 \\
    8     & 16    & 18.0    & {--}           & 3.45  & 3447  & 6833  & {--}  & 10280 \\
    \bottomrule%
    \end{tabular}%
    }
  \label{tab:inst_info}%
\end{table}%

%% file: Tables/table_stats_regular.tex
\begin{table}[htbp]
  \centering
  \caption{Statistics from the regular version of our two-phase method.}
  \scalebox{0.7}{
    \begin{tabular}{cc S[table-format=5.0] S[table-format=2.0] S[table-format=2.1] S[table-format=2.1] >{\columncolor[gray]{0.95}}S[table-format=2.1] >{\columncolor[gray]{0.95}}S[table-format=2.1] >{\columncolor[gray]{0.95}}S[table-format=2.1] S[table-format=3.1] S[table-format=5.1]}
    \toprule \\[-0.7cm]
       \parbox{1cm}{\centering \phantom{o}}& 
       \parbox{1.75cm}{\centering Case} & 
       \parbox{1.75cm}{\centering \# \\ $\eps$-vectors} & 
       \parbox{1.5cm}{\centering \# \\ sol} & 
       \parbox{1.5cm}{\centering \% \\ sol} & 
       \parbox{1.5cm}{\centering \% \\infeas} & 
       \parbox{1.75cm}{\centering \% omit in infeas} & 
       \parbox{1.75cm}{\centering \% omit for infeas} & 
       \parbox{1.85cm}{\centering \% omit for repeat} & 
       \parbox{2cm}{\centering Unit sol time (sec)} & 
       \parbox{2cm}{\centering Total sol time (sec)} \\[0.5cm]
    \midrule
    \multirow{8}[0]{*}{Phase I} 
          & 1  & 1000  & 15    & 1.5   & 83.1  & 98.7  & 82.0   & 15.4  & \crevi 10.9  & \crevi 163.5 \\
          & 2  & 10000 & 87    & 0.9   & 75.5  & 99.8  & 75.4   & 23.6  & \crevi 16.3  & \crevi 1419.8 \\
          & 3  & 1000  & 77    & 7.7   & 57.7  & 98.1  & 56.6   & 34.6  & \crevi 2.4   & \crevi 183.0 \\
          & 4  & 1000  & 57    & 5.7   & 78.0  & 98.6  & 76.9   & 16.3  & \crevi 108.2 & \crevi 6165.2 \\
          & 5 & 1000  & 44    & 4.4   & 61.0  & 98.2  & 59.9   & 34.6  & \crevi 0.5   & \crevi 23.7 \\
          & 6 & 1000  & 71    & 7.1   & 57.1  & 98.1  & 56.0   & 35.8  & \crevi 4.5  & \crevi 321.7 \\
          & 7 & 1000  & 52    & 5.2   & 57.9  & 98.4  & 57.0   & 36.9  & \crevi 10.3  & \crevi 537.0 \\
          & 8 & 10000 & 103   & 1.0   & 67.9  & 99.9  & 67.8   & 31.1  & \crevi 2.3   & \crevi 241.1 \\
    \cmidrule(lr){1-11}
Average   &   &       & 63.3  & 4.2   & 67.3  & 98.7   & 66.4  & 28.5  & \crevi 19.4  & \crevi 1131.9 \\
    \midrule
    \multirow{8}[0]{*}{Phase II} 
          & 1  & 25    & 2     & 8.0   & 0.0   & 0.0  & 0.0   & 92.0  & \crevi 16.1  & \crevi 32.2 \\
          & 2  & 125   & 50    & 40.0  & 0.0   & 0.0  & 0.0   & 60.0  & \crevi 0.8  & \crevi 41.6 \\
          & 3  & 25    & 19    & 76.0  & 4.0   & 0.0  & 0.0   & 20.0  & \crevi 2.1  & \crevi 39.7 \\
          & 4  & 25    & 15    & 60.0  & 0.0   & 0.0  & 0.0   & 40.0  & \crevi 33.6  & \crevi 504.4 \\
          & 5 & 25    & 17    & 68.0  & 24.0  & 66.7  & 16.0  & 8.0   & \crevi 0.4   & \crevi 7.0 \\
          & 6 & 25    & 21    & 84.0  & 16.0  & 50.0  & 8.0   & 0.0   & \crevi 3.9  & \crevi 82.2 \\
          & 7 & 25    & 17    & 68.0  & 16.0  & 50.0  & 8.0   & 16.0  & \crevi 9.2  & \crevi 155.8 \\
          & 8 & 125   & 17    & 13.6  & 24.0  & 93.3  & 22.4  & 62.4  & \crevi 3.1 & \crevi 52.3 \\
    \cmidrule(lr){1-11}
Average   &   &       & 19.8  & 52.2  & 10.5  & 32.5  & 6.8   & 37.3  & \crevi 8.7  & \crevi 114.4 \\
    \bottomrule
    \end{tabular}%
    }
  \label{tab:stats_regular}%
\end{table}%

%% file: Tables/table_comparison.tex
\afterpage{
\clearpage
\begin{landscape}
\begin{table}[htbp]
  \centering
  \caption{Performance summary of our two-phase method in comparison to {\cevikModel} and clinical results.}
  \scalebox{0.74}{
    \begin{tabular}{c cccccc cccc cccc}
    \toprule\\[-0.75cm]
    & \multicolumn{6}{c}{\twophase}   
    & \multicolumn{4}{c}{\cevikModel} & \multicolumn{4}{c}{Clinical} \\
    \cmidrule(lr){2-7} \cmidrule(lr){8-11} \cmidrule(lr){12-15} \\[-0.5cm]
    \parbox{1.5cm}{\centering Case} & 
     \parbox{1.25cm}{\centering } &
    \parbox{1.75cm}{\centering \# sol} & 
    \parbox{2cm}{\centering Cov} & 
    \parbox{2cm}{\centering PCI} & 
    \parbox{2.25cm}{\centering BOT} & 
    \parbox{2.5cm}{\centering Max OAR dose (Gy)} & 
    \parbox{1.25cm}{\centering Cov} & 
    \parbox{1.25cm}{\centering PCI} & 
    \parbox{1.25cm}{\centering BOT} & 
    \parbox{2cm}{\centering Max OAR dose (Gy)} & 
    \parbox{1.25cm}{\centering Cov} & 
    \parbox{1.25cm}{\centering PCI} & 
    \parbox{1.25cm}{\centering BOT} & 
    \parbox{2.5cm}{\centering Max OAR dose (Gy)} \\[0.3cm]
    \midrule \\[-0.5cm]
    %
    1  & R & 17 (6) & [0.809, 1] & [0.81, 1] & [10.2, 24.8]  & \{15.0, 4.0\} & 0.988 & 0.94  & 57.6  & \{16.4, 3.3\}   & 0.977 & 0.91  & 56.0    & \{15.0, 15.0\} \\
    & ML & 12 (6) & [0.871, 1] & [0.87, 1] & [10.2, 24.8] & \{15.0, 4.0\}  \vspace*{0.3cm} \\[0.6cm]
    2  & R & 137 (17) & [0.762, 1] & [0.71, 0.94] & [47.0, 180.0]  & \{15.0, 8.0\} & 0.816 & 0.79  & 100.8 & \{17.5, 7.8\}   & 0.820  & 0.56  & 156.0   & \{9.7, 7.9\} \\
    & ML & 19 (11) & [0.831, 1] & [0.76, 0.94] & [47.0, 180.0] & \{15.0, 8.0\}  \vspace*{0.3cm} \\[0.6cm]
    3  & R & 96 (16) & [0.679, 1]  & [0.68, 0.96] & [17.9, 128.2]  & \{13.7,  12.5\}  & 0.968 & 0.75  & 80.1  & \{15.2, 10.7\}  & 0.970  & 0.71  & 97.0  & \{13.6, 11.7\} \\
    & ML & 27 (9) & [0.733, 1] & [0.73, 0.96] & [17.9, 128.1] & \{13.8, 12.5\}  \vspace*{0.3cm}\\[0.6cm]
    4  & R & 72 (20) & [0.706, 1]  & [0.66, 0.91]  & [64.1, 180.0]  & {--}  & 0.995 & 0.76  & 54.4  & {--} & 0.996 & 0.68  & 138.0   &  {--} \\
    & ML & 30 (12) & [0.748, 1]  & [0.74, 0.90]  & [64.1, 180.0] & {--}  \vspace*{0.3cm} \\[0.6cm]
    5 & R & 61 (18) & [0.767, 1]  & [0.65, 0.95]  & [12.5, 78.2]  & {--}  & 0.985 & 0.84  & 20.1  & {--}  & 0.992 & 0.75  & 28.2  &  {--}  \\
    & ML & 32 (14) & [0.767, 1] & [0.77, 0.95] & [12.5, 78.2] & {--}  \vspace*{0.3cm} \\[0.6cm]
    6 & R & 92 (16) & [0.740, 1]  & [0.74, 0.96] & [12.0, 88.8]  & \{12.9,  10.8\}   & 0.985 & 0.84  & 32.0 & \{13.6, 6.9\} & 0.990  & 0.8   & 64.8  & \{15, 11.5\} \\
    & ML & 31 (11) & [0.779, 1] & [0.78, 0.96] & [12.1, 88.8]  & \{12.9, 11.1\}  \vspace*{0.3cm} \\[0.6cm]
    7 & R & 69 (17) & [0.827, 1] & [0.77, 0.99] & [14.9, 114.3]  & \{12.8, 6.7\} & 0.987 & 0.85  & 43.6  & \{13.7, 1.3\} & 0.991 & 0.83  & 64.4  & \{15.0, 7.0\} \\
    & ML & 31 (13) & [0.827, 1] & [0.83, 0.99] & [15.0, 114.1] & \{12.6, 6.8\}  \vspace*{0.3cm} \\[0.6cm]
    8 & R & 120 (17) & [0.802, 1]  & [0.77, 0.96] & [17.4, 127.5] & {--}   & 0.985 & 0.80  & 42.7 & {--} & 0.989 & 0.81  & 46.9  & {--} \\
    & ML & 51 (16) & [0.832, 1] & [0.77, 0.96] & [17.4, 127.5] & {--}   \vspace*{0.2cm} \\[0.5cm]
    \bottomrule
    \end{tabular}%
    }
  \label{tab:comparison}%
\end{table}%
\end{landscape}
}

%% file: Figures/fig_plots_cases1to4.tex
\afterpage{
\clearpage
\begin{figure}[!p]
    \centering
	\begin{subfigure}{0.4\textwidth}
		\centering
	    \includegraphics[scale=0.475]{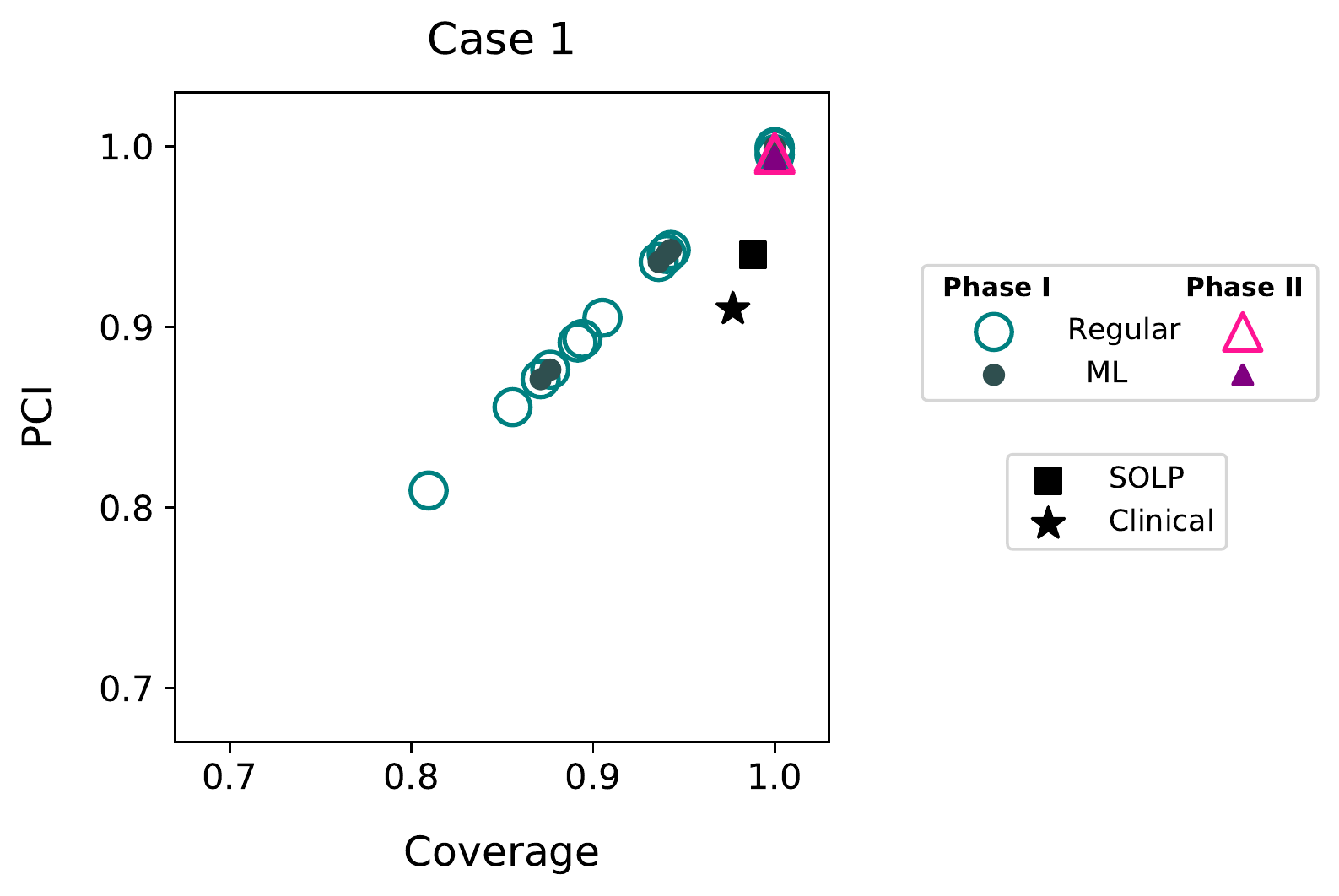}
	\end{subfigure}
	~\medskip
	\begin{subfigure}{0.4\textwidth}
		\centering
        \includegraphics[scale=0.475]{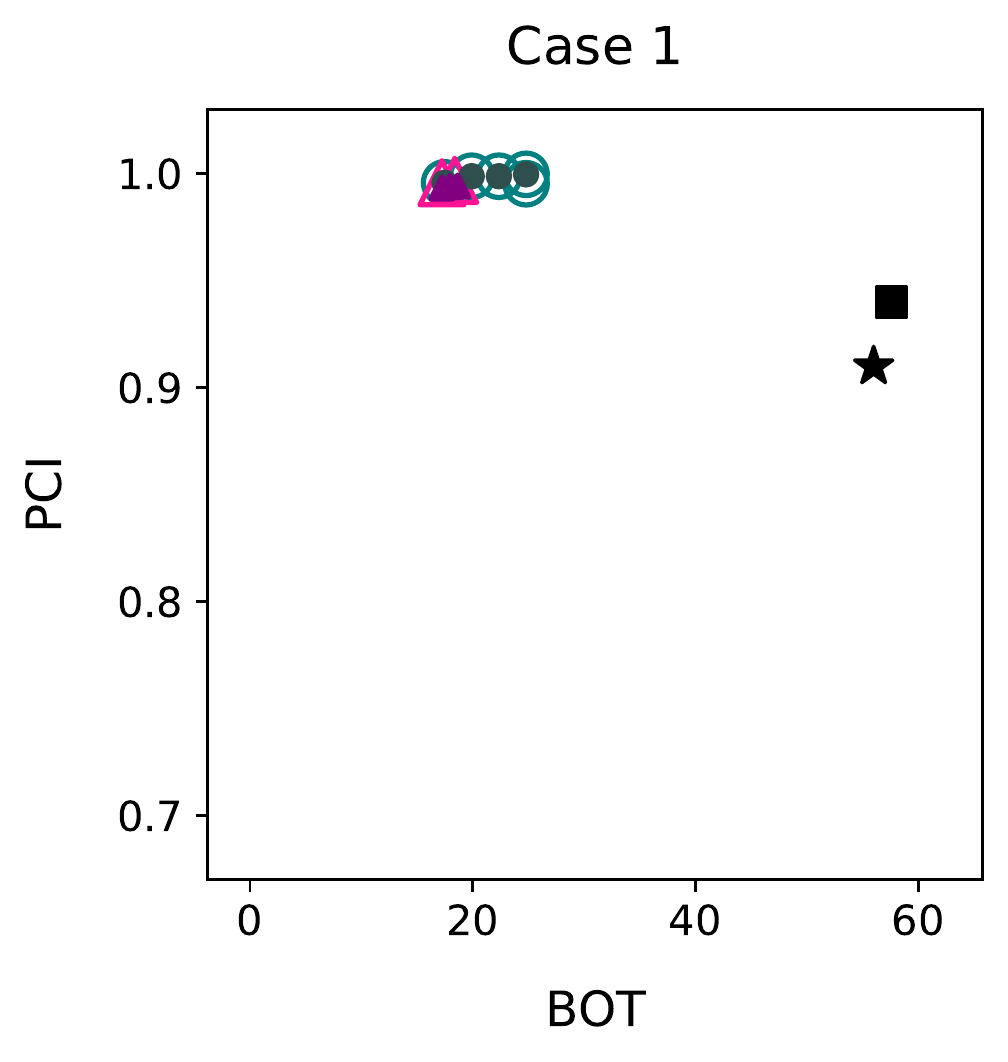}
	\end{subfigure}
    ~\medskip
	\begin{subfigure}{0.4\textwidth}
		\centering
	    \includegraphics[scale=0.475]{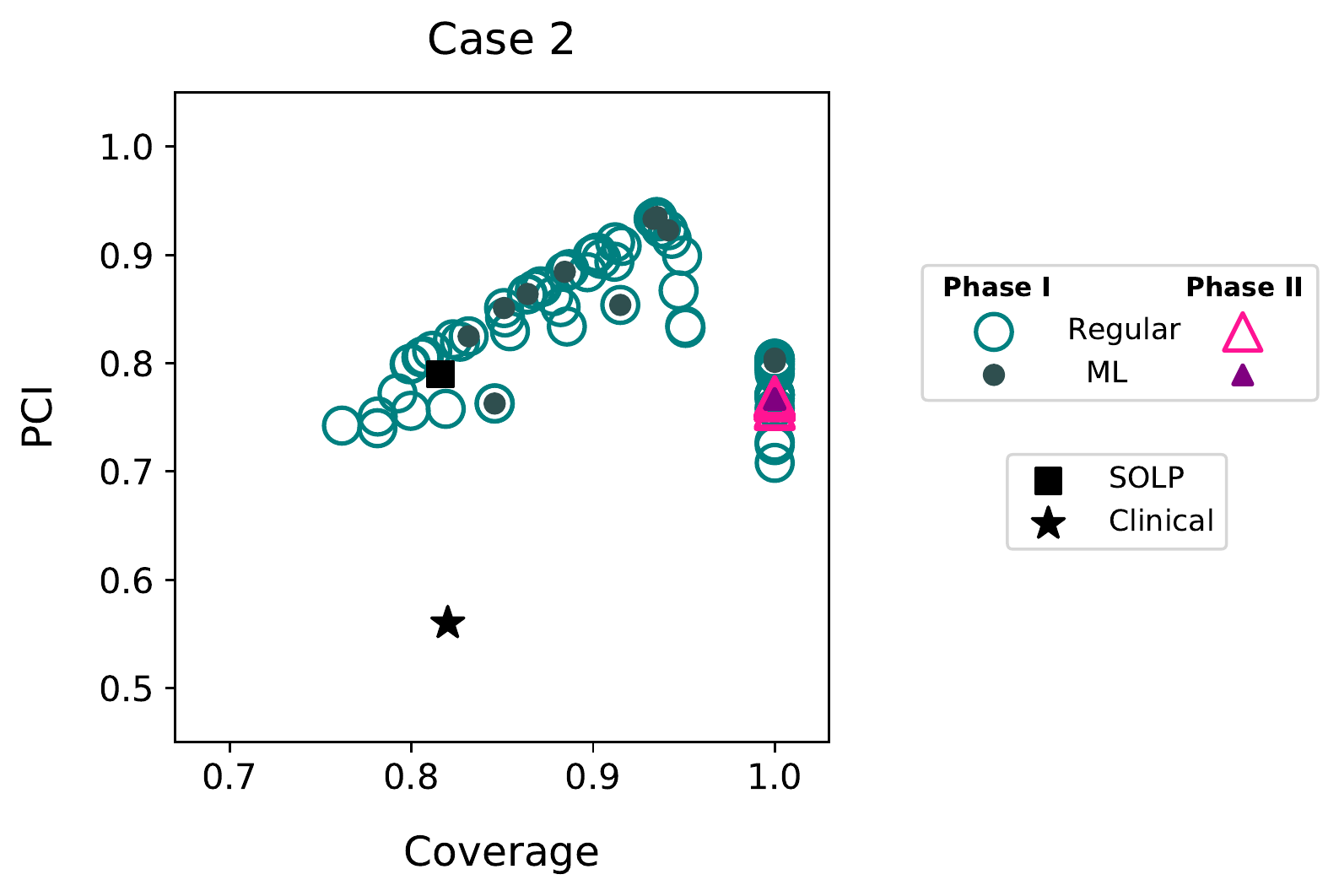}
	\end{subfigure}
	~
	\begin{subfigure}{0.4\textwidth}
		\centering
        \includegraphics[scale=0.475]{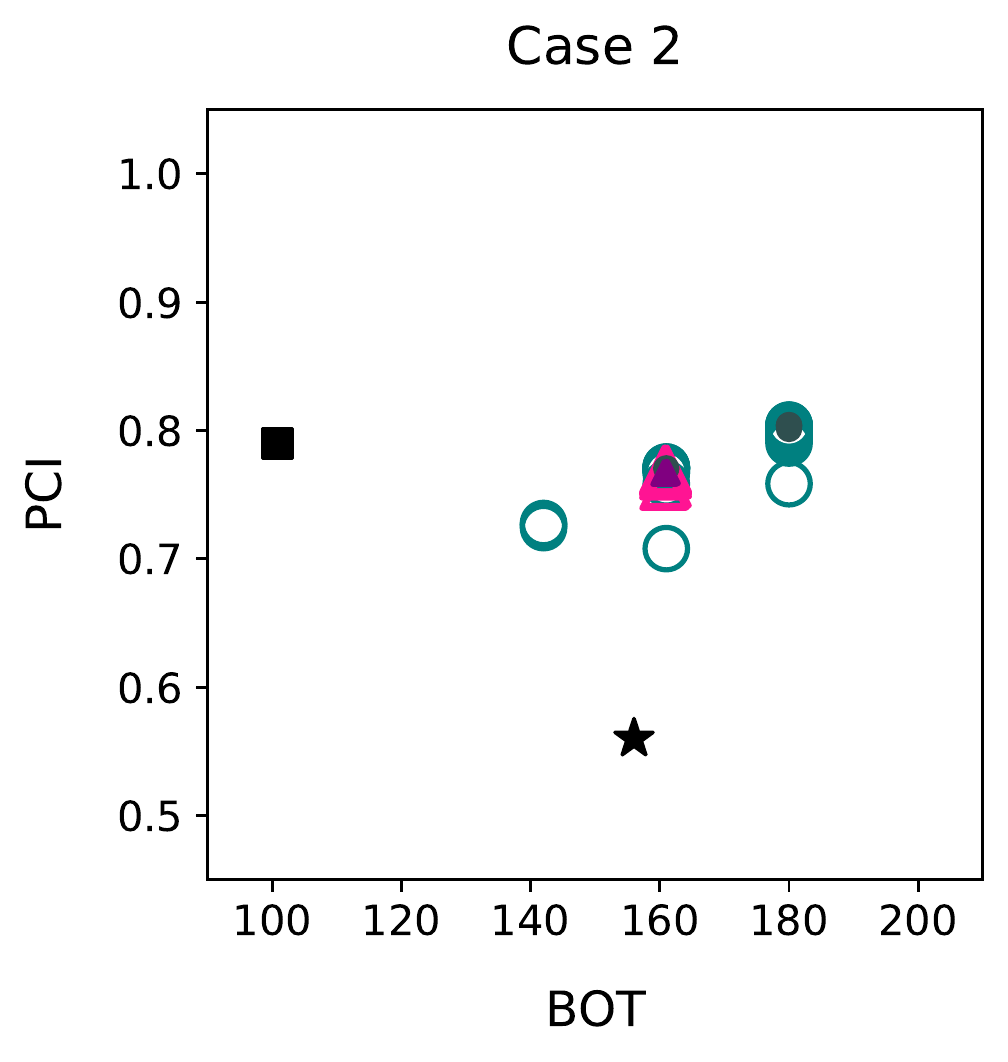}
	\end{subfigure}
	~\medskip
	\begin{subfigure}{0.4\textwidth}
		\centering
	    \includegraphics[scale=0.475]{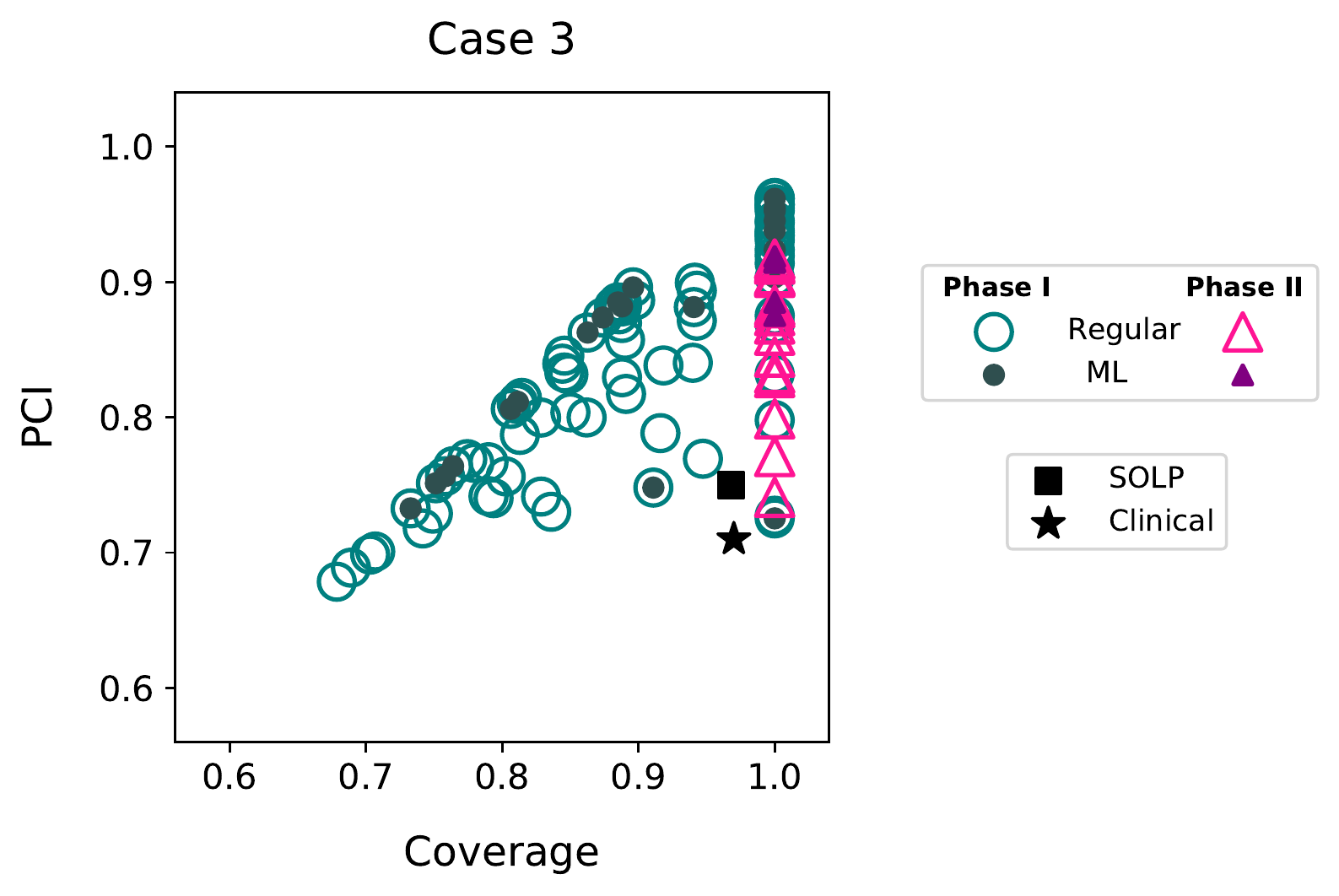}
	\end{subfigure}
	~
	\begin{subfigure}{0.4\textwidth}
		\centering
        \includegraphics[scale=0.475]{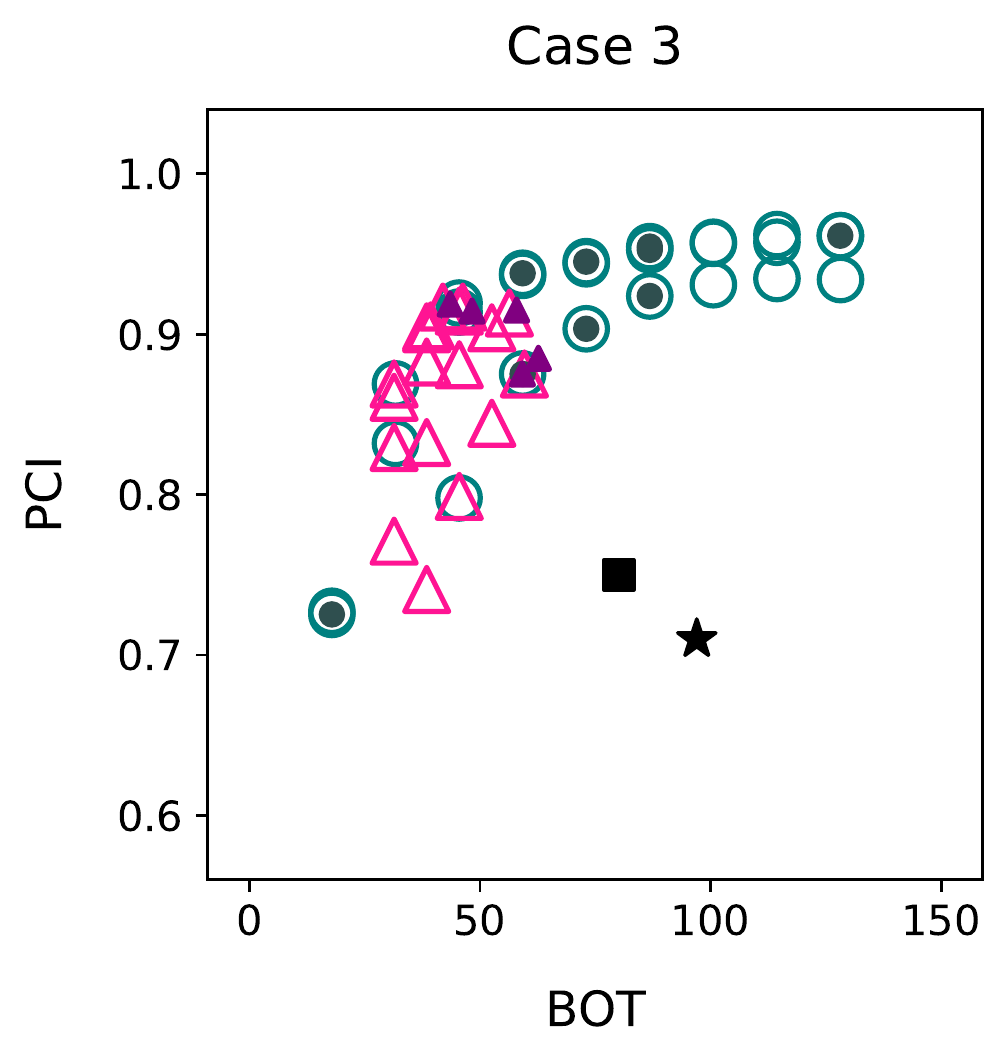}
	\end{subfigure}
	~
	\begin{subfigure}{0.4\textwidth}
		\centering
	    \includegraphics[scale=0.475]{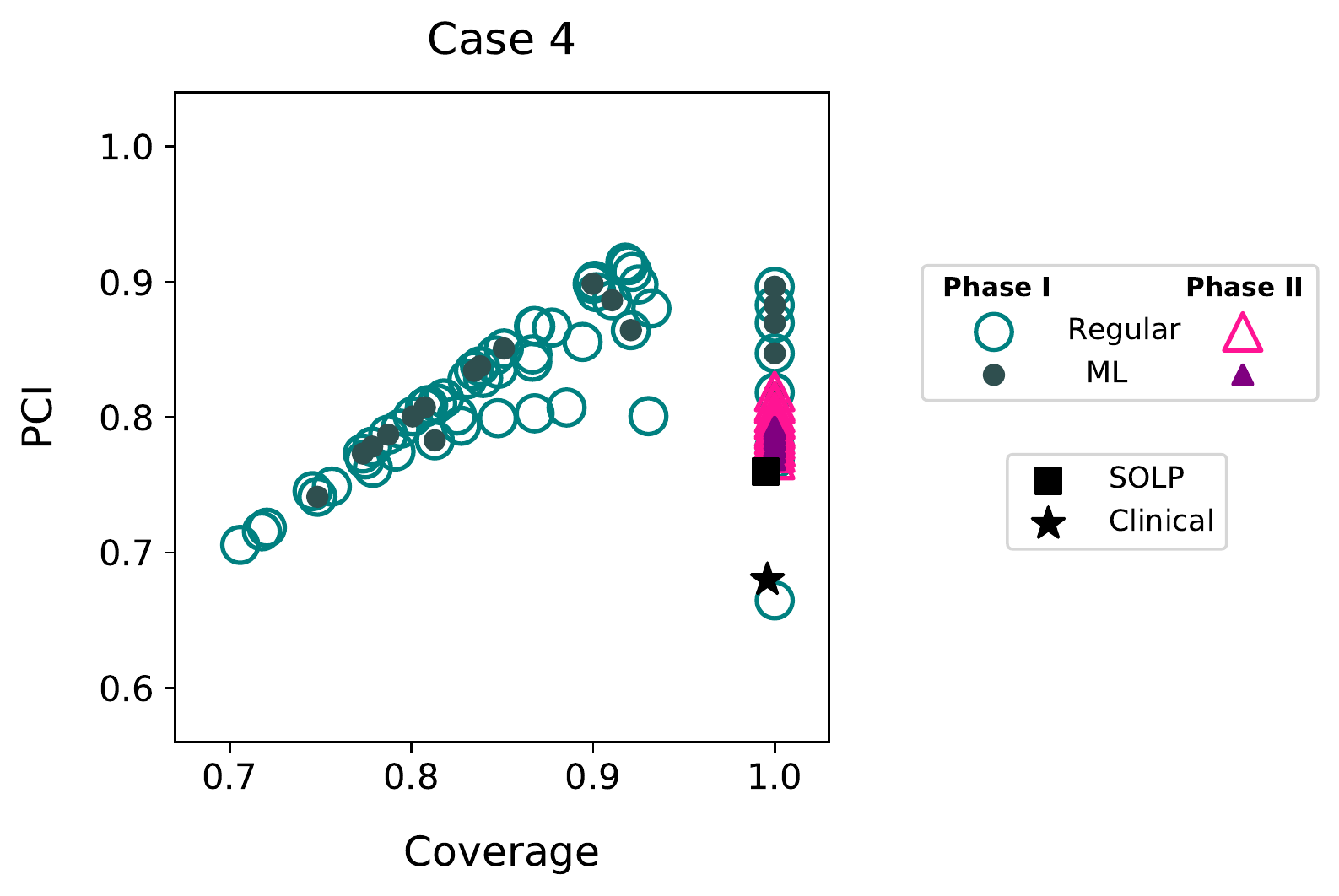}
	\end{subfigure}
	~
	\begin{subfigure}{0.4\textwidth}
		\centering
        \includegraphics[scale=0.475]{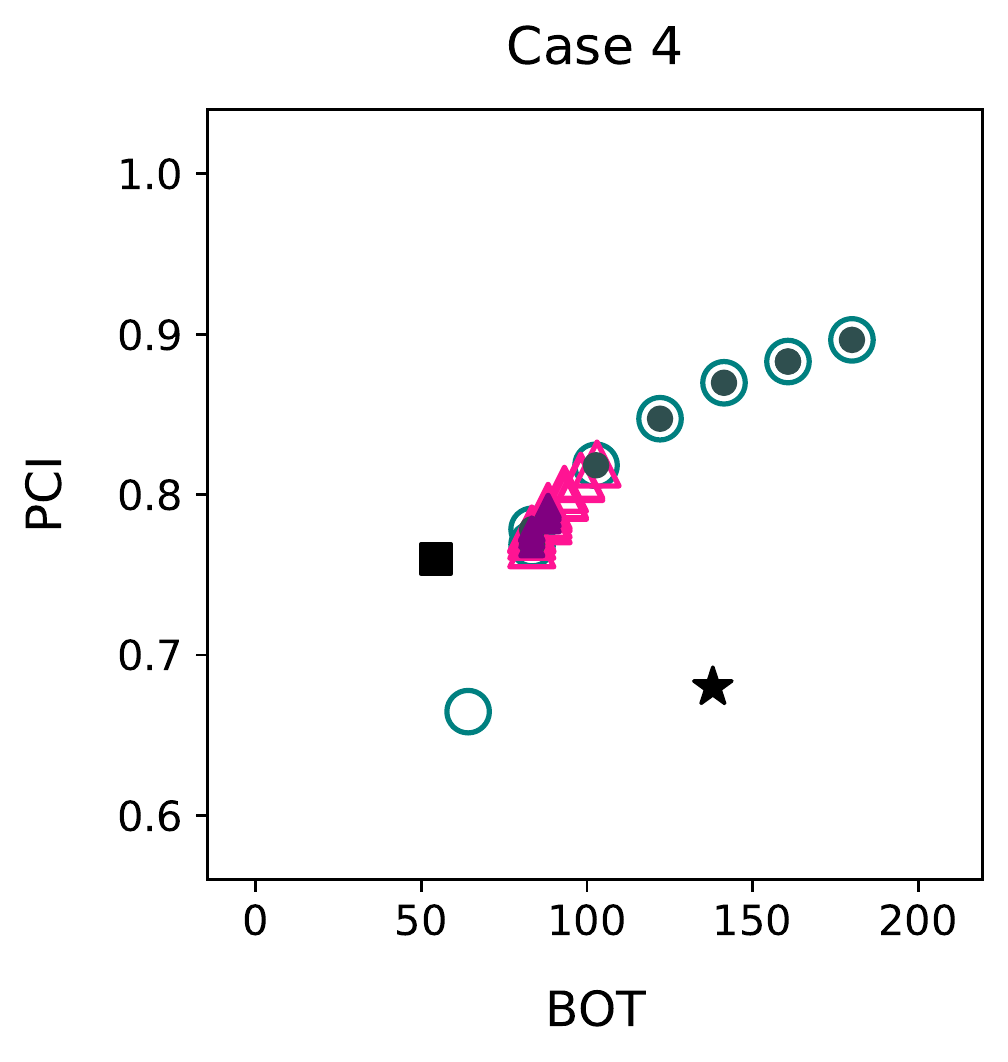}
	\end{subfigure}
	\caption{Coverage vs. PCI values of all generated solutions (left), and BOT vs. PCI values of solutions with at least 99.7\% coverage (right), for Cases 1 to 4. 
	}
	\label{fig:cases1to4}
\end{figure}%
}

%% file: Figures/fig_DVH_isodose.tex
\begin{figure}[!h]
    \centering
	\begin{subfigure}{0.43\textwidth}
	\centering
	    \includegraphics[scale=0.56]{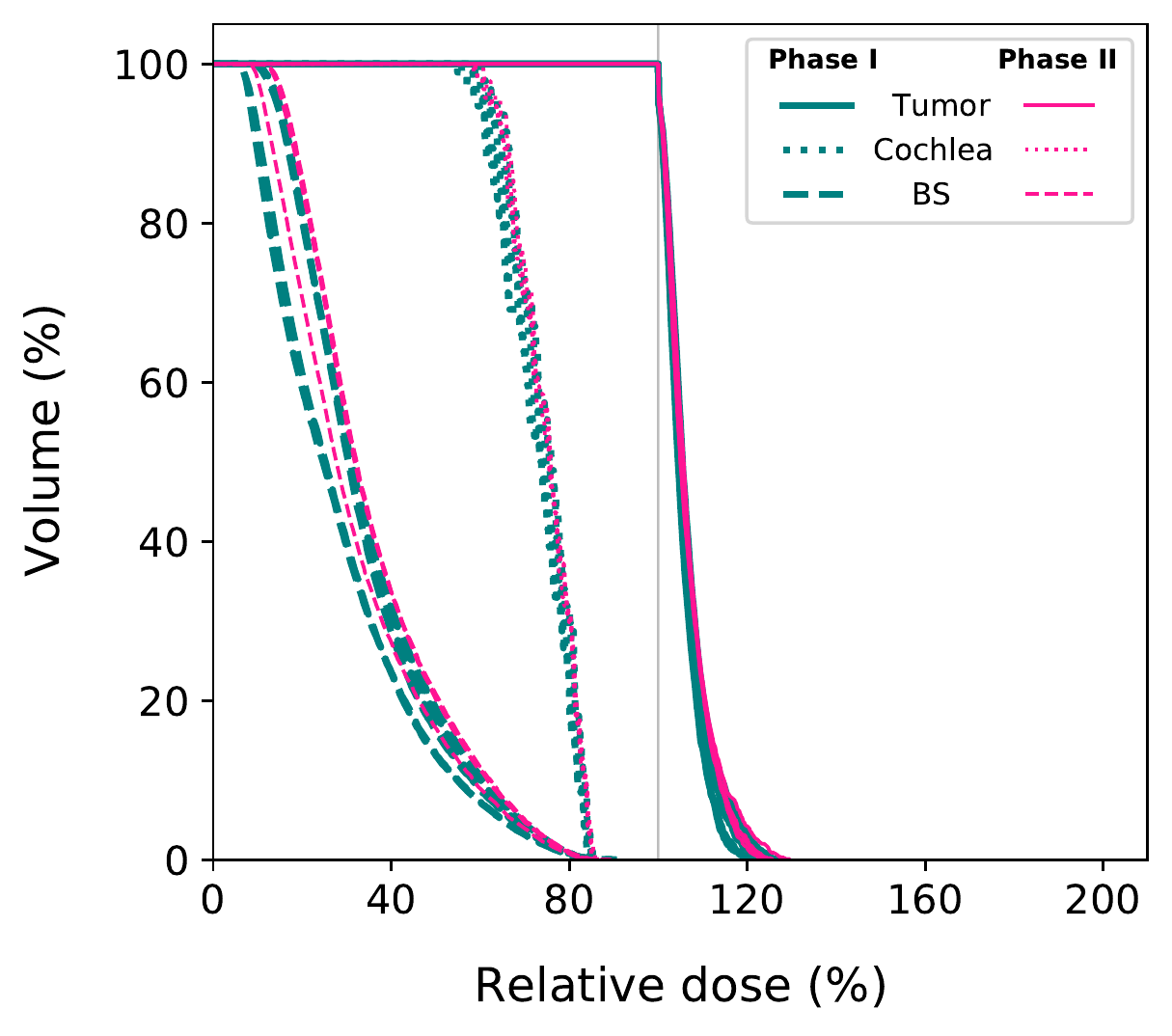} 
	\end{subfigure}
	~\hspace{0.9cm}
	\begin{subfigure}{0.43\textwidth}
		\centering
        \includegraphics[scale=0.52]{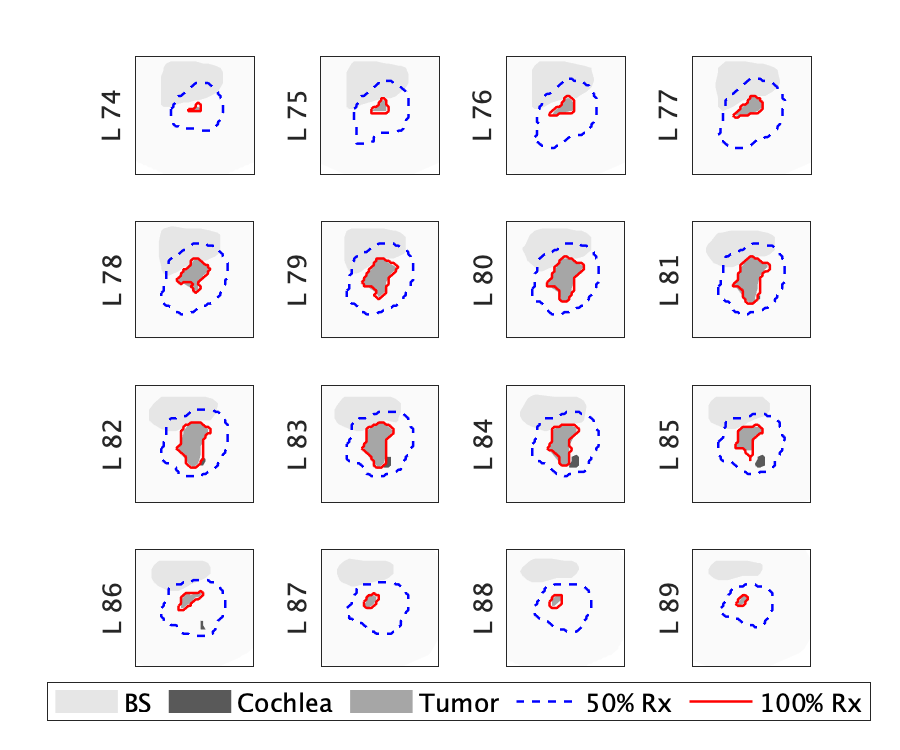} 
	\end{subfigure}
	\caption{Dose-volume histograms (left) and isodose curves (right) for Case 3. Relative dose percentages are computed with respect to the prescribed dose for the tumor (Rx), and maximum allowed doses for OARs, i.e., brainstem (BS) and cochlea. }
	\label{fig:DVH_isodose}
\end{figure}

%% file: conclusion.tex
\section{Conclusion}
\label{sec:conc}

In this study, we consider the sector duration optimization problem arising in stereotactic radiosurgery with {\crevi {\gk} delivery systems}.
We have proposed two-phase $\eps$-constraint methods to help address the trade-offs among different clinical goals.
A significant portion of our solutions 
outperformed clinical outcomes and those achieved with a single-objective model from the literature.
Machine learning algorithms reduced the overall computation time by almost half, while still capturing a sufficient amount of clinically desirable solutions.
{\crevi Our two-phase strategy is generic; it can be modified by incorporating other clinical metrics, and it can be applied to other radiotherapy or general multiobjective problems.}
Future research directions include development of multiobjective models for simultaneous optimization of isocenter selection and sector duration optimization, as well as applications of our generic computational schemes for other radiotherapy problems.

%% file: appendix_proofs_OO.tex
\section{Proofs}
\label{app:proofs}

In this section, we provide the proofs of all the propositions we present in the paper, as well as their statements for the sake of completeness. We number the equations in this section with a prefix of “A.”.

\subsection{Proof of Proposition \ref{prop:efficiency}}\label{app:efficiency_proof}
{\it \input{Propositions/prop_eps_cons_correctness}}
\input{Propositions/proof_eps_cons_correctness}

\subsection{Proof of Proposition \ref{prop:underdose_UB}}
{\it \input{Propositions/prop_underdose_UB_appendix}}
\input{Propositions/proof_underdose_UB}

\subsection{Proof of Proposition \ref{prop:BOT_LB}}
{\it \input{Propositions/prop_BOT_LB_appendix}}
\input{Propositions/proof_BOT_LB}

\subsection{Proof of Proposition \ref{prop:equiv1}}
{\it \input{Propositions/prop_repeating_sols_1}}
\input{Propositions/proof_repeating_sols_1}

%% file: Propositions/proof_eps_cons_correctness.tex
\proof{Proof.}
Let $(\optSol, \optSlack)$ be an optimal solution to \eqref{model:eps_cons_eff}. 
Assume for a contradiction that $\optSol$ is not an efficient solution of {\mop}. 
Then, there exists another feasible solution $(\allVars^{\prime}, \slack^{\prime})$ of \eqref{model:eps_cons_eff} such that 
$\objVar_{\objIndex}^{\prime} \coloneqq \objVar_{\objIndex}(\allVars^{\prime}) \leq \objVar_{\objIndex}(\optSol) \eqqcolon \objVarOptVal_{\objIndex}$ for all $\objIndex \in \{1, \ldots, \numObjs\}$ 
with at least one strict inequality. 
Using the equality constraints of \eqref{model:eps_cons_eff}, we equivalently write our assumption as
\begin{alignat*}{5}
&& \objVar_{1}^{\prime} \ && \ \leq \ & \ \objVarOptVal_{1} \\
& \eps_{2} \ & - \ \slack_{2}^{\prime} \ && \ \leq \ & \ \eps_{2} \ - \ & \optSlack_{2} \\
&&&& \vdotswithin{\leq} &\\
& \eps_{\numObjs} \ & - \ \slack_{\numObjs}^{\prime} \ && \ \leq \ & \ \eps_{\numObjs} \ - \ & \optSlack_{\numObjs},
\end{alignat*}
where at least one of the inequalities is strict. We obtain by summing these inequalities that 
\begin{alignat}{5}
 \objVar_{1}^{\prime} \ - \ \sum\limits_{\objIndex=2}^{\numObjs} \slack_{\objIndex}^{\prime}  \ < \ \objVarOptVal_{1}  \ - \ \sum\limits_{\objIndex=2}^{\numObjs} \optSlack_{\objIndex}. \label{eq:strict_ineq_sum}
\end{alignat}

\noindent Since $\optSol$ is optimal and $\allVars^{\prime}$ is feasible for \eqref{model:eps_cons_eff}, we have
\begin{alignat*}{5}
 \objVar_{1}^{\prime} \ - \ \sum\limits_{\objIndex=2}^{\numObjs} \epsCoefObj_{\objIndex} \slack_{\objIndex}^{\prime}  \ \geq \ \objVarOptVal_{1}  \ - \ \sum\limits_{\objIndex=2}^{\numObjs} \epsCoefObj_{\objIndex} \optSlack_{\objIndex}. 
\end{alignat*}

\noindent Rearranging the terms, we get
\begin{alignat}{5}
 \objVar_{1}^{\prime} \ - \ \objVarOptVal_{1}  \ \geq \ \sum\limits_{\objIndex=2}^{\numObjs} \epsCoefObj_{\objIndex} \slack_{\objIndex}^{\prime} \ - \ \sum\limits_{\objIndex=2}^{\numObjs} \epsCoefObj_{\objIndex} \optSlack_{\objIndex} \ = \ \sum\limits_{\objIndex=2}^{\numObjs} \epsCoefObj_{\objIndex} (\slack_{\objIndex}^{\prime} - \optSlack_{\objIndex}). \label{eq:objval_diff}
\end{alignat}

\noindent Since $\objVar_{1}^{\prime} \leq  \objVarOptVal_{1}$ by assumption, we have the following two cases to consider:
\medskip
\begin{compactenum}[ (i) ]
\item If $\objVar_{1}^{\prime} = \objVarOptVal_{1}$, then $ \sum\limits_{\objIndex=2}^{\numObjs} (\slack_{\objIndex}^{\prime} - \optSlack_{\objIndex}) > 0$ must be true as we need to have at least one strict inequality, necessitating inequality \eqref{eq:strict_ineq_sum} to hold. 
This implies that $ {\sum\limits_{\objIndex=2}^{\numObjs} \epsCoefObj_{\objIndex} (\slack_{\objIndex}^{\prime} - \optSlack_{\objIndex}) > 0} $. 
On the other hand, we have $0 = \objVar_{1}^{\prime} - \objVarOptVal_{1} \geq \sum\limits_{\objIndex=2}^{\numObjs} \epsCoefObj_{\objIndex} (\slack_{\objIndex}^{\prime} - \optSlack_{\objIndex}) > 0$ by inequality \eqref{eq:objval_diff}, which is a contradiction. 
\item If $\objVar_{1}^{\prime} < \objVarOptVal_{1}$, the condition that at least one of the inequalities is strict is already satisfied and the rest do not necessarily need to contain a strict inequality.
We, therefore, have $ \sum\limits_{\objIndex=2}^{\numObjs} \slack_{\objIndex}^{\prime} \geq \sum\limits_{\objIndex=2}^{\numObjs} \optSlack_{\objIndex}$, which implies $\  \sum\limits_{\objIndex=2}^{\numObjs} (\slack_{\objIndex}^{\prime} - \optSlack_{\objIndex}) \geq 0$.  
Since $\epsCoefObj_{\objIndex} > 0$ for all $\objIndex \in \{2, \ldots, \numObjs\}$, $ \sum\limits_{\objIndex=2}^{\numObjs} \epsCoefObj_{\objIndex} (\slack_{\objIndex}^{\prime} - \optSlack_{\objIndex}) \geq 0$. We also have ${0 > \objVar_{1}^{\prime} - \objVarOptVal_{1} \geq \sum\limits_{\objIndex=2}^{\numObjs} \epsCoefObj_{\objIndex} (\slack_{\objIndex}^{\prime} - \optSlack_{\objIndex}) \geq 0}$, which is a contradiction again. 
\end{compactenum}

\medskip
\noindent Hence, $\optSol$ is an efficient solution for {\mop} in \eqref{model:mop}. 
\Halmos

%% file: Propositions/prop_underdose_UB_appendix.tex
Any feasible solution to the MOLP in \eqref{model:5obj} yielding $\emph{\cov}_{\tumor}$'s as the tumor coverage values satisfies the following inequalities:
\begin{align}
\sum\limits_{\voxel \in \setVoxel_{\tumor}} \doseUnder_{\tumor \voxel} \ \leq \  \targetDose_{\tumor} \ |\setVoxel_{\tumor}| \ (1-\emph{\cov}_{\tumor}) \qquad  \tumor \in \setTumor \tag{\ref{eq:underdose_extra_cons}}
\end{align}
%
%

%% file: Propositions/proof_underdose_UB.tex
\proof{Proof.}
Suppose that we have a solution to the MOLP in \eqref{model:5obj}, with $\doseUnderSol_{\tumor \voxel}$ being the underdose values, $\setVoxel_{\tumor}^{\doseUnderSol}$ being the set of voxels in tumor $t$ that is underdosed, i.e., that did not fully receive the required dose amount $\targetDose_{\tumor}$, and $\covSol_{\tumor}$ being the resulting coverage of tumor $\tumor$. 
Then,
\begin{align}
    \frac{|\setVoxel_{\tumor}^{\doseUnderSol}|}{|\setVoxel_{\tumor}|} = 1 - \covSol_{\tumor}, \qquad  \tumor \in \setTumor. \label{eq:cov_setvoxel_link}
\end{align}
\noindent Moreover,
\begin{align}
0 \ < \ \doseUnderSol_{\tumor \voxel} \ \leq \ \targetDose_{\tumor} \qquad  \tumor \in \setTumor, \ \voxel \in \setVoxel_{\tumor}^{\doseUnderSol}.  \label{eq:underdose_bounds}
\end{align}
\noindent Summing over $\voxel \in \setVoxel_{\tumor}^{\doseUnderSol}$ on each side of \eqref{eq:underdose_bounds}, we get
\begin{align}
0 \ < \ \sum\limits_{\voxel \in \setVoxel_{\tumor}^{\doseUnderSol}}\doseUnderSol_{\tumor \voxel} \ \leq \ \sum\limits_{\voxel \in \setVoxel_{\tumor}^{\doseUnderSol}} \targetDose_{\tumor} \ = \ \targetDose_{\tumor} \ |\setVoxel_{\tumor}^{\doseUnderSol}|,  \qquad  \tumor \in \setTumor. \label{eq:underdose_bounds_sum}
\end{align}

\noindent Using \eqref{eq:cov_setvoxel_link} and \eqref{eq:underdose_bounds_sum}, we obtain
\begin{align}
\sum\limits_{\voxel \in \setVoxel_{\tumor}^{\doseUnderSol}} \doseUnderSol_{\tumor \voxel} \ \leq \  \targetDose_{\tumor} \ |\setVoxel_{\tumor}| \ (1-\covSol_{\tumor}) \qquad  \tumor \in \setTumor.  \label{eq:underdose_bounds_sum_2}
\end{align}

\noindent Since $\doseUnderSol_{\tumor \voxel} = 0, \ \ \forall \voxel \in \setVoxel_{\tumor} \setminus \setVoxel_{\tumor}^{\doseUnderSol}$, we
can expand the domain of the summation on the left of \eqref{eq:underdose_bounds_sum_2} to $\voxel \in \setVoxel_{\tumor}$, which reduces \eqref{eq:underdose_bounds_sum_2} to \eqref{eq:underdose_extra_cons}.
\Halmos

%% file: Propositions/prop_BOT_LB_appendix.tex
Any feasible solution to the MOLP in \eqref{model:5obj} yielding $\emph{\cov}_{\tumor}$'s as the coverage values satisfies the following inequality:
\begin{align}
\max\limits_{\tumor \in \setTumor} \left\{ \frac{\targetDose_{\tumor} \ \emph{\cov}_{\tumor}}{\max\limits_{\voxel, \isocenter, \sector, \sectorSize} \{ \doseMat_{\tumor \voxel \isocenter \sector \sectorSize} \} \ |\setSector| \ } \right\} \ \leq \ \objVarFive_5  \tag{\ref{eq:BOT_LB}}
\end{align}
%
%

%% file: Propositions/proof_BOT_LB.tex
\proof{Proof.}
Suppose that we have a solution to the MOLP in \eqref{model:5obj}, with $\beamWeightSol_{\isocenter \sector \sectorSize}$'s being the irradiation times, $\doseUnderSol_{\tumor \voxel}$'s being the underdose values, $\BOTSol_{\isocenter}$ being the BOT value of isocenter $\isocenter$, $\objVarFiveSol_{5}$ being the total BOT, and $\covSol_{\tumor}$ the resulting coverage of tumor $\tumor$.
Let us first rewrite the original dose requirement constraint for tumors, using \eqref{cons:dose} and \eqref{cons:underdose} as
\begin{align}
\doseUnderSol_{\tumor \voxel} \ \geq \ \targetDose_{\tumor} \ - \ \sum_{\isocenter \in \setIsocenter} \sum_{\sector \in \setSector} \sum_{\sectorSize \in \setSectorSize} \doseMat_{\tumor \voxel \isocenter \sector \sectorSize} \ \beamWeightSol_{\isocenter \sector \sectorSize},  && \qquad \tumor \in \setTumor, \ \voxel \in \setVoxel_{\tumor}. \label{eq:underdose}
\end{align}

\noindent Summing over $ \voxel \in \setVoxel_{\tumor} $ in \eqref{eq:underdose} and combining with \eqref{eq:underdose_extra_cons}, we obtain
\begin{align}
\targetDose_{\tumor} \ |\setVoxel_{\tumor}| \ (1-\covSol_{\tumor}) \ \geq \ \sum\limits_{\voxel \in \setVoxel_{\tumor}} \doseUnderSol_{\tumor \voxel} \ \geq \ |\setVoxel_{\tumor}| \targetDose_{\tumor} \ - \ \sum\limits_{\voxel \in \setVoxel_{\tumor}} \sum_{\isocenter \in \setIsocenter} \sum_{\sector \in \setSector} \sum_{\sectorSize \in \setSectorSize} \doseMat_{\tumor \voxel \isocenter \sector \sectorSize} \ \beamWeightSol_{\isocenter \sector \sectorSize}  && \qquad \tumor \in \setTumor.  \label{eq:underdose_sum}
\end{align}

\noindent 
We can remove the intermediate underdose term in \eqref{eq:underdose_sum} and manipulate the remaining terms to isolate the BOT objective,
\begin{align*}
\targetDose_{\tumor} \ |\setVoxel_{\tumor}| \ \covSol_{\tumor} \ \leq \ \sum\limits_{\voxel \in \setVoxel_{\tumor}} \sum_{\isocenter \in \setIsocenter} \sum_{\sector \in \setSector} \sum_{\sectorSize \in \setSectorSize} \doseMat_{\tumor \voxel \isocenter \sector \sectorSize} \ \beamWeightSol_{\isocenter \sector \sectorSize} 
\ & \leq \ 
\max\limits_{\voxel, \isocenter, \sector, \sectorSize} \{ \doseMat_{\tumor \voxel \isocenter \sector \sectorSize}  \} \ \sum\limits_{\voxel \in \setVoxel_{\tumor}} \sum_{\isocenter \in \setIsocenter} \sum_{\sector \in \setSector} \sum_{\sectorSize \in \setSectorSize} \beamWeightSol_{\isocenter \sector \sectorSize} & \\[0.2cm]
& = \ \max\limits_{\voxel, \isocenter, \sector, \sectorSize} \{ \doseMat_{\tumor \voxel \isocenter \sector \sectorSize}  \} \ |\setVoxel_{\tumor}| \ \sum_{\isocenter \in \setIsocenter} \sum_{\sector \in \setSector} \sum_{\sectorSize \in \setSectorSize} \beamWeightSol_{\isocenter \sector \sectorSize} & \\[0.2cm]
& \leq \ \max\limits_{\voxel, \isocenter, \sector, \sectorSize} \{ \doseMat_{\tumor \voxel \isocenter \sector \sectorSize}  \} \ |\setVoxel_{\tumor}| \ \sum_{\sector \in \setSector } \sum_{\isocenter \in \setIsocenter} \BOTSol_{\isocenter} &\\[0.2cm]
& = \ \max\limits_{\voxel, \isocenter, \sector, \sectorSize} \{ \doseMat_{\tumor \voxel \isocenter \sector \sectorSize}  \} \ |\setSector| \ |\setVoxel_{\tumor}| \ \sum_{\isocenter \in \setIsocenter} \BOTSol_{\isocenter}, &  \tumor \in \setTumor.  
\end{align*}

\noindent Taking the first and last expressions, we obtain 
\begin{align}
\targetDose_{\tumor} \ \covSol_{\tumor} \ \leq \ \max\limits_{\voxel, \isocenter, \sector, \sectorSize} \{ \doseMat_{\tumor \voxel \isocenter \sector \sectorSize}  \} \ |\setSector| \ \ \sum_{\isocenter \in \setIsocenter} \BOTSol_{\isocenter}, \qquad  \tumor \in \setTumor.  \label{eq:BOT_lb_prelim}
\end{align}

\noindent We rearrange the terms in \eqref{eq:BOT_lb_prelim} to leave only the BOT term on the right-hand side of the inequality,
\begin{align}
\frac{\targetDose_{\tumor} \ \covSol_{\tumor}}{\max\limits_{\voxel, \isocenter, \sector, \sectorSize} \{ \doseMat_{\tumor \voxel \isocenter \sector \sectorSize} \} \ |\setSector| \ } \ \leq \  \ \sum_{\isocenter \in \setIsocenter} \BOTSol_{\isocenter} \ = \ \objVarFiveSol_5, \qquad  \tumor \in \setTumor.  \label{eq:BOT_lb_prelim_2}
\end{align}
\noindent The left-hand side of inequality \eqref{eq:BOT_lb_prelim_2} gives a different lower bound for each tumor $\tumor \in \setTumor$. 
Since the right-hand side is independent of index $\tumor$, the maximum of these bounds over $\tumor \in \setTumor$ serves as a single valid lower bound for the BOT objective $\objVarFive_5$, yielding \eqref{eq:BOT_LB}. 
\Halmos

%% file: Propositions/proof_repeating_sols_1.tex
\proof{Proof.}
Let $\optSlack_{\objIndex}$ be the slack amount for the $\varepsilon_{\objIndex}$ constraint (on $\objVar_{\objIndex}(\allVars) $).
Then, ${\optSlack_{\objIndex} > 0} \ {\forall \objIndex \in \nonbindingSubsetEpsi(\optSoli)}$. 
The feasibility and optimality of a solution to a linear programming problem remains if the right-hand side of a nonbinding (and $\leq$ type) constraint is decreased by at most its slack amount or increased as much as desired. 
Then, with any $\eps_{\objIndex}$ satisfying $ \epsi_{\objIndex} - \optSlack_{\objIndex} \leq \eps_{\objIndex} $ for $\objIndex \in \nonbindingSubsetEpsi(\optSoli) $, the optimality of $\objVarOptVali$ remains.
Since $ \epsi_{\objIndex} - \optSlack_{\objIndex} = \objVarOptVali_{\objIndex} \leq  \epsii  \ \ \forall \objIndex \in \nonbindingSubsetEpsi(\optSoli) $, we conclude that $\objVarOptVali$ remains optimal if $\epsi$ is changed to $\epsii$, and hence $\objVarOptVali \in \optObjSet(\epsii)$.  
\Halmos

%% file: appendix_results.tex
\section{Experimental Results}
\label{app:results}

\input{Tables/table_ML_models_comparison}

\input{Tables/table_prediction_scores}

\input{Tables/table_stats_ML}

\input{Figures/fig_plots_cases5to8}

%% file: Tables/table_ML_models_comparison.tex
\begin{table}[htbp]
  \centering
  \caption{\crevii Comparison of machine learning models using the results from 10 replications of 5-fold cross validation. (For each model, average values are given in the first row, minimum and maximum values are given in the second row. ROC AUC stands for the area under the receiver operating characteristic curve, and abs stands for absolute.)}
    \scalebox{0.66}{
    {\crevii
    \begin{tabular}{l ccccc}
	\toprule
    & \multicolumn{2}{c}{\parbox{4.75cm}{\centering Feasibility classification}} & \multicolumn{3}{c}{\parbox{6cm}{\centering Regression}} \\
	\cmidrule(lr){2-3} \cmidrule(lr){4-6}
	\parbox{3cm}{Model} 	&  \parbox{3.75cm}{\centering Accuracy} & \parbox{3.75cm}{\centering ROC AUC} &  \parbox{3.75cm}{\centering Cov abs error} 	&  \parbox{3.75cm}{\centering PCI abs error}   	&  \parbox{3.75cm}{\centering BOT abs error} \\
	\midrule	   
    	Random forest    			& \bf 0.99511               &	\bf 0.99807 	        & \bf 0.00030 	            & \bf 0.00027 	            & \bf 0.14772 \\
    	                            & \bf [0.97580, 1]          &	\bf [0.98079, 1] 	    & \bf [0.00000, 0.00135] 	& \bf [0.00000, 0.00147] 	& \bf [0.00000, 0.65693] \\[0.3cm]
    	Decision tree    			&	0.99497             &	0.98370 	        & 0.00030 	            & 0.00028 	            & 0.16255 \\
    	                            & [0.97483, 1]          &	[0.90792, 1] 	    & [0.00000, 0.00151] 	& [0.00000, 0.00147] 	& [0.00000, 0.78643] \\[0.3cm]
    	Nearest neighbor   			&	0.98482             &	0.97032 	        & 0.00155 	            & 0.00128 	            & 0.57949 \\
    	                            & [0.95454, 1]          &	[0.87918, 1] 	    & [0.00018, 0.00411] 	& [0.00016, 0.00348] 	& [0.05127, 1.66809] \\[0.3cm]
    	Ridge 					    &	0.96790             &	0.98639 	        & 0.00620 	            & 0.00725 	            & 1.26637 \\
    	                            & [0.92559, 0.99420]    &	[0.94245, 1] 	    & [0.00454, 0.00822] 	& [0.00554, 0.00928] 	& [0.77902, 1.95676] \\[0.3cm]
    	Support vector machine   	&	0.97634             &	0.97051 	        & 0.02505 	            & 0.03960 	            & 3.08289 \\
    	                            & [0.93426, 0.99710]    &	[0.89769, 1] 	    & [0.01997, 0.03099] 	& [0.03596, 0.04394] 	& [1.84248, 4.52256] \\[0.3cm]
    	Logistic/linear regression	&	0.95354             &	0.94504 	        & 0.00620 	            & 0.00725 	            & 1.26638 \\
    	                            & [0.89942, 0.99710]    &	[0.80405, 0.99992] 	& [0.00454, 0.00822] 	& [0.00554, 0.00928] 	& [0.77903, 1.95677] \\
	\bottomrule
    \end{tabular}%
    }
    }
  \label{tab:ml_models_comparison}%
\end{table}%

%% file: Tables/table_prediction_scores.tex
\begin{table}[htbp]
  \centering
  \caption{Prediction performance of random forest models.}
  \scalebox{0.67}{
    \begin{tabular}{c  c S[table-format=2.1] S[table-format=1.3] S[table-format=1.3]  c S[table-format=3.1] S[table-format=1.3] S[table-format=1.3] S[table-format=2.2]}
    \toprule
    & \multicolumn{4}{c}{Phase I} & \multicolumn{5}{c}{Phase II} \\
    \cmidrule(lr){2-5} \cmidrule(lr){6-10}\\[-0.6cm]
   \parbox{1.5cm}{\centering Case} 
   & \parbox{2.75cm}{\centering Train set \\sizes} 
   & \parbox{2cm}{\centering \% correct feas class} 
   & \parbox{2cm}{\centering Cov avg abs error} 
   & \parbox{2cm}{\centering PCI avg abs error} 
   & \parbox{2.75cm}{\centering Train set \\sizes} 
   & \parbox{2cm}{\centering \% correct feas class} 
   & \parbox{2cm}{\centering Cov avg abs error} 
   & \parbox{2cm}{\centering PCI avg abs error} 
   & \parbox{2cm}{\centering BOT avg abs error} \\[0.3cm]
    \midrule
    1  &  (185, 102)    & 98.2  &  0.021  &  0.021  &  (190, 105)    &  100.0  &  0.000  &  0.003  &   0.97 \\
    2  &  (1284, 1205)  & 92.4  &  0.027  &  0.041  &  (1288, 1208)  &  100.0  &  0.000  &  0.068  &   0.84 \\
    3  &  (376, 318)    & 96.5  &  0.028  &  0.038  &  (380, 322)    &   80.0  &  0.000  &  0.025  &  19.06 \\
    4  &  (206, 133)    & 94.7  &  0.028  &  0.024  &  (215, 140)    &  100.0  &  0.000  &  0.036  &   6.47 \\
    5  &  (324, 264)    & 96.3  &  0.010  &  0.017  &  (340, 277)    &  100.0  &  0.000  &  0.013  &   2.70 \\
    6  &  (355, 304)    & 96.4  &  0.015  &  0.021  &  (364, 313)    &  100.0  &  0.000  &  0.015  &  11.40 \\
    7  &  (394, 336)    & 97.5  &  0.009  &  0.010  &  (413, 355)    &  100.0  &  0.000  &  0.011  &   9.23 \\
    8  &  (2048, 1979)  & 93.6  &  0.019  &  0.014  &  (2224, 2127)  &   88.0  &  0.000  &  0.014  &   3.47 \\
    \midrule
    Average  &  & 95.7  & 0.020 & 0.023 &  & 98.5  & 0.000 & 0.023 & 6.77 \\
    \bottomrule
    \end{tabular}%
    }
  \label{tab:prediction_scores}%
\end{table}%

%% file: Tables/table_stats_ML.tex
\begin{table}[htbp]
  \centering
  \caption{Statistics from the ML-guided version of our two-phase method.}
  \scalebox{0.65}{
    \begin{tabular}{cc c S[table-format=2.1] S[table-format=2.0] S[table-format=2.1] >{\columncolor[gray]{0.95}}S[table-format=2.1] >{\columncolor[gray]{0.95}}S[table-format=2.1] >{\columncolor[gray]{0.95}}S[table-format=2.1] S[table-format=2.1] S[table-format=3.1] S[table-format=5.1] S[table-format=3.1]}
    \toprule \\[-0.7cm]
      & \parbox{1.5cm}{\centering Case} 
      & \parbox{1.75cm}{\centering \# \\ $\eps$-vectors} 
      & \parbox{2.3cm}{\centering \% {\twophasereg} $\eps$-vectors} 
      & \parbox{1.3cm}{\centering \# \\sol} 
      & \parbox{1.3cm}{\centering \% \\sol} 
      & \parbox{1.5cm}{\centering \% infeas} 
      & \parbox{1.75cm}{\centering \% omit in infeas} 
      & \parbox{1.75cm}{\centering \% omit for infeas} 
      & \parbox{1.5cm}{\centering \% omit for uniq} 
      & \parbox{2cm}{\centering Unit sol time (sec)} 
      & \parbox{2cm}{\centering Total sol time (sec)} 
      & \parbox{2.5cm}{\centering \% {\twophasereg} time} \\[0.3cm]
  \midrule
 \multirow{8}[0]{*}{Phase I}
          & 1     & {100+5}     & 10.5  & 9     & 8.6     & 81.0  & 85.9  & 69.5  & 10.5  & \crevi 27.2  & \crevi 245.0   & \crevi 149.8 \\
          & 2     & {100+4}     & 1.0   & 14    & 13.5    & 76.9  & 72.5  & 55.8  & 9.6   & \crevi 28.3  & \crevi 395.9  & \crevi 27.9 \\
          & 3     & {100+4}     & 10.4  & 22    & 21.2    & 55.8  & 84.5  & 47.1  & 23.1  & \crevi 4.4  & \crevi 96.0  & \crevi 52.5 \\
          & 4     & {100+9}     & 10.9  & 21    & 19.3    & 68.8  & 85.3  & 58.7  & 11.9  & \crevi 154.6 & \crevi 3247.4 & \crevi 52.7 \\
          & 5     & {100+10}    & 11.0  & 17    & 15.5    & 57.3  & 87.3  & 50.0  & 27.3  & \crevi 1.0   & \crevi 16.5    & \crevi 69.6 \\
          & 6     & {100+9}     & 10.9  & 27    & 24.8    & 46.8  & 84.3  & 39.4  & 28.4  & \crevi 6.6  & \crevi 178.2   & \crevi 55.4 \\
          & 7     & {100+19}    & 11.9  & 30    & 25.2    & 48.7  & 77.6  & 37.8  & 26.1  & \crevi 13.8  & \crevi 414.3  & \crevi 77.2 \\
          & 8     & {100+116}   & 2.2   & 43    & 19.9    & 44.9  & 77.3  & 34.7  & 35.2  & \crevi 3.6  & \crevi 153.7   & \crevi 63.7 \\
  \midrule                                                                                
  Average &       &  122.0      & 8.6   & 22.9  & 18.5    & 60.0  & 81.8  & 49.1 & 21.5  & \crevi 29.9 & \crevi 593.4  & \crevi 52.4 \\
  \midrule
  \multirow{8}[0]{*}{Phase II}
          & 1     & 15    & 60.0    & 3     & 12.0     & 0.0    & 0.0    & 0.0   & 12.0   & \crevi 7.9   & \crevi 23.8   & \crevi 73.7 \\
          & 2     & 25    & 100.0   & 5     & 20.0     & 0.0    & 0.0    & 0.0   & 20.0   & \crevi 3.2    & \crevi 15.8   & \crevi 38.0 \\
          & 3     & 12    & 48.0    & 5     & 60.0     & 0.0    & 0.0    & 0.0   & 20.0   & \crevi 2.6    & \crevi 12.9   & \crevi 32.5 \\
          & 4     & 10    & 40.0    & 9     & 80.0     & 0.0    & 0.0    & 0.0   & 36.0   & \crevi 39.2  & \crevi 352.8  & \crevi 69.9 \\
          & 5     & 20    & 80.0    & 15    & 75.0     & 0.0    & 0.0    & 0.0   & 60.0   & \crevi 0.4    & \crevi 5.3   & \crevi 75.7 \\
          & 6     & 8     & 32.0    & 4     & 100.0    & 0.0    & 0.0    & 0.0   & 16.0   & \crevi 4.3    & \crevi 17.4   & \crevi 21.2 \\
          & 7     & 5     & 20.0    & 1     & 60.0     & 0.0    & 0.0    & 0.0   & 4.0    & \crevi 37.9   & \crevi 37.9   & \crevi 24.3 \\
          & 8     & 21    & 84.0    & 8     & 48 .0    & 12.0   & 66.7   & 8.0   & 32.0   & \crevi 3.0    & \crevi 23.7   & \crevi 45.3 \\
  \midrule                                                                                 
  Average &       & 14.5  & 63.0    & 7.4   & 56.9     & 1.5    & 66.7   & 1.0   & 25.0   & \crevi 12.3   & \crevi 61.2  & \crevi 53.5 \\
  \bottomrule
    \end{tabular}%
    }
  \label{tab:stats_ml}%
\end{table}%

%% file: Figures/fig_plots_cases5to8.tex
\begin{figure}[!h]
    \centering
	\begin{subfigure}{0.4\textwidth}
		\centering
	    \includegraphics[scale=0.475]{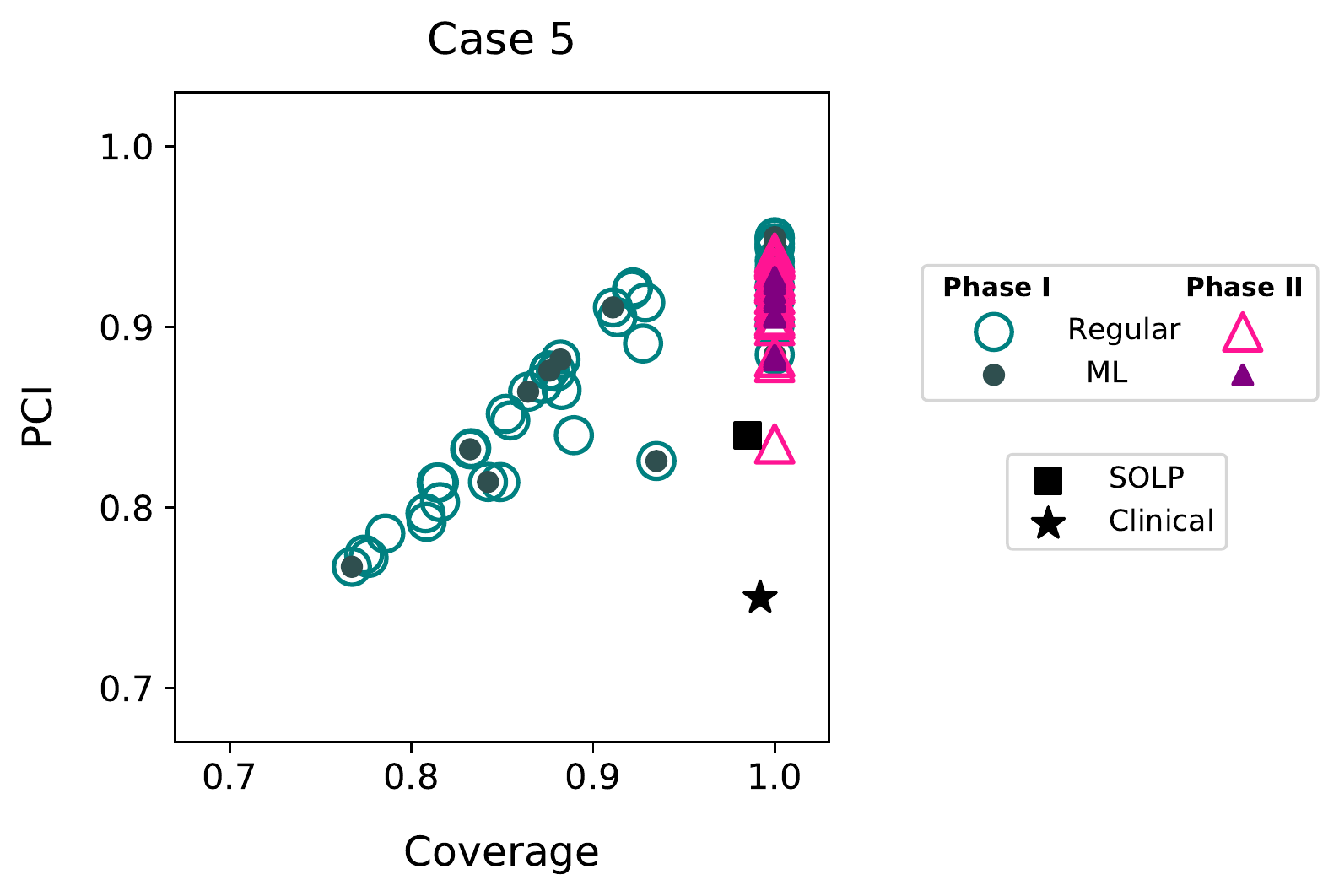}
	\end{subfigure}
	~\medskip
	\begin{subfigure}{0.4\textwidth}
		\centering
        \includegraphics[scale=0.475]{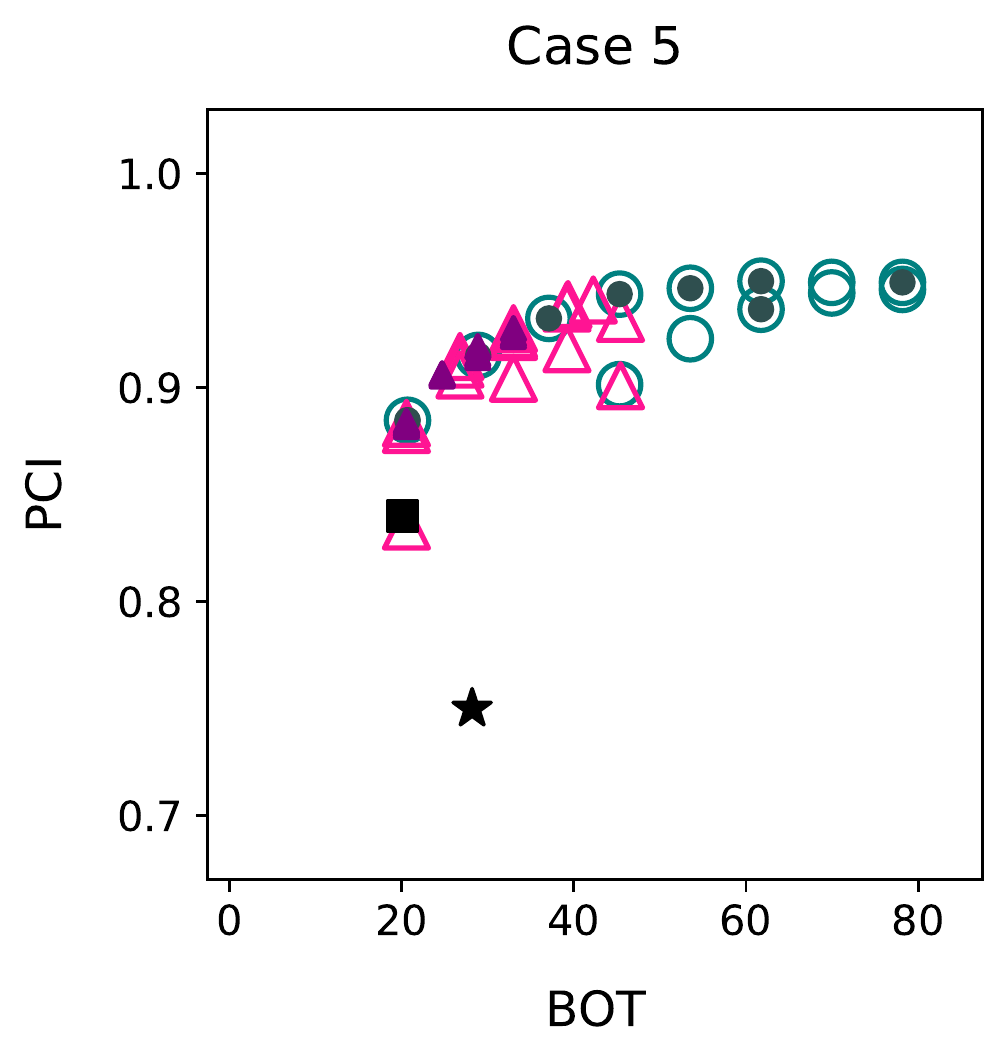}
	\end{subfigure}
	~\medskip
	\begin{subfigure}{0.4\textwidth}
		\centering
	    \includegraphics[scale=0.475]{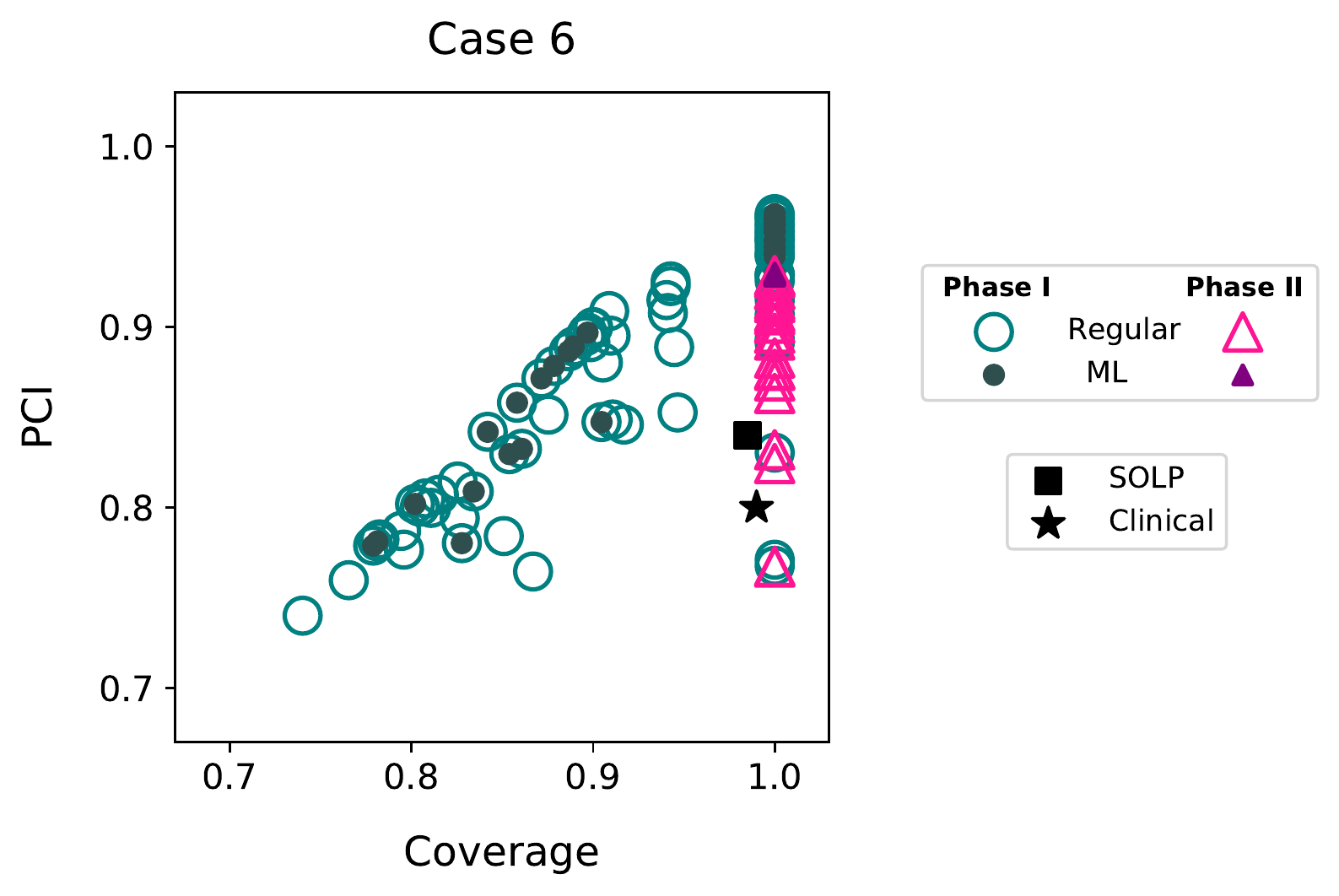}
	\end{subfigure}
	~
	\begin{subfigure}{0.4\textwidth}
		\centering
        \includegraphics[scale=0.475]{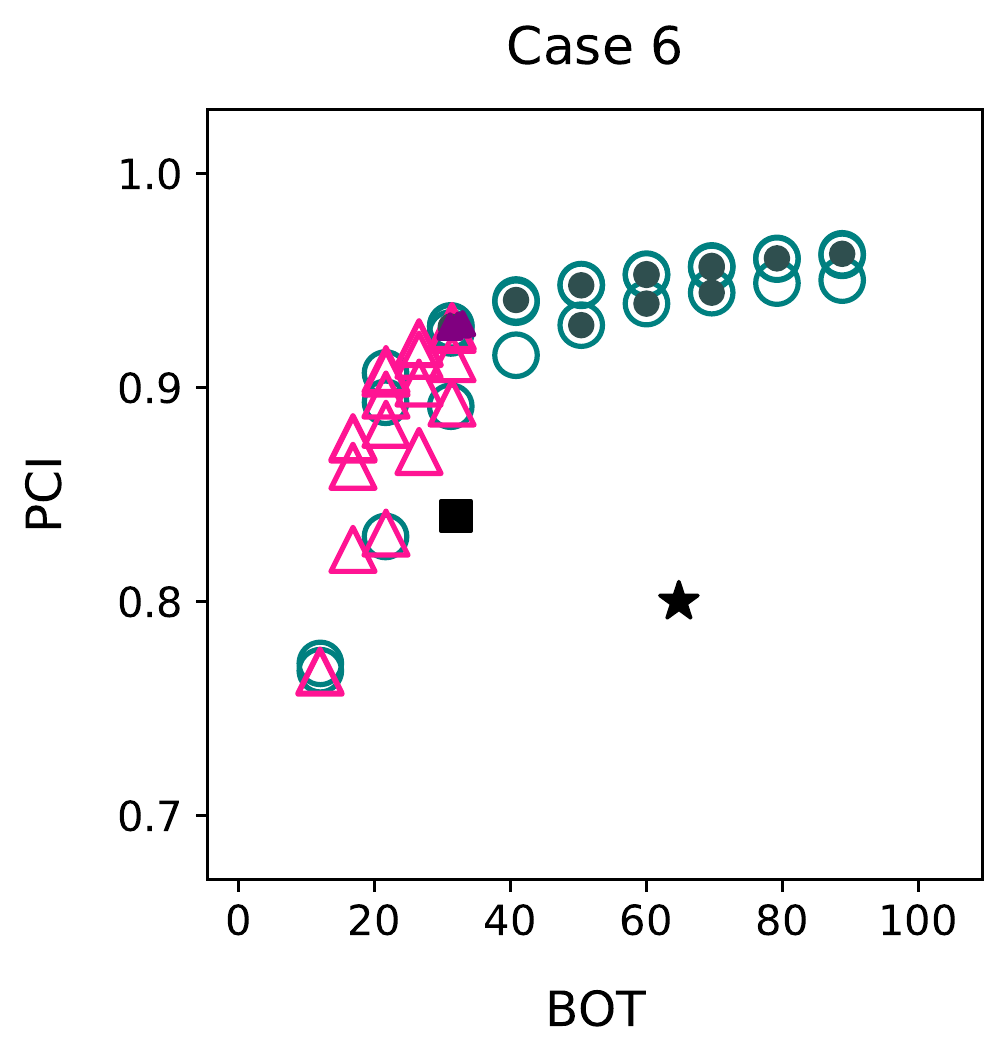}
	\end{subfigure}
	~\medskip
	\begin{subfigure}{0.4\textwidth}
		\centering
	    \includegraphics[scale=0.475]{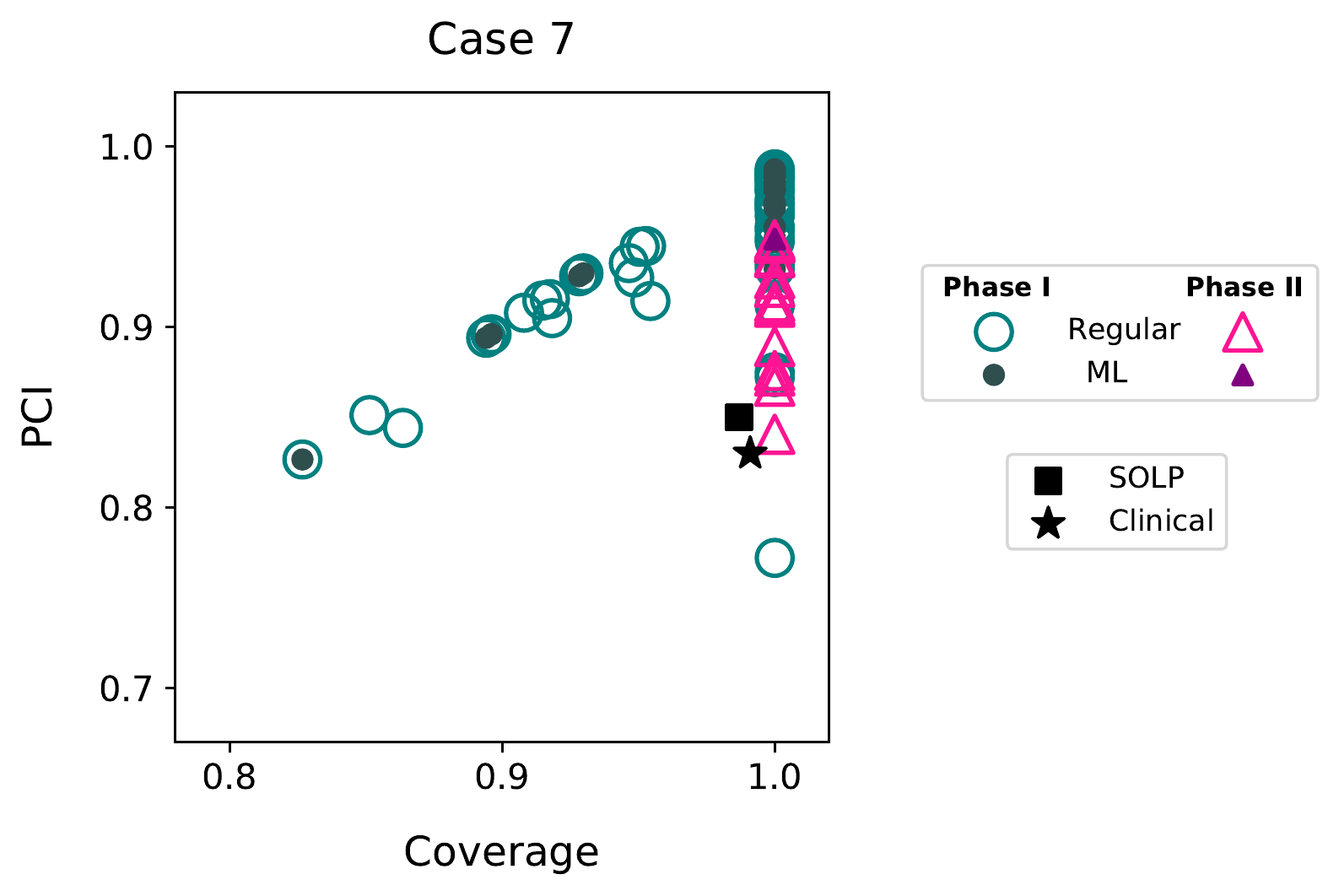}
	\end{subfigure}
	~
	\begin{subfigure}{0.4\textwidth}
		\centering
        \includegraphics[scale=0.475]{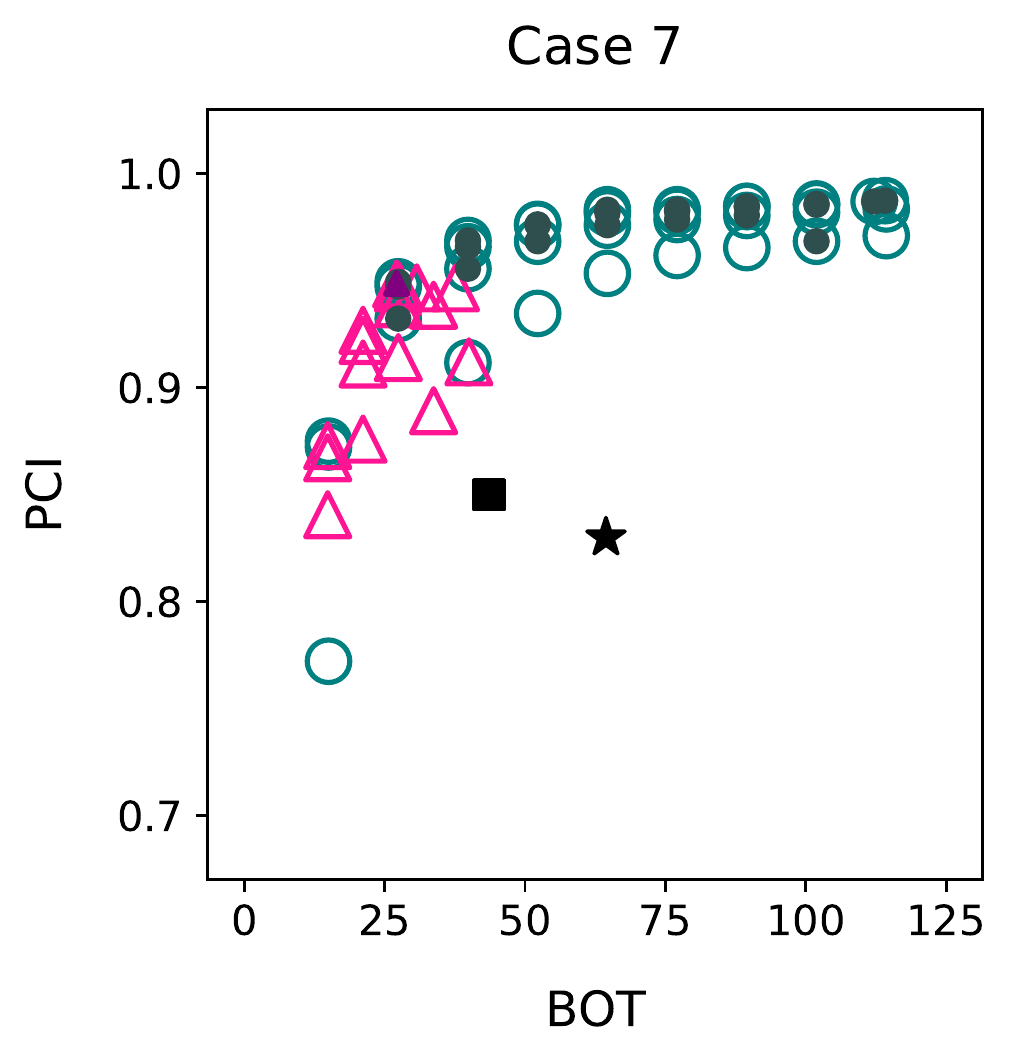}
	\end{subfigure}
	~\medskip
	\begin{subfigure}{0.4\textwidth}
		\centering
	    \includegraphics[scale=0.475]{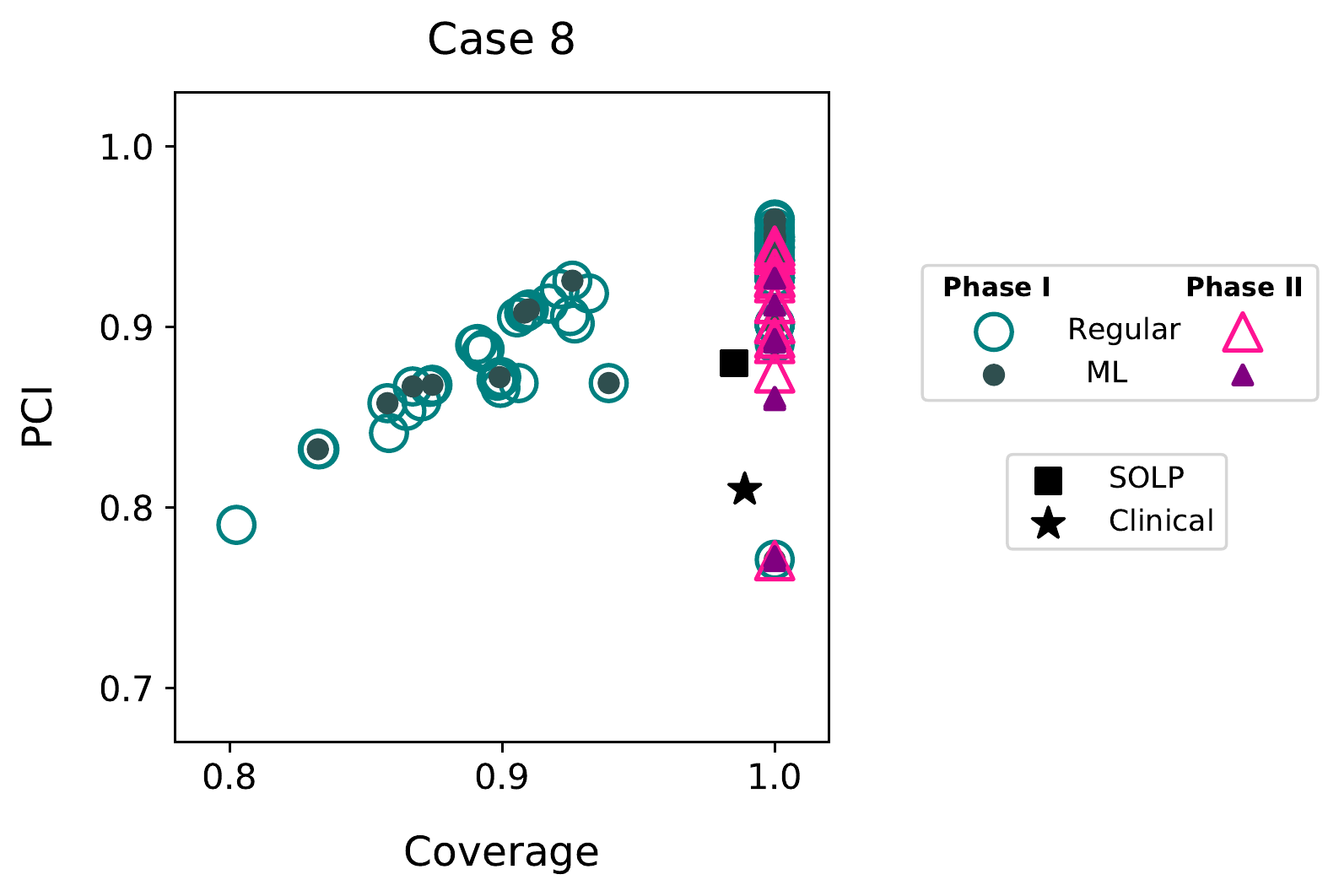}
	\end{subfigure}
	~
	\begin{subfigure}{0.4\textwidth}
		\centering
        \includegraphics[scale=0.475]{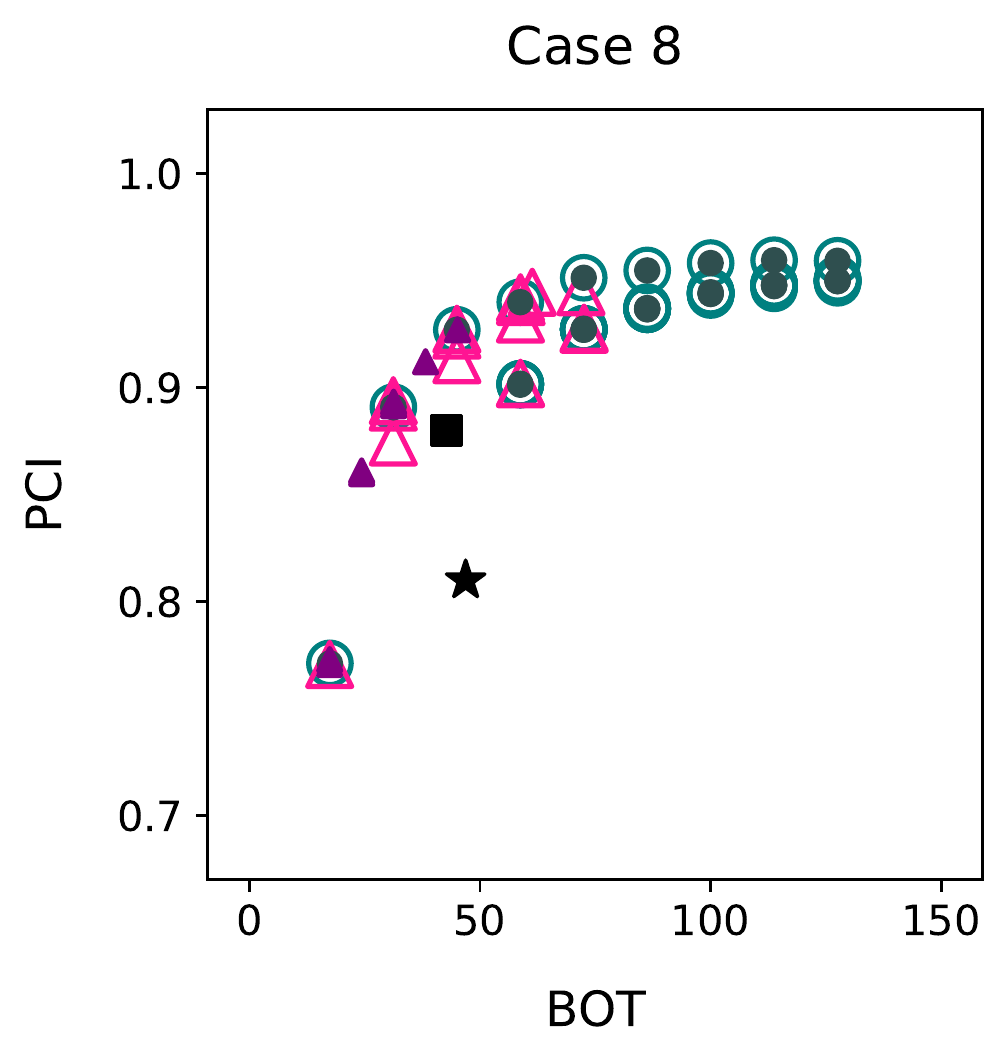}
	\end{subfigure}
	\caption{Coverage vs. PCI values of all generated solutions (left), and BOT vs. PCI values of solutions with at least 99.7\% coverage (right), for Cases 5 to 8. 
	}
	\label{fig:cases5to8}
\end{figure}

%% file: main_OO.bbl
\begin{thebibliography}{65}
\providecommand{\natexlab}[1]{#1}
\providecommand{\url}[1]{\texttt{#1}}
\providecommand{\urlprefix}{URL }

\bibitem[{Babier et~al.(2018)Babier, Boutilier, Sharpe, McNiven,
  \protect\BIBand{} Chan}]{babier2018inverse}
Babier A, Boutilier JJ, Sharpe MB, McNiven AL, Chan TC (2018) Inverse
  optimization of objective function weights for treatment planning using
  clinical dose-volume histograms. \emph{Physics in Medicine \& Biology}
  63(10):105004.

\bibitem[{Breedveld et~al.(2019)Breedveld, Craft, Van~Haveren,
  \protect\BIBand{} Heijmen}]{breedveld2019multi}
Breedveld S, Craft D, Van~Haveren R, Heijmen B (2019) Multi-criteria
  optimization and decision-making in radiotherapy. \emph{European Journal of
  Operational Research} 277(1):1--19.

\bibitem[{Breedveld et~al.(2012)Breedveld, Storchi, Voet, \protect\BIBand{}
  Heijmen}]{breedveld2012icycle}
Breedveld S, Storchi PR, Voet PW, Heijmen BJ (2012) {iCycle}: Integrated,
  multicriterial beam angle, and profile optimization for generation of
  coplanar and noncoplanar {IMRT} plans. \emph{Medical Physics} 39(2):951--963.

\bibitem[{Cabrera-Guerrero et~al.(2018{\natexlab{a}})Cabrera-Guerrero, Ehrgott,
  Mason, \protect\BIBand{} Raith}]{cabrera2018matheuristic}
Cabrera-Guerrero G, Ehrgott M, Mason AJ, Raith A (2018{\natexlab{a}}) A
  matheuristic approach to solve the multiobjective beam angle optimization
  problem in intensity-modulated radiation therapy. \emph{International
  Transactions in Operational Research} 25(1):243--268.

\bibitem[{Cabrera-Guerrero et~al.(2018{\natexlab{b}})Cabrera-Guerrero, Mason,
  Raith, \protect\BIBand{} Ehrgott}]{cabrera2018pareto}
Cabrera-Guerrero G, Mason AJ, Raith A, Ehrgott M (2018{\natexlab{b}}) Pareto
  local search algorithms for the multi-objective beam angle optimisation
  problem. \emph{Journal of Heuristics} 24(2):205--238.

\bibitem[{Cevik et~al.(2019)Cevik, Aleman, Lee, Berdyshev, Nordstr{\"o}m, Riad,
  Sahgal, \protect\BIBand{} Ruschin}]{cevik2019simultaneous}
Cevik M, Aleman D, Lee Y, Berdyshev A, Nordstr{\"o}m H, Riad S, Sahgal A,
  Ruschin M (2019) Simultaneous optimization of isocenter locations and sector
  duration in radiosurgery. \emph{Physics in Medicine \& Biology} 64(2):025010.

\bibitem[{Cevik et~al.(2018)Cevik, Ghomi, Aleman, Lee, Berdyshev, Nordstrom,
  Riad, Sahgal, \protect\BIBand{} Ruschin}]{cevik2018modeling}
Cevik M, Ghomi PS, Aleman D, Lee Y, Berdyshev A, Nordstrom H, Riad S, Sahgal A,
  Ruschin M (2018) Modeling and comparison of alternative approaches for sector
  duration optimization in a dedicated radiosurgery system. \emph{Physics in
  Medicine \& Biology} 63(15):155009.

\bibitem[{Chankong \protect\BIBand{} Haimes(1983)}]{chankong1983multiobjective}
Chankong V, Haimes YY (1983) \emph{Multiobjective decision making theory and
  methodology} (Elsevier).

\bibitem[{Cheek et~al.(2005)Cheek, Holder, Fuss, \protect\BIBand{}
  Salter}]{cheek2005relationship}
Cheek D, Holder A, Fuss M, Salter B (2005) The relationship between the number
  of shots and the quality of gamma knife radiosurgeries. \emph{Optimization
  and Engineering} 6(4):449--462.

\bibitem[{Chen \protect\BIBand{} Girvigian(2005)}]{chen2005stereotactic}
Chen JC, Girvigian MR (2005) Stereotactic radiosurgery: instrumentation and
  theoretical aspects—part 1. \emph{The Permanente Journal} 9(4):23.

\bibitem[{Cohon(1978)}]{cohon1978multiobjective}
Cohon JL (1978) \emph{Multiobjective programming and planning}, volume 140 of
  \emph{Mathematics in Science and Engineering} (Academic Press).

\bibitem[{Doudareva et~al.(2015)Doudareva, Ghobadi, Aleman, Ruschin,
  \protect\BIBand{} Jaffray}]{doudareva2015skeletonization}
Doudareva E, Ghobadi K, Aleman DM, Ruschin M, Jaffray DA (2015) Skeletonization
  for isocentre selection in {Gamma Knife\textsuperscript{\textregistered}
  Perfexion\textsuperscript{\texttrademark}}. \emph{{TOP}} 23(2):369--385.

\bibitem[{Ehrgott(2005)}]{ehrgott2005multicriteria}
Ehrgott M (2005) \emph{Multicriteria optimization}, volume 491 (Springer
  Science \& Business Media).

\bibitem[{Ehrgott(2006)}]{ehrgott2006discussion}
Ehrgott M (2006) A discussion of scalarization techniques for multiple
  objective integer programming. \emph{Annals of Operations Research}
  147(1):343--360.

\bibitem[{Ehrgott et~al.(2010)Ehrgott, G{\"u}ler, Hamacher, \protect\BIBand{}
  Shao}]{ehrgott2010mathematical}
Ehrgott M, G{\"u}ler {\c{C}}, Hamacher HW, Shao L (2010) Mathematical
  optimization in intensity modulated radiation therapy. \emph{Annals of
  Operations Research} 175(1):309--365.

\bibitem[{Ehrgott \protect\BIBand{} Ruzika(2008)}]{ehrgott2008improved}
Ehrgott M, Ruzika S (2008) Improved $\varepsilon$-constraint method for
  multiobjective programming. \emph{Journal of Optimization Theory and
  Applications} 138(3):375.

\bibitem[{{Elekta}(2010)}]{leksellGammaPlan2010}
{Elekta} (2010) {Inverse planning in {Leksell GammaPlan} (white paper)}.

\bibitem[{Eus{\'e}bio et~al.(2014)Eus{\'e}bio, Figueira, \protect\BIBand{}
  Ehrgott}]{eusebio2014finding}
Eus{\'e}bio A, Figueira JR, Ehrgott M (2014) On finding representative
  non-dominated points for bi-objective integer network flow problems.
  \emph{Computers \& operations research} 48:1--10.

\bibitem[{Ferris et~al.(2002)Ferris, Lim, \protect\BIBand{}
  Shepard}]{ferris2002optimization}
Ferris MC, Lim J, Shepard DM (2002) An optimization approach for radiosurgery
  treatment planning. \emph{SIAM Journal on Optimization} 13(3):921--937.

\bibitem[{Ferris et~al.(2003)Ferris, Lim, \protect\BIBand{}
  Shepard}]{ferris2003radiosurgery}
Ferris MC, Lim J, Shepard DM (2003) Radiosurgery treatment planning via
  nonlinear programming. \emph{Annals of Operations Research} 119(1):247--260,
  \urlprefix\url{https://doi.org/10.1023/A:1022951027498}.

\bibitem[{Ferris \protect\BIBand{} Shepard(2000)}]{ferris2000optimization}
Ferris MC, Shepard DM (2000) Optimization of gamma knife radiosurgery.
  \emph{Discrete Mathematical Problems with Medical Applications} 55:27--44.

\bibitem[{Ghaffari et~al.(2017)Ghaffari, Aleman, Jaffray, \protect\BIBand{}
  Ruschin}]{ghaffari2017incorporation}
Ghaffari HR, Aleman DM, Jaffray DA, Ruschin M (2017) Incorporation of delivery
  times in stereotactic radiosurgery treatment optimization. \emph{Journal of
  Global Optimization} 69(1):103--115.

\bibitem[{Ghobadi et~al.(2012)Ghobadi, Ghaffari, Aleman, Jaffray,
  \protect\BIBand{} Ruschin}]{ghobadi2012automated}
Ghobadi K, Ghaffari HR, Aleman DM, Jaffray DA, Ruschin M (2012) Automated
  treatment planning for a dedicated multi-source intracranial radiosurgery
  treatment unit using projected gradient and grassfire algorithms.
  \emph{Medical Physics} 39(6Part1):3134--3141.

\bibitem[{Ghobadi et~al.(2013)Ghobadi, Ghaffari, Aleman, Jaffray,
  \protect\BIBand{} Ruschin}]{ghobadi2013automated}
Ghobadi K, Ghaffari HR, Aleman DM, Jaffray DA, Ruschin M (2013) Automated
  treatment planning for a dedicated multi-source intra-cranial radiosurgery
  treatment unit accounting for overlapping structures and dose homogeneity.
  \emph{Medical Physics} 40(9):091715.

\bibitem[{Giller(2011)}]{giller2011feasibility}
Giller CA (2011) Feasibility of identification of {Gamma Knife} planning
  strategies by identification of {Pareto} optimal {Gamma Knife} plans.
  \emph{Technology in Cancer Research \& Treatment} 10(6):561--574.

\bibitem[{Haimes et~al.(1971)Haimes, Lasdon, \protect\BIBand{}
  Wismer}]{haimes1971bicriterion}
Haimes YY, Lasdon LS, Wismer DA (1971) On a bicriterion formation of the
  problems of integrated system identification and system optimization.
  \emph{IEEE Transactions on Systems, Man and Cybernetics} 1(3):296--297.

\bibitem[{Harris \protect\BIBand{} Das(2020)}]{maxBOTref}
Harris L, Das JM (2020) \emph{Stereotactic Radiosurgery} (StatPearls
  Publishing, Treasure Island (FL)).

\bibitem[{{IBM ILOG}(2019)}]{cplex1210}
{IBM ILOG} (2019) {CPLEX} {Optimization Studio} 12.10.0 {U}ser manual .

\bibitem[{Isermann \protect\BIBand{} Steuer(1988)}]{isermann1988computational}
Isermann H, Steuer RE (1988) Computational experience concerning payoff tables
  and minimum criterion values over the efficient set. \emph{European Journal
  of Operational Research} 33(1):91--97.

\bibitem[{James et~al.(2013)James, Witten, Hastie, \protect\BIBand{}
  Tibshirani}]{statlearningbook2013}
James G, Witten D, Hastie T, Tibshirani R (2013) \emph{An introduction to
  statistical learning}, volume 112 (Springer).

\bibitem[{Jindal \protect\BIBand{} Sangwan(2017)}]{jindal2017multi}
Jindal A, Sangwan KS (2017) Multi-objective fuzzy mathematical modelling of
  closed-loop supply chain considering economical and environmental factors.
  \emph{Annals of Operations Research} 257(1):95--120.

\bibitem[{Kirlik \protect\BIBand{} Say{\i}n(2014)}]{kirlik2014new}
Kirlik G, Say{\i}n S (2014) A new algorithm for generating all nondominated
  solutions of multiobjective discrete optimization problems. \emph{European
  Journal of Operational Research} 232(3):479--488.

\bibitem[{Laumanns et~al.(2006)Laumanns, Thiele, \protect\BIBand{}
  Zitzler}]{laumanns2006efficient}
Laumanns M, Thiele L, Zitzler E (2006) An efficient, adaptive parameter
  variation scheme for metaheuristics based on the epsilon-constraint method.
  \emph{European Journal of Operational Research} 169(3):932--942.

\bibitem[{Lee et~al.(2000)Lee, Fox, \protect\BIBand{}
  Crocker}]{lee2000optimization}
Lee EK, Fox T, Crocker I (2000) Optimization of radiosurgery treatment planning
  via mixed integer programming. \emph{Medical physics} 27(5):995--1004.

\bibitem[{Lee et~al.(2006)Lee, Barber, \protect\BIBand{}
  Walton}]{lee2006automated}
Lee KJ, Barber DC, Walton L (2006) Automated {Gamma Knife} radiosurgery
  treatment planning with image registration, data-mining, and {Nelder-Mead}
  simplex optimization. \emph{Medical Physics} 33(7):2532--2540.

\bibitem[{Mavrotas(2009)}]{mavrotas2009effective}
Mavrotas G (2009) Effective implementation of the {$\varepsilon$}-constraint
  method in multi-objective mathematical programming problems. \emph{Applied
  Mathematics and Computation} 213(2):455--465.

\bibitem[{Mavrotas \protect\BIBand{} Florios(2013)}]{mavrotas2013improved}
Mavrotas G, Florios K (2013) An improved version of the augmented
  $\varepsilon$-constraint method {(AUGMECON2)} for finding the exact pareto
  set in multi-objective integer programming problems. \emph{Applied
  Mathematics and Computation} 219(18):9652--9669.

\bibitem[{Nikas et~al.(2020)Nikas, Fountoulakis, Forouli, \protect\BIBand{}
  Doukas}]{nikas2020robust}
Nikas A, Fountoulakis A, Forouli A, Doukas H (2020) A robust augmented
  $\varepsilon$-constraint method (augmecon-r) for finding exact solutions of
  multi-objective linear programming problems. \emph{Operational Research}
  1--42.

\bibitem[{Novotn{\`y}(2012)}]{novotny2012leksell}
Novotn{\`y} J (2012) {Leksell Gamma Knife} -- past, present and future.
  \emph{L{\'e}ka{\v{r}} a Technika-Clinician and Technology} 42(3):5--13.

\bibitem[{Oskoorouchi et~al.(2011)Oskoorouchi, Ghaffari, Terlaky,
  \protect\BIBand{} Aleman}]{oskoorouchi2011interior}
Oskoorouchi MR, Ghaffari HR, Terlaky T, Aleman DM (2011) An interior point
  constraint generation algorithm for semi-infinite optimization with
  health-care application. \emph{Operations Research} 59(5):1184--1197.

\bibitem[{{\"O}zlen \protect\BIBand{} Azizo{\u{g}}lu(2009)}]{ozlen2009multi}
{\"O}zlen M, Azizo{\u{g}}lu M (2009) Multi-objective integer programming: A
  general approach for generating all non-dominated solutions. \emph{European
  Journal of Operational Research} 199(1):25--35.

\bibitem[{Paddick(2000)}]{paddick2000simple}
Paddick I (2000) A simple scoring ratio to index the conformity of
  radiosurgical treatment plans. \emph{Journal of neurosurgery}
  93(supplement\_3):219--222.

\bibitem[{Padilla \protect\BIBand{} Palta(2019)}]{srsbook2019ch7}
Padilla L, Palta JR (2019) Overview of technologies for {SRS} and {SBRT}
  delivery. Heron DE, Huq MS, Herman JM, eds., \emph{Stereotactic Radiosurgery
  and Stereotactic Body Radiation Therapy (SBRT)}, chapter~7, 73--101
  (Springer).

\bibitem[{Pedregosa et~al.(2011)Pedregosa, Varoquaux, Gramfort, Michel,
  Thirion, Grisel, Blondel, Prettenhofer, Weiss, Dubourg, Vanderplas, Passos,
  Cournapeau, Brucher, Perrot, \protect\BIBand{} Duchesnay}]{scikit-learn}
Pedregosa F, Varoquaux G, Gramfort A, Michel V, Thirion B, Grisel O, Blondel M,
  Prettenhofer P, Weiss R, Dubourg V, Vanderplas J, Passos A, Cournapeau D,
  Brucher M, Perrot M, Duchesnay E (2011) Scikit-learn: Machine learning in
  {P}ython. \emph{Journal of Machine Learning Research} 12:2825--2830.

\bibitem[{Ripsman et~al.(2015)Ripsman, Aleman, \protect\BIBand{}
  Ghobadi}]{ripsman2015interactive}
Ripsman DA, Aleman DM, Ghobadi K (2015) Interactive visual guidance for
  automated stereotactic radiosurgery treatment planning. \emph{Expert Systems
  with Applications} 42(21):8337--8348.

\bibitem[{Romeijn et~al.(2004)Romeijn, Dempsey, \protect\BIBand{}
  Li}]{romeijn2004unifying}
Romeijn HE, Dempsey JF, Li JG (2004) A unifying framework for multi-criteria
  fluence map optimization models. \emph{Physics in Medicine \& Biology}
  49(10):1991.

\bibitem[{Sahebjamnia et~al.(2015)Sahebjamnia, Torabi, \protect\BIBand{}
  Mansouri}]{sahebjamnia2015integrated}
Sahebjamnia N, Torabi SA, Mansouri SA (2015) Integrated business continuity and
  disaster recovery planning: Towards organizational resilience. \emph{European
  Journal of Operational Research} 242(1):261--273.

\bibitem[{Say{\i}n(2000)}]{sayin2000measuring}
Say{\i}n S (2000) Measuring the quality of discrete representations of
  efficient sets in multiple objective mathematical programming.
  \emph{Mathematical Programming} 87(3):543--560.

\bibitem[{Schlaefer \protect\BIBand{} Schweikard(2008)}]{schlaefer2008stepwise}
Schlaefer A, Schweikard A (2008) Stepwise multi-criteria optimization for
  robotic radiosurgery. \emph{Medical Physics} 35(5):2094--2103.

\bibitem[{Shao \protect\BIBand{} Ehrgott(2016)}]{shao2016discrete}
Shao L, Ehrgott M (2016) Discrete representation of non-dominated sets in
  multi-objective linear programming. \emph{European Journal of Operational
  Research} 255(3):687--698.

\bibitem[{Shepard et~al.(2003)Shepard, Chin, DiBiase, Naqvi, Lim,
  \protect\BIBand{} Ferris}]{shepard2003clinical}
Shepard DM, Chin LS, DiBiase SJ, Naqvi SA, Lim J, Ferris MC (2003) Clinical
  implementation of an automated planning system for {Gamma Knife}
  radiosurgery. \emph{International Journal of Radiation Oncology, Biology,
  Physics} 56(5):1488--1494.

\bibitem[{Shepard et~al.(1999)Shepard, Ferris, Olivera, \protect\BIBand{}
  Mackie}]{shepard1999optimizing}
Shepard DM, Ferris MC, Olivera GH, Mackie TR (1999) Optimizing the delivery of
  radiation therapy to cancer patients. \emph{Siam Review} 41(4):721--744.

\bibitem[{Shepard et~al.(2000)Shepard, Ferris, Ove, \protect\BIBand{}
  Ma}]{shepard2000inverse}
Shepard DM, Ferris MC, Ove R, Ma L (2000) Inverse treatment planning for gamma
  knife radiosurgery. \emph{Medical Physics} 27(12):2748--2756.

\bibitem[{Sj{\"o}lund et~al.(2019)Sj{\"o}lund, Riad, Hennix, \protect\BIBand{}
  Nordstr{\"o}m}]{sjolund2019linear}
Sj{\"o}lund J, Riad S, Hennix M, Nordstr{\"o}m H (2019) A linear programming
  approach to inverse planning in {Gamma Knife} radiosurgery. \emph{Medical
  Physics} 46(4):1533--1544.

\bibitem[{Svensson(2014)}]{svensson2014multiobjective}
Svensson J (2014) \emph{Multiobjective optimization in radiosurgery: How to
  approximate and navigate on the Pareto surface}. Master's thesis, KTH Royal
  Institute of Technology.

\bibitem[{Tamby \protect\BIBand{} Vanderpooten(2021)}]{tamby2021enumeration}
Tamby S, Vanderpooten D (2021) Enumeration of the nondominated set of
  multiobjective discrete optimization problems. \emph{INFORMS Journal on
  Computing} 33(1):72--85.

\bibitem[{Tian et~al.(2020)Tian, Yang, Giles, Wang, Gao, Butker, Liu,
  \protect\BIBand{} Kahn}]{tian2020preliminary}
Tian Z, Yang X, Giles M, Wang T, Gao H, Butker E, Liu T, Kahn S (2020) A
  preliminary study on a multiresolution-level inverse planning approach for
  {Gamma Knife} radiosurgery. \emph{Medical Physics} 47(4):1523--1532.

\bibitem[{Wagner et~al.(2000)Wagner, Yi, Meeks, Bova, Brechner, Chen, Buatti,
  Friedman, Foote, \protect\BIBand{} Bouchet}]{wagner2000geometrically}
Wagner TH, Yi T, Meeks SL, Bova FJ, Brechner BL, Chen Y, Buatti JM, Friedman
  WA, Foote KD, Bouchet LG (2000) A geometrically based method for automated
  radiosurgery planning. \emph{International Journal of Radiation Oncology,
  Biology, Physics} 48(5):1599--1611.

\bibitem[{Weistroffer(1985)}]{weistroffer1985careful}
Weistroffer H (1985) Careful usage of pessimistic values is needed in multiple
  objectives optimization. \emph{Operations Research Letters} 4(1):23--25.

\bibitem[{Yu(1997)}]{yu1997multiobjective}
Yu Y (1997) Multiobjective decision theory for computational optimization in
  radiation therapy. \emph{Medical Physics} 24(9):1445--1454.

\bibitem[{Yu et~al.(2000)Yu, Zhang, Cheng, Schell, \protect\BIBand{}
  Okunieff}]{yu2000multi}
Yu Y, Zhang J, Cheng G, Schell M, Okunieff P (2000) Multi-objective
  optimization in radiotherapy: applications to stereotactic radiosurgery and
  prostate brachytherapy. \emph{Artificial Intelligence in Medicine}
  19(1):39--51.

\bibitem[{Zeverino et~al.(2017)Zeverino, Jaccard, Patin, Ryckx, Marguet,
  Tuleasca, Schiappacasse, Bourhis, Levivier, Bochud, \protect\BIBand{}
  Moeckli}]{zeverino2017}
Zeverino M, Jaccard M, Patin D, Ryckx N, Marguet M, Tuleasca C, Schiappacasse
  L, Bourhis J, Levivier M, Bochud FO, Moeckli R (2017) Commissioning of the
  {Leksell Gamma Knife\textsuperscript{\textregistered}
  Icon\textsuperscript{\texttrademark}}. \emph{Medical Physics} 44(2):355--363.

\bibitem[{Zhang et~al.(2001)Zhang, Dean, Metzger, \protect\BIBand{}
  Sibata}]{zhang2001optimization}
Zhang P, Dean D, Metzger A, Sibata C (2001) Optimization of gamma knife
  treatment planning via guided evolutionary simulated annealing. \emph{Medical
  Physics} 28(8):1746--1752.

\bibitem[{Zhang et~al.(2003)Zhang, Wu, Dean, Xing, Xue, Maciunas,
  \protect\BIBand{} Sibata}]{zhang2003plug}
Zhang P, Wu J, Dean D, Xing L, Xue J, Maciunas R, Sibata C (2003) Plug pattern
  optimization for gamma knife radiosurgery treatment planning.
  \emph{International Journal of Radiation Oncology* Biology* Physics}
  55(2):420--427.

\bibitem[{Zhang \protect\BIBand{} Reimann(2014)}]{zhang2014simple}
Zhang W, Reimann M (2014) A simple augmented $\epsilon$-constraint method for
  multi-objective mathematical integer programming problems. \emph{European
  Journal of Operational Research} 234(1):15--24.

\end{thebibliography}
